\documentclass[envcountsame,envcountsect]{svmult}

\usepackage{makeidx}         
\usepackage{multicol}        
\usepackage[bottom]{footmisc}

\usepackage{amssymb, amsmath}
\usepackage{ownmacs}

\newcommand{\der}{\mathop{{\rm der}}\nolimits}
\renewcommand{\div}{\mathop{{\rm div}}\nolimits}
\newcommand{\Sp}{\mathop{{\rm Sp}}\nolimits}
\newcommand{\HSp}{\mathop{{\rm HSp}}\nolimits}
\newcommand{\Bp}{{\overline \Psi}}
\newcommand{\fint}{=\!\!\!\!\!\!\!\int}

\begin{document}

\title*{Lie groups of bundle automorphisms and their extensions} 
\author{Karl-Hermann Neeb}
\institute{Department of Mathematics, 
Darmstadt University of Technology, 
Schlossgartenstrasse 7, 64289 Darmstadt, Germany 
\texttt{neeb@mathematik.tu-darmstadt.de}}%
\maketitle

\begin{abstract} We describe natural abelian extensions of the Lie algebra 
$\aut(P)$ of infinitesimal automorphisms of a principal 
bundle over a compact manifold $M$ 
and discuss their integrability to corresponding 
Lie group extensions. Already the case of a trivial 
bundle $P = M \times K$ is quite interesting. 
In this case, we show that essentially all central extensions 
of the gauge algebra $C^\infty(M,\fk)$ can be obtained from 
three fundamental types of cocycles with values in one of the spaces 
$\fz := C^\infty(M,V)$, $\Omega^1(M,V)$ and 
$\Omega^1(M,V)/\dd C^\infty(M,V)$. These cocycles extend 
to $\aut(P)$, and, under the assumption that $TM$ is trivial, we 
also describe the space $H^2({\cal V}(M),\fz)$ classifying 
the twists of these extensions. We then show that all fundamental 
types have natural generalizations to non-trivial bundles and 
explain under which conditions they extend to $\aut(P)$ and 
integrate to global Lie group extensions. 
\end{abstract}

\section*{Introduction} \label{sec:intro}

Two of the most important classes of infinite-dimensional Lie groups 
are groups of smooth maps, such as the Lie group 
$C^\infty(M,K)$ of smooth maps of a compact smooth manifold $M$ with 
values in a Lie group $K$, and groups of diffeomorphisms, 
such as the group $\Diff(M)$ of diffeomorphisms of a compact smooth 
manifold $M$. 

A case of particular importance arises for the circle $M = \bS^1$, 
where $LK := C^\infty(\bS^1,K)$ is called the {\it loop group} of $K$. If 
$K$ is a compact simple Lie group, then $LK$ has a universal central 
extension $\hat LK$ by the circle group $\T$. Furthermore, 
the group $\T_r \subeq \Diff(\bS^1)$ of rigid rotations of $\bS^1$ 
acts smoothly 
by automorphisms on $LK$, this action lifts to the central 
extension $\hat LK$, and we 
thus obtain the {\it affine Kac--Moody Lie groups}
$\hat LK \rtimes \T_r$.\begin{footnote} {Strictly speaking, these 
are the ``unitary forms'' of the affine Kac--Moody groups. Starting 
with a complex simple Lie group $K$ instead, we obtain 
complex versions.}
\end{footnote} 
The {\it twisted affine Kac--Moody Lie groups} can be realized 
in this picture as the fixed point groups for an automorphism 
$\sigma$ acting trivially on $\T_r$, and inducing on 
$LK$ an automorphism of the form 
$\sigma(f)(t) = \phi(f(\zeta t))$, where 
$\phi \in \Aut(K)$ is an automorphism of finite order $m$ 
and $\zeta\in \T_r$ satisfies $\zeta^m = 1$. 
On each affine Kac--Moody group, the group 
$\Diff_+(\bS^1)$ of orientation preserving diffeomorphisms of $\bS^1$  
acts by automorphism, so that 
$\hat LK \rtimes \T_r$ embeds into 
$\hat LK \rtimes \Diff_+(\bS^1)$. Actually, the latter group permits 
an interesting twist, where the subgroup $\T \times \Diff_+(\bS^1)$ 
is replaced by the {\it Virasoro group}, a non-trivial central extension 
of $\Diff_+(\bS^1)$ by $\T$. 

The purpose of this note is to discuss several 
infinite-dimensional Lie groups that are 
constructed in a similar fashion from higher dimensional compact manifolds $M$. Since 
the one-dimensional manifold $M = \bS^1$ is a rather simple object, 
the general theory leaves much more room for different Lie group constructions, 
extensions and twistings thereof. Here we shall focus on 
recent progress in several branches of this area, in particular relating 
the Lie algebra picture to global objects. We shall also explain how 
this relates to other structures, such as multiloop
algebras, which are currently 
under active investigation from the algebraic point of view 
(\cite{ABP06}, \cite{ABFP07}). 
We also present open problems and describe several areas where the 
present knowledge is far from satisfactory.

The analog of the loop algebra $C^\infty(\bS^1,\fk)$ of a Lie algebra $\fk$ 
is the Lie algebra of smooth maps $\g = C^\infty(M,\fk)$ 
on the compact manifold 
$M$. Here the compactness assumption is convenient if we want to deal 
with Lie groups, because it implies that $\g$ is 
the Lie algebra of any group $G = C^\infty(M,K)$, where $K$ is a Lie 
group with Lie algebra $\fk$. 
This group can also be identified with the gauge 
group of the trivial $K$-principal bundle $P = M \times K$, so that 
gauge groups of principal bundles are natural generalizations of 
mapping groups. Twisted loop groups are gauge groups of certain 
non-trivial principal bundles over $\bS^1$ (cf.\ Section~\ref{sec:gaugeext}). 

For a $K$-principal bundle $q \: P \to M$ over the compact manifold $M$,  
we write $\Aut(P) := \Diff(P)^K$ for the group of all diffeomorphisms 
of $P$ commuting with the $K$-action, i.e., the group of {\it bundle 
automorphisms}. 
The {\it gauge group} $\Gau(P)$ is the normal subgroup of all 
bundle isomorphisms inducing the identity on $M$. 
Writing $\Diff(M)_P$ for the set of all diffeomorphisms of $M$ 
that can be lifted to  
bundle automorphisms, which is an open subgroup of the Lie group $\Diff(M)$, 
we obtain a short exact sequence of Lie groups 
\begin{equation}
  \label{eq:seq1}
\1 \to \Gau(P) \to \Aut(P) \to \Diff(M)_P \to \1 
\end{equation}
(cf.\ \cite{ACM89}, \cite{Wo06}). 
On the Lie algebra level, we have a corresponding short exact sequence of Lie algebras 
of vector fields 
\begin{equation}
  \label{eq:seq2}
\0 \to \gau(P) \to \aut(P) := {\cal V}(P)^K \to {\cal V}(M)\to \1.
\end{equation}

It is an important problem to understand the central extensions of 
gauge groups by an abelian Lie group $Z$ on which $\Diff(M)_P$ acts 
naturally and the extent to which they can be enlarged to abelian extensions of 
the full group $\Aut(P)$, or at least its identity component. 
Whenever such 
an enlargement exists, one has to understand the set of all enlargements, 
the twistings, 
which leads to the problem to classify all abelian extensions of $\Diff(M)_P$ 
by $Z$. Below we shall 
discuss various partial results concerning these questions. 
We also describe several 
tools to address them and to explain what still remains to be done. 

The contents of the paper is as follows: 
After discussing old and new results on Lie group structures on 
gauge groups, mapping groups and diffeomorphisms groups in Section~\ref{sec:1}, 
we turn in Section~\ref{sec:2} to central extensions of mapping 
groups $G := C^\infty(M,K)$. Here we first exhibit three fundamental 
classes of cocycles over which all other non-trivial cocycles can be factored, 
up to cocycles vanishing on the commutator algebra. 
The fundamental cocycles are defined by an invariant symmetric bilinear form 
$\kappa\: \fk \times \fk \to V$ and an alternating map 
$\eta \: \fk \times \fk \to V$ as follows: 
$$ \omega_\kappa(\xi_1, \xi_2) := [\kappa(\xi_1, \dd \xi_2)] \in 
\oline\Omega^1(M,V) := \Omega^1(M,V)/\dd C^\infty(M,V), $$
$$ \omega_\eta(\xi_1, \xi_2) := \eta(\xi_1, \xi_2) \in  C^\infty(M,V)
\quad\mbox{ for } \quad \dd_\fk \eta = 0, $$
and if $\dd_\fk\eta(x,y,z) = \kappa([x,y],z)$, we also have 
$$ \omega_{\kappa,\eta}(\xi_1, \xi_2) := 
\kappa(\xi_1, \dd \xi_2) - \kappa(\xi_2, \dd \xi_1)- \dd(\eta(\xi_1, \xi_2)) 
\in \Omega^1(M,V). $$
We then discuss the action of $G_0 \rtimes \Diff(M)$ on the corresponding 
central Lie algebra extension and explain under which conditions it 
integrates to an extension of  Lie groups. 
Since we are also interested in the corresponding abelian extensions of 
the full automorphism group $\Aut(M \times K) \cong G \rtimes \Diff(M)$ and 
its Lie algebra, we turn in Section~\ref{sec:3} to abelian extensions of 
${\cal V}(M)$, resp., $\Diff(M)$ by the three types of target spaces $\fz$ 
from above. If $TM$ is trivial and $V = \R$, then a full description of 
$H^2({\cal V}(M),\fz)$ in all cases has been obtained recently in \cite{BiNe07} 
(Theorem~\ref{thm:h2vect}). 

In Section~\ref{sec:gaugeext} we then turn to gauge and automorphism groups of 
general principal bundles, where the first step consists in 
the construction of analogs of the 
three fundamental types of cocycles. A crucial 
difficulty arises from the fact that one has to consider principal bundles 
with non-connected structure groups to realize natural classes of Lie algebras 
such as twisted affine algebras as gauge algebras. We therefore have to 
consider target spaces $V$ on which the quotient group $\pi_0(K) = K/K_0$ 
acts non-trivially and consider $K$-invariant forms $\kappa \: \fk \times \fk \to V$. 
A typical example is the target space $V(\fk)$ of the universal invariant 
symmetric bilinear form of a semisimple Lie algebra $\fk$ on which 
$K = \Aut(\fk)$ acts non-trivially.  
Replacing the exterior differential by a covariant derivative, we thus obtain 
cocycles 
$$ \omega_\kappa(\xi_1, \xi_2) := [\kappa(\xi_1, \dd \xi_2)] \in 
\oline\Omega^1(M,\bV) := \Omega^1(M,\bV)/\dd \Gamma \bV, $$
where $\bV \to M$ is the flat vector bundle associated to $P$ via the 
$K$-module structure on $V$. 
For the cocycles $\omega_\kappa$ there are natural criteria 
for integrability whenever $\pi_0(K)$ is finite (\cite{NeWo07}), but 
it is not so clear how to extend $\omega_\kappa$ to $\aut(P)$.  
However, the situation simplifies considerably 
for the cocycle $\dd \circ \omega_\kappa$ with values in $\Omega^2(M,\bV)$ 
which extends naturally to a cocycle on $\aut(P)$ that integrates to a smooth 
group cocycle on $\Aut(P)$. 

The analogs of 
$C^\infty(M,V)$-valued cocycles are determined by a central extension 
$\hat\fk$ of $\fk$ by $V$ to which the adjoint action of $K$ lifts. 
Here the central Lie algebra extension $\hat\gau(P)$ is a space of sections 
of an associated Lie algebra bundle with fiber $\hat\fk$. 
We then show that for any central extension $\hat K$ of $K$ by $V$ 
there  exists a $\hat K$-bundle $\hat P$ with $\hat P/V \cong P$, 
so that $\Aut(\hat P)$ is an abelian extension of $\Aut(P)$ by 
$C^\infty(M,V)$ containing the central extension $\Gau(\hat P)$ of $\Gau(P)$ 
with Lie algebra $\hat\gau(P)$. 
There are also analogs of the $\Omega^1(M,V)$-valued cocycles which exist 
whenever the $3$-cocycle $\kappa([x,y],z)$ is a coboundary. 
For non-trivial bundles the fact that the extension $\aut(P)$ of 
${\cal V}(M)$ does in general not split makes the analysis of the 
abelian extensions of $\aut(P)$ considerably more difficult than in the 
trivial case, where we have a semidirect product Lie algebra. 
As a consequence, the results on group actions on the Lie algebra extensions 
and on integrability are much less complete than for 
trivial bundles. 

To connect our geometric setting with the algebraic setup, we describe in 
Section~\ref{sec:loopalg} the connection 
between multiloop algebras and gauge algebras of flat bundles 
over tori which can be trivialized by finite coverings. 
The paper concludes with some remarks on bundles with infinite-dimensional 
structure groups and appendices on 
integrability of abelian Lie algebra extensions, 
extensions of semidirect products (such as 
$C^\infty(M,\fk) \rtimes {\cal V}(M)$), and the triviality of the 
action of a connected Lie group on the corresponding continuous 
Lie algebra cohomology. 

Throughout we only consider Lie algebras and groups of 
smooth maps endowed with the compact open $C^\infty$-topology. 
To simplify matters, we focus on compact manifolds to work 
in a convenient Lie theoretic setup for the corresponding groups. 
For non-compact manifolds, the natural 
setting is provided by spaces of compactly supported smooth maps 
with the direct limit 
$LF$-topology (cf.\ \cite{KM97}, \cite{ACM89}, \cite{Mi80}, 
\cite{Gl02}). However, there are some classes of 
non-compact manifolds on which we still have Lie group structure on 
the full group of smooth maps (Theorem~\ref{thm:real-lie}). 
For groups of sections of group bundles one can also develop 
a theory for $C^k$- and Sobolev sections, 
but this has the disadvantage that 
smooth vector fields do not act as derivations (\cite{Sch04}).

\subsection*{Notation and basic concepts} 

A {\it Lie group} $G$ is a group equipped with a 
smooth manifold structure modeled on a locally convex space 
for which the group multiplication and the 
inversion are smooth maps. We write $\1 \in G$ for the identity element and 
$\lambda_g(x) = gx$, resp., $\rho_g(x) = xg$ for the left, resp.,
right multiplication on $G$. Then each $x \in T_\1(G)$ corresponds to
a unique left invariant vector field $x_l$ with 
$x_l(g) := d\lambda_g(\1).x, g \in G.$
The space of left invariant vector fields is closed under the Lie
bracket of vector fields, hence inherits a Lie algebra structure. In
this sense we obtain on $T_\1(G)$ a continuous Lie bracket which
is uniquely determined by $[x,y]_l = [x_l, y_l]$ for $x,y \in T_\1(G)$. 
We write $\L(G) = \g$ for the so obtained locally convex Lie algebra and 
note that for morphisms $\phi \: G \to H$ of Lie groups we obtain 
with $\L(\phi) := T_\1(\phi)$ a functor from the category of Lie groups 
to the category of locally convex Lie algebras. 
We write $q_G \: \tilde G_0 \to G_0$ for the universal covering map of the identity 
component $G_0$ of $G$ and identify 
the discrete central subgroup $\ker q_G$ of $\tilde G_0$ with $\pi_1(G) \cong \pi_1(G_0)$. 

In the following we always write $I = [0,1]$ for the unit interval in $\R$. 
A Lie group $G$ is called {\it regular} if for each 
$\xi \in C^\infty(I,\g)$, the initial value problem 
$$ \gamma(0) = \1, \quad \gamma'(t) = \gamma(t).\xi(t) = T(\lambda_{\gamma(t)})\xi(t) $$
has a solution $\gamma_\xi \in C^\infty(I,G)$, and the
evolution map 
$$ \evol_G \: C^\infty(I,\g) \to G, \quad \xi \mapsto \gamma_\xi(1) $$
is smooth (cf.\ \cite{Mil84}). 
For a locally convex space $E$, the regularity of the Lie group $(E,+)$ is equivalent 
to the Mackey completeness of $E$, i.e., to the existence of integrals 
of smooth curves $\gamma \: I \to E$. We also recall that for each regular 
Lie group $G$, its Lie algebra $\g$ is Mackey complete and 
that all Banach--Lie groups are regular (\cite{GN07}). 
For a smooth map $f \: M \to G$ we define 
the left, resp., right logarithmic derivative in 
$\Omega^1(M,\g)$ by 
$\delta^l(f)v := f(m)^{-1} \cdot T_m(f)v$  and 
$\delta^r(f)v := T_m(f)v\cdot f(m)^{-1}$, where $\cdot$ refers to the 
two-sided action of $G$ on its tangent bundle $TG$. 

A smooth map $\exp_G \: \L(G) \to G$ is said to be an {\it exponential 
function} if for each $x \in \L(G)$, the curve $\gamma_x(t) := \exp_G(tx)$ is a 
homomorphism $\R \to G$ with $\gamma_x'(0) = x$. Presently, all known 
Lie groups modelled on complete locally convex spaces are regular, hence 
possess an exponential function. For Banach--Lie groups, its existence follows from the 
theory of ordinary differential equations in Banach spaces. 
A Lie group $G$ is called {\it locally exponential}, if 
it has an exponential function mapping an open $0$-neighborhood in 
$\L(G)$ diffeomorphically onto an open neighborhood of $\1$ in $G$. 
For more details, we refer to Milnor's lecture notes \cite{Mil84}, 
the survey \cite{Ne06a} or the monograph \cite{GN07}. 

If $q \: E \to B$ is a smooth fiber bundle, then we write 
$\Gamma E := \break \{ s \in C^\infty(B,E) \: q \circ s = \id_B\}$ for its space 
of smooth sections. 

If $\g$ is a topological Lie algebra and $V$ a topological $\g$-module, 
we write $(C^\bullet(\g,V),\dd_\g)$ for the corresponding Lie algebra 
complex (\cite{ChE48}). 

\section{Lie group structures on mapping groups and automorphism groups of bundles} 
\label{sec:1}

Before we turn to extensions, we collect in this short section some 
results on Lie group structures on gauge groups and automorphism group of 
bundles and introduce notation and conventions used below. 
Some of the material is classical, but there are also some new 
interesting results. 

\subsection{Automorphism groups of bundles} 

\begin{theorem} Let $M$ be a compact manifold and 
$K$ a Lie group with Lie algebra $\fk$. Then the following assertions hold: 
\begin{enumerate}
\item[\rm(1)] $\Diff(M)$ carries the structure of a Fr\'echet--Lie group 
whose Lie algebra  is the space ${\cal V}(M)$ of smooth vector fields on $M$. 
\item[\rm(2)] $G := C^\infty(M,K)$ carries a natural Lie group structure 
with Lie algebra $\g := C^\infty(M,\fk)$, endowed with the pointwise bracket. 
If $(\phi,U)$ is a $\fk$-chart of $K$, hen 
$(\phi_G,U_G)$ with $U_G := \{f \in G \: f(M) \subeq U\}$ and 
$\phi_G(f) = \phi_K \circ f$ is a $\g$-chart of $G$.  
\end{enumerate}
\end{theorem}

For (1) we refer to \cite{Les67}, \cite{Omo70}, and \cite{Ha82}, and for (2) 
we refer to \cite{Mil84} and \cite{Mi80} (see also \cite{GN07} for all that). 

\begin{theorem} If $q \: P \to M$ is a smooth principal bundle over a compact 
manifold $M$ with locally exponential structure group $K$, then its 
gauge group $\Gau(P)$ and its automorphism group $\Aut(P)$ 
carry Lie group structures turning 
(\ref{eq:seq1}) into a Lie group extension, where 
the group $\Diff(M)_P$ is open in $\Diff(M)$. 
In particular, the conjugation action of $\Aut(P)$ on $\Gau(P)$ is smooth.
\end{theorem}

For the case where $K$ is finite-dimensional this can be found in 
\cite{KYMO85}, \cite{ACM89} and \cite{KM97}. 
For general locally exponential Lie groups $K$, this has recently 
been proved by Ch.~Wockel \cite{Wo06}.  

\begin{remark} \label{rem:2.3} (a) If 
$P = M \times K$ is a trivial bundle, the group extension $\Aut(P)$ 
splits and we have 
$$ \Aut(P) \cong \Gau(P) \rtimes \Diff(M) \quad \mbox{ where } \quad 
\Gau(P) \cong C^\infty(M,K). $$
Here $\Diff(M)$ acts on $P$ by $\phi.(m,k) := (\phi(m),k)$ and 
$C^\infty(M,K)$ acts by $f.(m,k) = (m,f(m)k)$. 

(b) For any principal $K$-bundle $P$, each gauge transformation 
$\phi \in \Gau(P)$ determines a unique smooth function 
$f \: P \to K$ by $\phi(p) = p.f(p)$ for $p \in P$. 
Then $f$ satisfies $f(p.k) = k^{-1}f(p)k$, 
and any such function corresponds to a gauge transformation $\phi_f$, 
so that 
$$ \Gau(P) \cong \{ f \in C^\infty(P,K) \: 
(\forall p \in P)(\forall k \in K)\ f(p.k) = k^{-1}f(p)k\}. $$ 
Accordingly, we have on the level of vector fields 
$$ \gau(P) \cong \{ \xi \in C^\infty(P,\fk) \: 
(\forall p \in P)(\forall k \in K)\ \xi(p.k) = \Ad(k)^{-1}\xi(p)\}. $$

(c) For any connected manifold $M$, the universal covering map \break 
$q_M \: \tilde M \to M$ defines on $\tilde M$ the structure of 
a $\pi_1(M)$-principal bundle, 
where $\pi_1(M)$ denotes the group of deck transformations of this 
bundle, considered as a discrete group. Then 
$$ \Gau(\tilde M) = Z(\pi_1(M))\quad \hbox{ and } \quad 
\Aut(\tilde M) = C_{\Diff(\tilde M)}(\pi_1(M)). $$
If $M$ is compact, then $\Aut(\tilde M)$ carries a natural Lie group 
structure turning it into a covering of the open subgroup 
$\Diff(M)_{[\tilde M]}$ of $\Diff(M)$ with kernel $Z(\pi_1(M))$. 
The normalizer $\tilde\Diff(M)$ of $\pi_1(M)$ 
in $\Diff(\tilde M)$ carries a Lie group structure that leads to a 
short exact sequence 
$$ \1 \to \pi_1(M) \to \tilde\Diff(M) \to \Diff(M) \to \1. $$
This is due to the fact that each diffeomorphism $\phi$ of $M$ lifts 
to some diffeomorphism of $\tilde M$ normalizing $\pi_1(M)$, but in general 
it does not centralize $\pi_1(M)$. 

(d) Let $\rho \: \pi_1(M) \to K$ be a homomorphism and 
$$P_\rho := \tilde M \times_\rho K :=  (\tilde M \times K)/\pi_1(M)$$ 
be the corresponding flat $K$-bundle associated to the $\pi_1(M)$-bundle 
$\tilde M$ via $\rho$. Then $P_\rho$ is the 
set of $\pi_1(M)$-orbits in $\tilde M \times K$ 
for the right action $(\tilde m,k).\gamma := 
(\tilde m.\gamma, \rho(\gamma)^{-1}k)$. We write 
$[(\tilde m,k)]$ for the orbit of the pair $(\tilde m,k)$. 
Then \break $s \: \tilde M \to P, 
s(\tilde m) := [(\tilde m,\1)]$ 
is a smooth function with $q \circ s = q_M$, satisfying 
$$s(\tilde m.\gamma) = s(\tilde m).\rho(\gamma) \quad \mbox{ for } \quad 
\tilde m \in \tilde M, \gamma \in \pi_1(M).$$ 

In view of $P = s(\tilde M).K$, any bundle automorphism 
$\hat\phi$ is determined by its values on the image of $s$. 
If $\phi \in \Diff(M)$ is the corresponding diffeomorphism of $M$ 
and $\tilde \phi$ is a lift of $\phi$ to a diffeomorphism of $\tilde M$, 
then $\hat \phi$ can be written 
\begin{equation}
  \label{eq:hatphi}
\hat\phi(s(\tilde m)) = s(\tilde\phi(\tilde m)).f(\tilde m) 
\end{equation}
for some smooth function $f \: \tilde M \to K$, satisfying 
\begin{equation}
  \label{eq:hatphif}
f(\tilde m) \rho(\gamma) = \rho(\tilde\phi \gamma \tilde\phi^{-1}) 
f(\tilde m.\gamma) \quad \mbox{ for } \quad \tilde m \in \tilde M, \gamma \in \pi_1(M).
\end{equation}
Conversely, if (\ref{eq:hatphif}) it satisfied, then (\ref{eq:hatphi}) 
defines an automorphism of $P$. 
We thus obtain a gauge transformation if and only if 
$\phi = \id_M$, which leads to 
$$ \Gau(P) \cong \{ f \in C^\infty(\tilde M,K) \: 
(\forall \gamma \in \pi_1(M)) \ f \circ \gamma = 
c_{\rho(\gamma)}^{-1} \circ f \}, $$
where $c_k(h) := khk^{-1}$. The 
group $\Aut(\tilde M)$ defined under (c) acts naturally on 
$P$, preserving $s(\tilde M)$, via $\tilde\phi(s(\tilde m)) 
:= s(\tilde\phi(\tilde m)).$ An element of this group induces the 
identity on $M$ if and only if it comes from some $\gamma \in Z(\pi_1(M))$, 
and this corresponds to the constant function $\tilde M \to K$ with 
the value $\rho(\gamma)$. We thus 
obtain an open subgroup of $\Aut(P)$ as the quotient 
$$ (\Gau(P) \rtimes \Aut(\tilde M))/Z(\pi_1(M)). $$

Accordingly, the horizontal lift $X \mapsto \tilde X$ of vector 
fields on $M$ to vector fields on $P$ yields an isomorphism 
$$ \aut(P) \cong \gau(P) \rtimes {\cal V}(M), $$
where 
$\gau(P) \cong \{ f \in C^\infty(\tilde M,\fk) \: 
(\forall \gamma \in \pi_1(M))
\ f\circ \gamma = \Ad(\rho(\gamma))^{-1} \circ f\}.$
\end{remark}

\subsection{Mapping groups on non-compact manifolds} 

For a non-compact smooth manifold $M$, the above constructions of the 
atlas on $C^\infty(M,K)$ does no longer work because the sets $U_G$ are 
not open. However, it turns out that in some interesting cases there are 
other ways to obtain charts. To formulate the results, we 
say that a Lie group structure on $G = C^\infty(M,K)$ 
is {\it compatible with evaluations} if 
$\g := C^\infty(M,\fk)$ is its Lie algebra and all evaluation maps 
$\ev_{m}^K \: G \to K$, $m \in M$, are smooth with 
$\L(\ev_m^K) = \ev_m^\fk.$ We then have 
(cf.\ \cite{NeWa07}): 

\begin{theorem} \label{thm:real-lie} Let $K$ be a connected 
regular real Lie group and $M$ 
a real finite-dimensional connected manifold (with boundary). 
Then the group \break $C^\infty(M,K)$ carries a Lie group structure compatible with 
evaluations if one of the following conditions hold: 
\begin{enumerate}
\item[\rm(1)] The universal covering group $\tilde K$ of $K$ 
is diffeomorphic to a locally convex space; which is the case if 
$K$ is finite-dimensional solvable. 
If, in addition, $\pi_1(M)$ is finitely generated, the Lie group structure 
is compatible with the smooth compact open topology. 
\item[\rm(2)] $\dim M = 1$. 
\item[\rm(3)] $M \cong \R^k \times C$, where $C$ is compact. 
\end{enumerate}
\end{theorem}
   
For complex groups we have (\cite{NeWa07}, Thms.~III.12, IV.3): 
\begin{theorem}  \label{thm:cplx-lie} 
Let $K$ be a connected regular complex Lie group and $M$ 
a finite-dimensional connected complex manifold without boundary. 
Then the group ${\cal O}(M,K)$, endowed with the compact open topology, 
carries a Lie group structure with  
Lie algebra ${\cal O}(M,\fk)$ compatible with evaluations if
\begin{enumerate}
\item[\rm(1)] $\tilde K$ is diffeomorphic to a locally convex space. If,  
in addition, $\pi_1(M)$ is finitely generated, the Lie group structure 
is compatible with the compact open topology. 
\item[\rm(2)] $\dim_\C M = 1$, $\pi_1(M)$ is finitely generated 
and $K$ is a Banach--Lie group
\end{enumerate} 
\end{theorem}

\begin{problem} (1) Do the results of Section 3 below on central extensions 
of mapping group $C^\infty(M,K)$ extend to corresponding groups 
of smooth and holomorphic maps on non-compact manifolds 
whenever Theorem~\ref{thm:real-lie} or Theorem~\ref{thm:cplx-lie} provide a 
suitable Lie group structure? For corresponding results on groups of compactly 
supported maps we refer to \cite{Ne04b}. 

(2) Find an example of a connected non-compact smooth 
manifold $M$ and a finite-dimensional Lie group $K$ 
for which $C^\infty(M,K)$ does not carry a Lie group structure compatible 
with evaluations. The first candidate causing problems is the compact 
group $K = \SU_2(\C)$ and $M$ should be a manifold that is not 
a product of some $\R^n$ with a compact manifold 
(cf.\ Theorem~\ref{thm:real-lie}(3)); 
which excludes all 
Lie groups and all Riemannian symmetric spaces. 

(3) Do the preceding theorems generalize in a natural way to gauge 
groups? Which restrictions does one have to impose on the bundles? 
\end{problem}

\section{Central extensions of mapping groups} 
\label{sec:2}

In this section we discuss several issues concerning central extensions of 
mapping groups. We start with a description of the fundamental types of 
$2$-cocycles on the Lie algebra $\g := C^\infty(M,\fk)$, where $M$ is a 
finite-dimensional smooth manifold and $\fk$ a finite-dimensional Lie algebra 
over $\K \in \{\R,\C\}$. These cocycles have values in spaces like 
$C^\infty(M,V)$, $\Omega^1(M,V)$ and 
$\oline\Omega^1(M,V) = 
\Omega^1(M,V)/\dd C^\infty(M,V).$ In particular we describe 
what is known about their integrability to corresponding central Lie group extensions. 
This is best understood for the identity component 
$C^\infty(M,K)_0$, but the extendibility to all of 
$C^\infty(M,K)$ is not clear (cf.\ Problem~\ref{prob:8.3} below). 
In Section~\ref{sec:gaugeext} below we outline some generalizations 
of the results in the present section to gauge groups of non-trivial 
bundles. 

\subsection{Central extensions of $C^\infty(M,\fk)$} 

We write $[\alpha]\in \oline\Omega^1(M,V)$ for 
the image of a $1$-form $\alpha\in \Omega^1(M,V)$ in this space. 
The subspace $\dd C^\infty(M,V)$ of exact $V$-valued $1$-forms 
is characterized by the vanishing of all integrals over loops in $M$, 
hence a closed subspace, and therefore $\oline\Omega^1(M,V)$ carries 
a natural Hausdorff locally convex topology. 

For a trivial $\fk$-module $V$, we write $\Sym^2(\fk,V)^\fk$ for the 
space of $V$-valued symmetric invariant bilinear forms, and recall the 
{\it Cartan map} 
$$ \Gamma \: \Sym^2(\fk,V)^\fk \to Z^3(\fk,V), \quad 
\Gamma(\kappa)(x,y,z) := \kappa([x,y],z). $$
We call $\kappa$ {\it exact} if $\Gamma(\kappa)$ is a coboundary.

\begin{theorem} \label{thm:cocyctypes} Each continuous $2$-cocycle 
with values in the trivial $\g$-module $\K$ is a sum 
of cocycles that factor through a cocycle of one of the following types: 
\begin{description}
\item[\rm(I)] Cocycles $\omega_\kappa \in Z^2(\g, \oline\Omega^1(M,V))$ of the 
form 
$$ \omega_\kappa(\xi_1, \xi_2) = [\kappa(\xi_1,\dd \xi_2)], \quad 
\kappa \in \Sym^2(\fk,V)^\fk, $$
where $\kappa(\xi_1,\dd \xi_2)$ is considered as a $V$-valued 
$1$-form on $M$. 
\item[\rm(II)] Cocycles $\omega_\eta \in Z^2(\g, C^\infty(M,V))$ of the 
form 
$$ \omega_\eta(\xi_1, \xi_2) = \eta(\xi_1, \xi_2) := \eta \circ (\xi_1, \xi_2), \quad 
\eta \in Z^2(\fk,V). $$
\item[\rm(III)] Cocycles $\omega_{\kappa,\eta} \in Z^2(\g, \Omega^1(M,V))$ of the 
form 
$$ \omega_{\kappa,\eta}(\xi_1, \xi_2) 
= \kappa(\xi_1, \dd \xi_2) - \kappa(\xi_2, \dd \xi_1) 
- \dd(\eta(\xi_1, \xi_2)), $$
where $\eta \in C^2(\fk,V)$ and $\kappa \in \Sym^2(\fk,V)^\fk$ 
are related by $\Gamma(\kappa) = \dd_\fk \eta$. 
\item[\rm(IV)] Cocycles vanishing on $\g' \times \g$, 
where $\g' = C^\infty(M,\fk')$ is the commutator algebra, 
i.e., pull-backs of cocycles of an abelian quotient of~$\g$. 
\end{description}
\end{theorem}

\begin{proof}  First we apply 
\cite{NeWa06} (Lemma~1,1, Thm.~3.1, Cor.~3.5 and 
Section~6) to $\g = C^\infty(M,\K)\otimes \fk$ to see that 
after subtracting a cocycle of type (IV), each 
continuous cocycle $\omega \in Z^2(\g,\K)$ can be written as 
\begin{equation}
  \label{eq:coc-deco} 
 \omega(f \otimes x, g \otimes y) =  \beta_a(f\dd g - g \dd f)(x,y) 
-  \beta_s(fg)(x,y), 
\end{equation}
where 
$$\beta_a \: \Omega^1(M,\K) \to \Sym^2(\fk,\K)^\fk \quad \mbox{ and } \quad 
\beta_s \: C^\infty(M,\K) \to C^2(\fk,\K)$$
are continuous liner maps coupled by the condition 
\begin{equation}
  \label{eq:couple}
\Gamma(\beta_a(\dd f)) = \dd_\fk(\beta_s(f)) \quad \mbox{ for } \quad 
f \in C^\infty(M,\K). 
\end{equation}
Conversely, any pair $(\beta_a, \beta_s)$, satisfying 
(\ref{eq:couple}) defines a cocycle via (\ref{eq:coc-deco}): 
\begin{align*}
& (\dd_\g\omega)(f \otimes x, f'\otimes x', f'' \otimes x'')= - \sum_\cyc \omega(ff' \otimes [x,x'], f'' \otimes x'') \\ 
&= - \sum_\cyc \beta_a(ff'\dd f'' - f'' \dd(ff'))([x,x'], x'') + 
\sum_\cyc \beta_s(ff'f'')([x,x'],x'') \\ 
&= \Gamma(\beta_a(\dd(ff'f'')))(x,x',x'') - \dd_\fk(\beta_s(ff'f''))(x,x',x'').
\end{align*}

Let $V(\fk):= S^2(\fk)/\fk.S^2(\fk)$ be the target space of the 
universal invariant symmetric bilinear form 
$\kappa_u(x,y) := [x \vee y]$
 and observe that $\Sym^2(\fk,\K)^\fk \cong V(\fk)^*$. 
We write $\Sym^2(\fk,\K)^\fk_{\rm ex}$ for the subspace of 
$\Sym^2(\fk,\K)^\fk$ consisting of exact forms, and 
$V(\fk)_0 \subeq V(\fk)$ for its annihilator. 
Then there exists a linear map 
$$ \chi \: \Sym^2(\fk,\K)^\fk_{\rm ex} \to C^2(\fk,\K) 
\quad \mbox{with} \quad \dd_\fk(\chi(\beta)) = \Gamma(\beta). $$

As $\beta_a(\dd C^\infty(M,\K))\subeq \Sym^2(\fk,\K)^\fk_{\rm ex}$, we 
may use the Hahn--Banach Extension Theorem to extend $\beta_a$ from 
$\dd C^\infty(M,\K)$ to a 
continuous linear map $\tilde \beta_a \: \Omega^1(M,\K) \to 
\Sym^2(\fk,\K)^\fk_{\rm ex}$ and put 
$\tilde \beta_s(f) := \chi(\tilde \beta_a(\dd f)$, so that 
$$ \tilde\omega(f \otimes x, f' \otimes x') := 
\tilde \beta_a(f\dd f' - f' \dd f)(x,x') - \tilde \beta_s(ff')(x,x') $$
is a $2$-cocycle. 
Further, $\hat \beta_a := \beta_a - \tilde \beta_a$ vanishes on 
$\dd C^\infty(M,\K)$, so that the values of 
$\hat \beta_s := \beta_s - \tilde \beta_s$ are cocycles. Hence both pairs 
$(\hat\beta_a, 0)$, resp., $(0,-\hat\beta_s)$ define cocycles, 
$\omega_2$, resp., $\omega_3$, satisfying 
$$ \omega = \tilde\omega + \omega_2 - \omega_3. $$

We now show that each of these summands factors through a 
cocycle of the form (I)-(III). 
For $\omega_2$ we put $V := V(\fk)$, so that we may write 
$$ \omega_2(f \otimes x, f' \otimes x') 
= \hat\beta_a(f\dd f' - f'\dd f)(x,x') 
= \oline\beta_a(f\dd f' - f'\dd f \otimes \kappa_u(x,x')), $$
where $\oline\beta_a \: \Omega^1(M, V) \cong \Omega^1(M,\K) \otimes V \to \K$ 
is a continuous linear map vanishing on $\dd C^\infty(M,V)$. 
This shows that $\omega_2$ factors through a cocycle of type (I) 
with values in $\oline\Omega^1(M,V)$. 

For $\omega_3$ we put $V := Z^2(\fk,\K)^*$, and write 
$\eta_u \in Z^2(\fk,V)$ for the $2$-cocycle  
defined by $\eta_u(x,x)(f) = f(x,x')$ for $f \in Z^2(\fk,\K)$.  
Then 
$$ \omega_3(f \otimes x, f' \otimes x') 
= \hat\beta_s(ff')(x,x') 
= \eta_u(x,x')(\hat\beta_s(ff')) 
= \oline\beta_s(ff' \otimes \eta_u(x,x')), $$
where $\oline\beta_s \: C^\infty(M,V) \cong C^\infty(M,\K) \otimes V \to \K$ 
is a continuous linear map. 
Hence  $\omega_3$ factors through a cocycle of type (II) 
with $\eta = \eta_u$. 

Finally, we turn to $\tilde\omega$. Now we put $V := V(\fk)/V(\fk)_0$, 
so that $V^* \cong \Sym^2(\fk,\K)^\fk_{\rm ex}$ and write 
$\kappa \: \fk \times \fk \to V$ for the symmetric bilinear map 
obtained from $\kappa_u$. We note that 
$\chi \: V^* \to C^2(\fk,\K)$ 
defines a map $\chi^\sharp \in C^2(\fk,V)$ by 
$$ \chi(\lambda)(x,y) = \lambda(\chi^\sharp(x,y)) \quad \mbox{ for } \quad 
x,y \in \fk, \lambda \in V^*. $$
For each $\lambda \in V^*$ we then identify $\lambda$ with the 
corresponding bilinear form $\lambda \circ \kappa$ and obtain 
$$ \lambda \circ \dd_\fk(\chi^\sharp) 
= \dd_\fk(\lambda \circ \chi^\sharp) 
= \dd_\fk(\chi(\lambda)) = \Gamma(\lambda) 
= \Gamma(\lambda \circ \kappa) 
= \lambda \circ \Gamma(\kappa), $$
showing that $\dd_\fk(\chi^\sharp) = \Gamma(\kappa)$. 

Let $\oline\beta_a \: \Omega^1(M,V) \cong \Omega^1(M,\K) \otimes V \to \K$ 
be the continuous linear functional defined by 
$\oline\beta_a(\alpha \otimes v) = \tilde\beta_a(\alpha)(v).$
Then 
\begin{align*}
& \oline\beta_a\big(\kappa(f \otimes x, \dd f' \otimes x') - 
\kappa(f' \otimes x', \dd f \otimes x) 
- \dd(\chi^\sharp(f \otimes x, f' \otimes x'))\big) \\
&= \oline\beta_a\big((f\dd f' - f' \dd f) \otimes \kappa(x,x') 
- \dd(ff' \otimes \chi^\sharp(x, x'))\big) \\
&= \tilde\beta_a(f\dd f' - f' \dd f)(x,x') 
- \tilde\beta_a(\dd(ff'))(\chi^\sharp(x, x'))\\
&= \tilde\beta_a(f\dd f' - f' \dd f)(x,x') 
- \chi(\tilde\beta_a(\dd(ff')))(x, x')\\
&= \tilde\beta_a(f\dd f' - f' \dd f)(x,x') 
- \tilde\beta_s(ff')(x, x') = \tilde\omega(f \otimes x, f' \otimes x').
\end{align*}
This completes the proof. 
\smartqed\qed
\end{proof}

\begin{remark} \label{rem:3.2} (a) The central extension 
$\hat\g := C^\infty(M,V) \oplus_{\omega_\eta} \g$ defined 
by a cocycle $\omega_\eta$ of type (II) can also be described more directly: If 
$\hat\fk = V \oplus_\eta \fk$ is the central extension defined by $\eta$, then 
$\hat\g \cong C^\infty(M,\hat\fk) \cong C^\infty(M,\K) \otimes \hat\fk$, 
i.e., the $C^\infty(M,\K)$-Lie algebra $\hat\g$ 
is obtained from $\hat\fk$ by extension of scalars from $\K$ to 
the ring $C^\infty(M,\K)$. 

(b) For any $\hat\eta \in C^\infty(M,Z^2(\fk,V))$, 
$$\omega_{\hat\eta}(\xi_1, \xi_2)(m) := \hat\eta(m)(\xi_1(m), \xi_2(m)) $$
also defines a continuous $C^\infty(M,V)$-valued $2$-cocycle on $\g$. 
To see how these cocycles fit into the scheme of the preceding 
theorem, we observe that $Z^2(\fk,V)$ is finite-dimensional 
(if $\fk$ and $V$ are finite-dimensional), so that 
$\hat\eta = \sum_{i = 1}^n \alpha_i \cdot \eta_i$ for 
some $\alpha_i \in C^\infty(M,\R)$ and 
$\eta_i \in Z^2(\fk,V)$. This leads to 
$\omega_{\hat\eta} =\break  \sum_i \alpha_i \cdot \omega_{\eta_i}$, showing 
that $\omega_{\hat\eta}$ factors through the cocycle 
$\omega_{(\eta_1,\ldots, \eta_n)}$ with values in 
$C^\infty(M,V)^n \cong C^\infty(M,V^n)$. 

All the cocycles $\omega_{\hat\eta}$ are 
$C^\infty(M,\R)$-bilinear and if such a cocycle 
is a coboundary, there exists a 
continuous linear map $\beta \: \g \to C^\infty(M,V)$ with 
$\omega_\eta(\xi_1, \xi_2) = \beta([\xi_2, \xi_1])$ for 
$\xi_1, \xi_2 \in \g$. This relation easily implies that the 
restriction of $\beta$ to the commutator algebra 
$\g' = C^\infty(M,\fk')$ is 
$C^\infty(M,\R)$-bilinear so that we may w.l.o.g.\ assume that 
$\beta$ itself has this property. 
Then there exists a smooth map 
$\hat\beta \: M \to \Hom(\fk,V)$ with $\beta(\xi)(m) 
= \hat\beta(m)(\xi(m))$ and we obtain 
$$ \dd_\fk\big(\hat\beta(m)\big) = \hat\eta(m)\quad \mbox{ for } \quad 
m \in M. $$
Thus $\omega_{\hat\eta}$ is a coboundary 
if and only if $\hat\eta(M)$ consists of coboundaries. 
In particular, $[\omega_\eta]= 0$ is equivalent to $[\eta] = 0$ 
if $\hat\eta = \eta$ is constant, which corresponds to cocycles 
of type (II).

(c) Each cocycle $\omega_{\kappa,\eta}$ of type (III) 
may be composed with the quotient map 
$q \: \Omega^1(M,V) \to \oline\Omega^1(M,V)$, which leads to the cocycle 
$q \circ \omega_{\kappa,\eta} = 2\omega_\kappa$ 
of type (I). Here the exactness of $\kappa$ ensures the 
existence of a lift of $\omega_\kappa$ to an $\Omega^1(M,V)$-valued 
cocycle. If $\tilde\eta \in C^2(\fk,V)$ also satisfies 
$\dd_\fk \tilde\eta = \Gamma(\kappa)$, then 
$\tilde\eta - \eta \in Z^2(\fk,V)$ and 
\begin{equation}
  \label{eq:eta-change}
\omega_{\kappa,\tilde\eta} = \omega_{\kappa,\eta} + \dd \circ \omega_{\tilde\eta-\eta} 
\end{equation}
is another lift of $2\omega_\kappa$. 

(d) If $\fk$ is a semisimple Lie algebra and $\kappa \: 
\fk \times \fk \to V(\fk)$ is the universal invariant symmetric bilinear 
form, then the corresponding cocycle $\omega_\kappa$ with values 
in $\oline\Omega^1(M,V(\fk))$ is universal, i.e., up to coboundaries 
each $2$-cocycle with values in a trivial module 
can be written as $f \circ \omega_\kappa$ 
for a continuous linear map $f \: \oline\Omega^1(M,V(\fk)) \to V$ 
(cf.\ \cite{Ma02}, \cite{PS86}). 

(e) Cocycles of type (III) exist if and only if $\fk$ carries an exact 
invariant symmetric bilinear form $\kappa$ with  $\Gamma(\kappa)\not=0$.  
Typical examples of Lie algebras with this 
property are cotangent bundles $\fk = T^*(\fh) := \fh^* \rtimes \fh$ 
with $\kappa((\alpha,x),(\alpha',x')) := \alpha(x') + \alpha'(x)$ 
and $\eta((\alpha,x),(\alpha',x')) := \alpha'(x) - \alpha'(x)$ 
(cf.\ \cite{NeWa06}, Example~5.3). 

(f) Cocycles of the types (I)-(III) can also be defined if $\fk$ is an 
infinite-dimensional locally convex 
Lie algebra. In this case we require $\kappa$ and 
$\eta$ to be continuous. The only point where we have used the finite 
dimension of $\fk$ in the proof of Theorem~\ref{thm:cocyctypes} 
is to show that every $\K$-valued cocycle is a sum of 
cocycles factoring through a cocycle of type (I)-(IV); the corresponding 
arguments do not carry over to the infinite dimensional case. 
\end{remark}

\begin{remark} \label{rem:redux} Many 
results concerning cocycles of type (I)-(III) can be 
reduced from general manifolds to simpler ones as follows. 

For cocycles of type (I), the integration maps
$\int_\gamma \: \oline\Omega^1(M,V) \to V,$ $\gamma \in C^\infty(\bS^1,M),$
separate the points. Accordingly, we have pull-back homomorphisms of 
Lie groups $\gamma^* \: C^\infty(M,K) \to C^\infty(\bS^1,K), f \mapsto 
f \circ\gamma$ satisfying 
\begin{equation}
  \label{eq:red-i}
\int_\gamma \circ \omega_\kappa^M = \int_{\bS^1} \circ \omega_\kappa^{\bS^1} \circ (\L(\gamma^*) \times \L(\gamma^*)),
\end{equation}
which can be used to reduce many things to $M = \bS^1$ 
(cf.\ \cite{MN03}). 

For cocycles of type (II), the evaluation maps 
$$ \ev_m^V  \: C^\infty(M,V)  \to V,\quad 
f \mapsto f(m) $$ 
separate the points and the corresponding evaluation homomorphism 
of Lie groups $\ev_m^K  \: C^\infty(M,K) \to K, f \mapsto f(m),$ satisfies 
\begin{equation}
  \label{eq:red-ii}
\ev_m^V \circ \omega_\eta 
= \eta \circ  (\L(\ev_m^K) \times \L(\ev_m^K)) 
= \eta \circ  (\ev_m^\fk \times \ev_m^\fk). 
\end{equation}

Finally, we observe that for cocycles of type (III), 
the integration maps
$\int_\gamma \: \Omega^1(M,V) \to V, \gamma \in C^\infty(I,M)$, 
separate the points and that the 
pull-back homomorphism $\gamma^* \: C^\infty(M,K) \to C^\infty(I,K)$ 
satisfies 
\begin{equation}
  \label{eq:red-iii}
\int_\gamma \circ \omega_{\kappa,\eta}^M 
= \int_I \circ \omega_{\kappa,\eta}^I \circ 
(\L(\gamma^*) \times \L(\gamma^*)). 
\end{equation}
\end{remark}

\subsection{Covariance of the Lie algebra cocycles} 

In this subsection we discuss the covariance of the cocycles of 
type (I)-(III) under $\Diff(M)$ and $C^\infty(M,K)$, i.e., the 
existence of an actions of these groups on the centrally extended 
Lie algebras $\hat\g = \fz \oplus_\omega \g$, 
compatible with the actions on $\g$ and $\fz$. 
For $\Diff(M)$, the situation is simple because the cocycles 
are invariant, but for $C^\infty(M,K)$ several problems arise. 
For the identity component $C^\infty(M,K)_0$,  
the existence of a corresponding action on $\hat\g$ is equivalent 
to the vanishing of the flux homomorphism 
(Proposition~\ref{prop:fluxcrit}), but for the existence of an 
extension to the full group there is no such simple criterion.

\begin{remark} \label{rem:3.3} (Covariance under $\Diff(M)$) 
(a) The Lie algebra ${\cal V}(M)$ acts in a natural way 
on all spaces $\Omega^p(M,V)$ and $\oline\Omega^p(M,V)$ and it acts on 
the Lie algebra $\g = C^\infty(M,\fk)$ by derivations. With respect to 
this action, all cocycles 
$\omega_\kappa, \omega_\eta$ and $\omega_{\kappa,\eta}$ are 
${\cal V}(M)$-invariant, hence extend by 
$$ \tilde\omega((f,X),(f',X')) := \omega(f,f') $$
to cocycles of the semidirect sum $\g \rtimes {\cal V}(M)$. 

(b) These cocycles are invariant under the full diffeomorphism group 
$\Diff(M)$, which implies that the diagonal action of this group 
on the corresponding central extension defined by 
$\phi.(z,f) := ((\phi^{-1})^*z, \phi.f)$ is a smooth action by Lie algebra 
automorphisms (\cite{MN03}, Thm.~VI.3). 
\end{remark} 

\begin{remark} \label{rem:covmap} (Covariance under $C^\infty(M,K)$) (a) For 
$\omega_\kappa \in Z^2(\g,\oline\Omega^1(M,V))$ of type (I) 
we also have a natural action of the full group $C^\infty(M,K)$ on the 
extended Lie algebra $\hat\g = \oline\Omega^1(M,V) \oplus_{\omega_\kappa} \g$ 
(\cite{MN03}, Thm.~VI.3). 
If $\delta(f) := \delta^l(f)$ denotes the left logarithmic derivative of $f$, 
then $\theta(f)(x) := [\kappa(\delta(f),x)]$ defines a linear map 
$\g \to \oline\Omega^1(M,V)$, satisfying 
$$ \Ad(f)^*\omega_\kappa - \omega_\kappa = \dd_\g(\theta(f)), $$
so that 
$$ \Ad_{\hat\g}(f).(z,\xi) := (z - \theta(f)(\xi), \Ad(f)\xi) $$
defines an automorphism of $\hat\g$ (Lemma~\ref{lem:aut-lift}), 
and we thus obtain 
a smooth action of the Lie group $C^\infty(M,K)$ on $\hat\g$. 

(b) Let $K$ be a connected Lie group with Lie algebra $\fk$, $V$ a trivial 
$\fk$-module and $\eta \in Z^2(\fk,V)$ be a $2$-cocycle. 
For the cocycle $\omega_\eta \in Z^2(\g,C^\infty(M,V))$ of type (II), 
the situation is slightly more complicated. We recall from 
Proposition~\ref{prop:fluxcrit} below that the 
action of $K$ on $\fk$ lifts to an action on $\hat\fk = V \oplus_\eta \fk$ 
if and only if the flux 
$$F_\eta \: \pi_1(K)\to \Hom_{\rm Lie}(\fk,V) = H^1(\fk,V) $$
vanishes. 

Let $m_0 \in M$ be a base point and 
$C^\infty_*(M,K) \trile C^\infty(M,K)$ be the normal subgroup of all maps 
vanishing in $m_0$. Then this also is a Lie group, and 
$$C^\infty(M,K) \cong C^\infty_*(M,K) \rtimes K $$
as Lie groups. If $q_K \: \tilde K \to K$ denotes the universal covering group, 
then a smooth map $f \: M \to K$ vanishing in $m_0$ lifts to a map 
$M \to \tilde K$ if and only if the induced homomorphism 
$\pi_1(f) \: \pi_1(M)\to \pi_1(K)$ vanishes. We therefore have an exact 
sequence of groups 
\begin{equation}
  \label{eq:maplift}
\1 \to C^\infty_*(M,\tilde K) \to C^\infty_*(M,K)  
\to \Hom(\pi_1(M), \pi_1(K)),
\end{equation}
and since $\pi_1(M)$ is finitely generated because $M$ is compact, we may 
consider $C^\infty_*(M,\tilde K)$ as an open subgroup of $C^\infty(M,K)$. 
From the canonical action $\Ad_{\hat\fk}$ 
of $\tilde K$ on $\hat\fk$ we immediately 
get a smooth action of $C^\infty(M,\tilde K)$ on 
$\hat\g \cong C^\infty(M,\hat\fk)$ by 
$$(\Ad_{\hat\g}(f).\xi)(m) := \Ad_{\hat \fk}(f(m)).\xi(m). $$
This shows that the flux 
$$ F_{\omega_\eta} \: \pi_1(C^\infty(M,K)) 
\cong \pi_1(C^\infty_*(M,K)) \times \pi_1(K) \to 
H^1(\g, C^\infty(M,V)) $$
vanishes on the subgroup  $\pi_1(C^\infty_*(M,K))$, hence factors through 
the flux $F_\eta$ of $\eta$. 

For any $f \in C^\infty(M,K)$ and $m \in M$, the connectedness of $K$ 
implies that $\Ad(f(m))^*\eta - \eta$ is a $2$-coboundary 
(Theorem~\ref{thm:9.1}), and from that 
we derive the existence of a smooth map 
$\tilde\theta \: M \to C^1(\fk,V)$ with 
$\Ad(f(m)^{-1})^*\eta - \eta = \dd_\fk(\theta(m))$ for $m \in M$. 
Now the automorphism $\Ad(f)$ lifts to an 
automorphism of $\hat\g$ (Lemma~\ref{lem:aut-lift}).

(c) To determine the covariance of the cocycles 
$\omega_{\kappa,\eta}$ of type (III) under the group $C^\infty(M,K)$, 
we first observe that for $f \in C^\infty(M,K)$ and 
$\xi_1, \xi_2 \in \g$ we obtain similarly as in \cite{MN03}, Prop.~III.3: 
\begin{align*}
&\big(\Ad(f)^*\omega_{\kappa,\eta}\big)(\xi_1, \xi_2)=\\
& \kappa(\Ad(f).\xi_1, \dd(\Ad(f).\xi_2))-\kappa(\Ad(f).\xi_2, \dd(\Ad(f).\xi_1)) - \dd\big((\Ad(f)^*\eta)(\xi_1, \xi_2)\big) \\ 
&= \kappa(\xi_1,\dd\xi_2) - \kappa(\xi_2,\dd\xi_1) 
+ 2 \kappa(\delta^l(f), [\xi_2, \xi_1]) 
 - \dd\big((\Ad(f)^*\eta)(\xi_1, \xi_2)\big) \\ 
&= \omega_{\kappa,\eta}(\xi_1, \xi_2) 
+ 2 \kappa(\delta^l(f), [\xi_2, \xi_1])
 - \dd\big((\Ad(f)^*\eta-\eta)(\xi_1, \xi_2)\big). 
\end{align*}
It is clear that the 
term $\kappa(\delta^l(f), [\xi_2, \xi_1])$ is a Lie algebra
coboundary. To see that the second term also has this property, 
we note that $\Gamma(\kappa) = \dd_\fk\eta$ implies that for 
ech $x \in \fk$ the cochain $i_x(\dd_\fk\eta)(y,z) 
= \kappa([x,y],z) = \kappa(x,[y,z])$ 
is a coboundary. Hence ${\cal L}_x \eta = i_x\dd_\fk\eta + \dd_\fk(i_x\eta)$ 
also is a coboundary and since $K$ is connected, 
for each $m \in M$ the cocycle $\Ad(f(m))^*\eta - \eta$ is a coboundary. 
With Lemma~\ref{lem:aut-lift} we now see that for 
each smooth function $f \: M \to K$, the adjoint action on $\g$ 
lifts to an automorphism of the central extension $\hat\g$, defined by 
$\omega_{\kappa,\eta}$. If $f$ is not in $C^\infty(M,K)_0$, i.e., homotopic 
to a constant map, this does not follow from Theorem~\ref{thm:9.1}. 
\end{remark}

\begin{problem} \label{prob:2.7} For 
which groups $K$ and which connected smooth manifolds 
$M$ does the short exact sequence 
$$\1 \to C^\infty(M,K)_0  \to 
C^\infty(M,K) \to \pi_0(C^\infty(M,K)) \cong [M,K] \to \1 $$
split? 

If $K\cong \fk/\Gamma_K$ is abelian, then 
$C^\infty(M,K)_0$ is a quotient of $C^\infty(M,\fk)$, hence divisible, and therefore 
the sequence splits. 

If $K$ is non-abelian, the situation is more involved and even the case $M = \bS^1$ 
is non-trivial. Then $[M,K] \cong \pi_1(K)$ and if $K$ is abelian, then 
we obtain a splitting from the isomorphism $\pi_1(K) \cong \Hom(\T,K)$. 
If $K$ is solvable and $T \subeq K$ is a maximal torus, then 
$K \cong \R^m \times T$ as a manifold (\cite{Ho65}), 
so that $[M,K] \cong [M,T]$ and we obtain a splitting from the 
splitting $[M,T] \to C^\infty(M,T) \subeq C^\infty(M,K)$. 
The case of compact (semisimple) groups is the crucial one. 

If $M$ is a product of $d$ spheres (f.i.\ a $d$-dimensional torus), then 
$[M,K]$ is nilpotent of length $\leq d$ which is not always abelian 
(\cite{Whi78}, Th.~X.3.6).  
\end{problem}

\begin{problem} Can the automorphisms of the central extensions 
$\hat\g = \fz \oplus_\omega \g$ of type (I)---(III) 
corresponding to elements of $C^\infty(M,K)$ be chosen in such a way 
as to define a group action? 
If $\fk$ and hence $\g$ is perfect, this 
follows trivially from the uniqueness of lifts to $\hat\g$. 
In general, lifts are only unique up to an element of 
of the group $H^1(\g,\fz)$, so that only an abelian 
extension acts. Are the corresponding cohomology classes in 
$H^2(C^\infty(M,K), \fz)$ arising  this way always trivial? 
Since the action of $\g$ on $\hat\g$ is determined 
by the cocycle, we obtain an action of the identity component 
$C^\infty(M,K)_0$ whenever the flux vanishes, and then the 
restriction of the cocycle to the identity component is trivial. 

Here even the case $M = \bS^1$ of loop groups is of particular interest. 
Then \break $\pi_0(C^\infty(\bS^1,K)) \cong \pi_1(K)$, so that 
we have a natural inflation map 
$$ H^2(\pi_1(K),\fz) \to H^2(C^\infty(\bS^1,K), \fz),$$ 
and it would be nice if all classes under consideration can be determined 
this way. 

If we have a semidirect decomposition 
$C^\infty(M,K) \cong C^\infty(M,K)_0 \rtimes [M,K]$, 
then the present problem simplifies significantly because then   
semidirect product techniques similar to those discussed in Appendix B  
for Lie algebras apply. 
\end{problem}

\begin{problem} Generalize the description of $H^2(A \otimes \fk,\K)$ 
obtained in \cite{NeWa06} under the assumption that $A$ is unital 
to non-unital commutative associative algebras. 
This would be of particular interest for algebras of the form 
$C^\infty_*(M,\R)$ (functions vanishing in one point) or the algebra 
$C^\infty_c(M,\R)$ of compactly supported functions. 

If $A_+ := A \oplus \K$ is the algebra with unit $\1 = (0,1)$, 
then we have 
$$ A_+ \otimes \fk \cong (A \otimes \fk) \rtimes \fk, $$
so that we may use the semidirect product techniques in Appendix~B below. 
\end{problem}

\subsection{Corresponding Lie group extensions} 

It is a natural question to which extent 
the cocycles of the form $\omega_\kappa$, $\omega_\eta$ and 
$\omega_{\kappa,\eta}$ on $\g = C^\infty(M,\fk)$ 
actually define central extensions of the corresponding group 
$C^\infty(M,K)$. 

Cocycles of the form $\omega_\kappa$ for a vector-valued 
$\kappa \: \fk \times \fk \to V$ are treated in \cite{MN03}, 
where it is shown that the corresponding period homomorphism 
$$ \per_{\omega_\kappa} \: \pi_2(C^\infty(M,K)) \to \oline\Omega^1(M,V) $$
(cf.\ Appendix A) 
has values in the subspace $H^1_{\rm dR}(M,V)$, and the image consists 
of cohomology classes whose integrals over circles take values in the 
image of the homomorphism 
$$ \per_\kappa \: \pi_3(K) \to V, \quad 
[\sigma] \mapsto \int_\sigma \Gamma(\kappa)^{\rm eq},$$
where $\Gamma(\kappa)^{\rm eq} \in \Omega^3(K,V)$ is the left invariant 
closed $3$-form whose value in $\1$ is $\Gamma(\kappa)$. 
The flux homomorphism 
$$ F_\omega \:\pi_1(C^\infty(M,K)) \to H^1(\g, \oline\Omega^1(M,V)) $$
vanishes because the action of $\g$ on the central extension 
$\hat\g$ defined by $\omega_\kappa$ integrates to a smooth action of the 
group $C^\infty(M,K)$ on $\hat\g$ (Remark~\ref{rem:covmap}(a), 
\cite{MN03}, Prop.~III.3). 
This leads to the following theorem (\cite{MN03}, Thms.~I.6, II.9, Cor.~III.7): 

\begin{theorem} \label{thm:int-type1} 
For the cocycles $\omega_\kappa$ and the connected group 
$G := C^\infty(M,K)_0$, the following are equivalent: 
\begin{description}
\item[\rm(1)] $\omega_{\kappa}$ integrates for each compact manifold $M$ 
to a Lie group extension of $G$. 
\item[\rm(2)] $\omega_{\kappa}$ integrates for $M = \bS^1$ 
to a Lie group extension of $G$. 
\item[\rm(3)] The image of $\per_{\omega_{\kappa}}$ in 
$H^1_{\rm dR}(M,V) \subeq \oline\Omega^1(M,V)$ 
is discrete. 
\item[\rm(4)] The image of $\per_{\kappa}$ in $V$ is discrete. 
\end{description}
These conditions are satisfied if $\kappa$ is the universal invariant 
bilinear form with values in $V(\fk)$. 
\end{theorem}

\begin{remark} Composing a cocycle $\omega_\kappa$ with the exterior derivative, 
we obtain the $\Omega^2(M,V)$-valued cocycle 
$$ (\dd \circ \omega_\kappa)(\xi_1,\xi_2) 
= \kappa(\dd \xi_1, \dd \xi_2). $$

In view of \cite{MN03}, Thm.~III.9, we have a smooth $2$-cocycle 
$$ c(f_1, f_2) 
:=  \delta^l(f_1) \wedge_\kappa \delta^r(f_2)
=  \delta^r(f_1) \wedge_\kappa \Ad(f_1).\delta^r(f_2) $$
on the full group $C^\infty(M,K)$ which defines a central extension 
by $\Omega^2(M,V)$ whose corresponding Lie algebra cocycle is 
$2\dd \circ \omega_\kappa$. Here $\wedge_\kappa$ denotes the natural product 
$\Omega^1(M,\fk) \times \Omega^1(M,\fk) \to \Omega^2(M,V)$ 
defined by $\kappa\: \fk \times \fk \to V$. 
In Subsection~\ref{sec:4.2} below, this construction is generalized 
to gauge groups of non-trivial bundles. 
\end{remark}

For the cocycles of type (II), the situation is much simpler 
(Theorem~\ref{thm:abext}): 

\begin{theorem} \label{thm:int-type2} Let $K$ be a connected finite-dimensional 
Lie group with Lie algebra 
$\fk$, $\eta \in Z^2(\fk,V)$, 
and $\hat\fk = V \oplus_\eta \fk$. 
Then the following are equivalent: 
\begin{itemize}
\item[\rm(a)] The central Lie algebra extension 
$\hat\g = C^\infty(M,V) \oplus_{\omega_\eta} \g$ 
of $\g = C^\infty(M,\fk)$ integrates to a central Lie group extension 
of $C^\infty(M,K)_0$. 
\item[\rm(b)] $\hat\fk \to \fk$ integrates to a Lie group 
extension $V\into \hat K \to K$. 
\item[\rm(c)] The flux homomorphism $F_\eta \: \pi_1(K) \to H^1(\fk,V)$ 
vanishes. 
\end{itemize}
If this is the case, then 
$$ \1 \to C^\infty(M,Z) \to C^\infty(M,\hat K) \to C^\infty(M, K)\to \1 $$
defines a central Lie group extension of the full group 
$C^\infty(M,K)$ integrating~$\hat\g$. 
\end{theorem}

\begin{proof} For $m \in M$, the relation  
$\ev_m^V \circ \omega_{\eta}  = \L(\ev_m^K)^*\eta = (\ev_m^\fk)^*\eta$
((\ref{eq:red-ii}) in Remark~\ref{rem:redux}) 
leads to 
$\ev_m^V \circ \per_{\omega_{\eta}}  = \per_\eta \circ \pi_2(\ev_m^K)$
(see (\ref{eq: per-trans}) in Remark~\ref{rem:8.4}). 
Since $\pi_2(K)$ vanishes (\cite{CaE36}), the period map $\per_\eta$ 
is trivial and therefore all maps 
$\ev_m^V \circ \per_{\omega_{\eta}}$ vanish, which implies that 
$\per_{\omega_\eta} = 0$. 

A direct calculation further shows that the flux homomorphisms 
$F_\eta$ and $F_{\omega_\eta} \: \pi_1(G) \to H^1(\g,C^\infty(M,V))$ 
satisfy 
$$ F_{\omega_\eta}([\gamma])(\xi)(m) 
= F_\eta([\ev_m^K \circ \gamma])(\xi(m)), \quad \xi \in \g, [\gamma]\in \pi_1(G) $$
(see (\ref{eq: flux-trans}) in Remark~\ref{rem:8.4}). 
Therefore $F_\eta$ vanishes if and only if 
$F_{\omega_\eta}$ does. Now Theorem~\ref{thm:abext} below shows that 
(a) is equivalent to (b) which in turn is equivalent to (c) 
because all period maps vanish. 

If these conditions are satisfied, we apply loc.\ cit.\ 
with $V = A$ to obtain a central extension $V \into \hat K \to K$. 
This extension is a principal $V$-bundle, hence trivial as such, so that 
there exists a smooth section $\sigma \: K \to \hat K$. 
Therefore $C^\infty(M,\hat K)$ is a central Lie group extension of 
$C^\infty(M,K)$ with a smooth global section. 
\smartqed\qed
\end{proof}

\begin{remark} If, in the context of Theorem~\ref{thm:int-type2}, 
$q_K \: \hat K \onto K$ is a connected Lie group extension 
of $K$ by $Z = V/\Gamma_Z$, $\Gamma_Z \subeq V$ a discrete subgroup, 
and $\L(\hat K) = \hat\fk = V \oplus_\eta \fk$, then 
$\hat K$ does in general not have a smooth global section and 
the image of the canonical map 
$C^\infty(M,\hat K) \to C^\infty(M,K)$ only is an open subgroup. 

In our context, where $K$ is finite-dimensional, 
the obstruction for a map $f \: M \to K$ to lift to $\hat K$ 
can be made quite explicit. 
The existence of a smooth lift $\hat f \: M \to \hat K$ of $f$ is 
equivalent to the triviality of the smooth $Z$-bundle 
$f^*\hat K \to M$, and the equivalence classes of these bundles 
are parametrized by $H^2(M,\Gamma_Z)$ (cf.\ \cite{Bry93}). 
Describing bundles in terms of \v Cech cocycles, it is easy to see that 
we thus obtain a group homomorphism 
$$ C^\infty(M,K) \to H^2(M,\Gamma_Z), \quad f \mapsto [f^*\hat K] $$
which factors through an injective homomorphism 
$$ \pi_0(C^\infty(M,K)) \cong [M,K] \to H^2(M,\Gamma_Z), \quad 
[f] \mapsto [f^*\hat K] $$
of discrete groups. 

Since $\pi_2(K)$ and $F_\eta$ vanish, 
Remark~6.12 in \cite{Ne04a} implies that the homology class 
$[\hat K] \in H^2(K,\Gamma_Z)$ vanishes on $H_2(K)$, so that the 
Universal Coefficient Theorem shows that it is 
defined by an element of $\Ext(H_1(K),\Gamma_Z) = \Ext(\pi_1(K),\Gamma_Z)$, 
we may interprete $[f^*\hat K]$ as an element of 
the group \break $\Ext(H_1(M),\Gamma_Z)$. The long exact homotopy sequence of the  
$Z$-bundle $\hat K$ contains the short exact sequence 
$$ \1 \to \Gamma_Z = \pi_1(Z) \to \pi_1(\hat K) \to \pi_1(K) \to \1, $$
which exhibits $\pi_1(\hat K)$ as a central extension of $\pi_1(K)$ 
corresponding to the class in $\Ext(\pi_1(K), \Gamma_Z)$. 

We now see that the homomorphism from above yields an exact sequence 
of groups 
$$ C^\infty(M,\hat K) \to C^\infty(M,K) \to \Ext(H_1(M),\Gamma_Z). $$
Here the rightmost map can be calculated by first assigning to 
$f \in C^\infty(M,K)$ the homomorphism 
$H_1(f) \: H_1(M) \to H_1(K) \cong \pi_1(K)$ and then 
$H_1(f)^*[\hat K] \in \Ext(H_1(M),\Gamma_Z)$, which can be evaluated 
easily in concrete cases. 
\end{remark}

Now we turn to cocycles of type (III). We start with a key example: 

\begin{example} \label{ex:pathgrp} 
Let $\fk$ be a locally convex real Lie algebra 
and consider the locally convex Lie algebra $\g := C^\infty(I,\fk)$, where 
$I := [0,1]$ is the unit interval. 

Let $\kappa \: \fk \times \fk \to V$
be a continuous invariant symmetric bilinear form with values in the 
Mackey complete space $V$ and consider $V$ as a 
trivial $\g$-module. 
If $\Gamma(\kappa) = \dd_\fk \eta$ for some $\eta \in C^2(\fk,V)$, then the 
cocycle $\omega_{\kappa,\eta} \in Z^2(\g, \Omega^1(I,V))$ 
vanishes on the subalgebra $\fk \subeq \g$ of constant maps. 

Let $K$ be a connected Lie group with Lie algebra $\fk$ and  
$G := C^\infty(I,K)$. The map 
$H \: I \times G \to G, H(t,f)(s) := f(ts)$
is smooth because the corresponding map 
$\tilde H \: I \times G \times I \to K, \tilde H(t,f,s) := f(ts)$ 
is smooth (cf.\ \cite{NeWa07}, Lemma~A.2). 
Since $H_1 = \id_G$ and $H_0(f) = f(0)$, 
$H$ is a smooth retraction of $G$ to the subgroup 
$K$ of constant maps. Therefore the inclusion 
$j \: K \to G$ induces an isomorphism 
$\pi_2(j) \: \pi_2(K) \to \pi_2(G).$
The period map $\per_{\omega_{\kappa,\eta}} \: \pi_2(G) \to V$ satisfies 
$\per_{\omega_{\kappa,\eta}} \circ \pi_2(j) = \per_{j^*\omega_{\kappa,\eta}} = 0$
because the cocycle $j^*\omega_{\kappa,\eta} = \omega_{\kappa,\eta}\res_{\fk \times \fk}$ vanishes, so that the period group of $\omega_{\kappa,\eta}$ 
is trivial. 
\end{example}

\begin{theorem} Let $M$ be a smooth compact manifold, $\fk$ 
a locally convex 
Lie algebra, $\kappa \: \fk \times \fk \to V$ a 
continuous invariant symmetric bilinear 
form and $\eta \in C^2(\fk,V)$ with $\Gamma(\kappa) = \dd_\fk \eta.$ 
For a connected Lie group $K$ with Lie algebra 
$\fk$ and $G := C^\infty(M,K)$, the following assertions hold: 
\begin{description}
\item[\rm(a)] The period map 
$\per_{\omega_{\kappa,\eta}} \: \pi_2(G) \to \Omega^1(M,V)$ vanishes. 
\item[\rm(b)] If $K$ is simply connected, then 
the flux $F_{\omega_{\kappa,\eta}}$ vanishes. 
\item[\rm(c)] If $K$ is finite-dimensional, then there exists 
$\tilde\eta \in C^2(\fk,V)$ 
with $\dd_\fk\tilde\eta = \Gamma(\kappa)$ for which the 
flux $F_{\omega_{\kappa,\tilde\eta}}$ vanishes.
\item[\rm(d)] If $F_{\omega_{\kappa,\eta}}=0$, 
then $\omega_{\kappa,\eta}$ integrates to a 
central extension of $G_0$ by $\Omega^1(M,V)$. 
\end{description}
\end{theorem}

\begin{proof} (a) Combing (\ref{eq:red-iii}) in Remark~\ref{rem:redux} 
with (\ref{eq: per-trans}) in Remark~\ref{rem:8.4}, 
we obtain for each $\gamma \in C^\infty(I,M)$: 
$$ \int_\gamma \circ \per_{\omega_{\kappa,\eta}^M} 
= \per_{\int_\gamma \circ \omega_{\kappa,\eta}^M} 
= \per_{\int_I \circ \omega_{\kappa,\eta}^I} \circ \pi_2(\gamma^*)
= \int_I \circ \per_{\omega_{\kappa,\eta}^I} \circ \pi_2(\gamma^*) $$
and Example~\ref{ex:pathgrp} implies that 
$\per_{\omega^I_{\kappa,\eta}}$ vanishes. Since $\gamma$ was arbitrary, 
all periods of $\omega_{\kappa,\eta}^M$ vanish. 

(b) Next we assume that $K$ is $1$-connected. 
In Example~\ref{ex:pathgrp} we have seen that the group 
$C^\infty(I,K)$ is also $1$-connected, so that the flux of 
$\omega_{\kappa,\eta}^I$ vanishes. 
Now we apply (\ref{eq: flux-trans}) in Remark~\ref{rem:8.4} 
with (\ref{eq:red-iii}) in Remark~\ref{rem:redux} to see that 
the flux $F_{\omega_{\kappa,\eta}}$ also vanishes. 
 
(c) Assume that $K$ is finite-dimensional and pick a maximal compact 
subgroup $C \subeq K$. Since $\kappa$ is $K$-invariant, the affine action of 
the compact group $C$ on the affine space 
$\{ \eta \in C^2(\fk,V) \: \dd_\fk \eta = \Gamma(\kappa)\}$ 
has a fixed point $\tilde\eta$. Then $\omega_{\kappa,\tilde\eta}$ is $C$-invariant, 
and therefore $C$ acts diagonally on 
$\hat\g = \Omega^1(M,V) \oplus_{\omega_{\kappa,\tilde\eta}} \g$. It follows 
in particular that the flux vanishes on 
the image of $\pi_1(C) \cong \pi_1(K)$ in $\pi_1(G)$ 
(Proposition~\ref{prop:fluxcrit}). 
For $M = I = [0,1]$ we immediately derive that the flux vanishes and in the 
general case we argue as in (b) by reduction via pull-back maps 
$\gamma^* \: C^\infty(M,K) \to C^\infty(I,K)$. 

(d) follows  from Theorem~\ref{thm:abext}. 
\smartqed\qed
\end{proof}

\begin{remark} We have seen in the preceding theorem that if 
$K$ is simply connected, the period homomorphism and the flux of 
$\omega := \omega_{\kappa,\eta}$ vanish. Since $G := 
C^\infty(I,K)$ is also $1$-connected, 
one may therefore expect an explicit formula for a corresponding group cocycle 
integrating $\omega_{\kappa,\eta}$ because the cohomology class  
$[\omega_{\kappa,\eta}^{\rm eq}] 
\in H^2_{\rm dR}(C^\infty(I,K),C^\infty(I,\R))$ vanishes 
(cf.\ \cite{Ne02}, Thm.~8.8).  
Such a formula can be obtained by a method due to Cartan, 
combined with fiber integration and the homotopy $H \: I \times G \to G$ 
to $K$. Since $\omega$ vanishes on $\fk$, the fiber integral 
$\theta := \fint_I H^*\omega^{\rm eq} \in \Omega^1(G,V)$ satisfies 
$\dd\theta = \omega^{\rm eq}$, so that we may use 
Prop.~8.2 in \cite{Ne04a} to construct an explicit cocycle. 
For a more general method, based on path groups, we refer to \cite{Vi07}. 
\end{remark}

\begin{remark} In view of Remarks~\ref{rem:3.3} and 
\ref{rem:covmap}, for cocycles $\omega$ of type (I)-(III), the 
diffeomorphism group $\Diff(M)$ acts diagonally on $\hat\g$. 
If, in addition, $K$ is $1$-connected, then the flux 
$F_\omega$ vanishes, so that the identity component 
$G_0 := C^\infty(M,K)_0$ also acts naturally on $\hat\g$ 
(Proposition~\ref{prop:fluxcrit}), and this action is uniquely 
determined by the adjoint action of $\g$ on $\hat\g$ 
(cf.\ \cite{GN07}). From that is easily follows that 
these two actions combine 
to a smooth action of $C^\infty(M,K)_0 \rtimes \Diff(M)$ on $\hat\g$ 
and the Lifting Theorem \ref{thm:lifting-theorem} provides a smooth 
action of this group by automorphisms on any corresponding central 
extension $\hat G$ of the simply connected covering group $\tilde G_0$. 
\end{remark}

\begin{problem}
Do the central extensions of the connected group $C^\infty(M,K)_0$ 
constructed above extend to the full group? 

The answer is positive 
for $M = \bS^1$, $K$ compact simple and $\omega = \omega_\kappa$ if 
$\kappa$ is universal (\cite{PS86}, Prop.~4.6.9). 
See Problem~\ref{prob:8.3} 
for the general framework in which this problem can be investigated. 
\end{problem}

\section{Twists and the cohomology of vector fields} 
\label{sec:3}

We have seen above how to obtain central extensions of the Lie algebra  
$C^\infty(M,\fk)$ and that the three types of 
cocycles $\omega_\kappa$, $\omega_\eta$ and $\omega_{\kappa,\eta}$ 
are ${\cal V}(M)$-invariant, which leads to central extensions  
of the semidirect sum $C^\infty(M,\fk) \rtimes {\cal V}(M)$ by spaces 
of the type $\fz := C^\infty(M,V), \Omega^1(M,V)$, and 
$\oline\Omega^1(M,V)$. The bracket on such a central extension has the 
form 
\begin{align*}
&[(z_1,(\xi_1,X_1)), (z_2,(\xi_2,X_2))] \\
&= \big(X_1.z_2 - X_2.z_1 + \omega(\xi_1,\xi_2), 
(X_1.\xi_2 - X_2.\xi_1 + [\xi_1,\xi_2], [X_1,X_2])\big), 
\end{align*}
and there are natural twists of these central extensions that correspond 
to replacing the cocycle 
$\hat\omega((\xi_1,X_1), (\xi_2, X_2)) := \omega(\xi_1,\xi_2)$
by $\hat\omega + \hat\eta$, where 
$$ \hat\eta((\xi_1,X_1), (\xi_2, X_2)) := \eta(X_1,X_2) 
\quad \mbox{ for some } \quad 
\eta \in Z^2({\cal V}(M),\fz) $$
(cf.\ Appendix B on extensions of semidirect products). 
To understand the different types of twists, one has to 
determine the cohomology groups 
$$ H^2_c({\cal V}(M),\R), \quad H^2_c({\cal V}(M),\Omega^1(M,\R)) \quad 
\mbox{ and } \quad H^2_c({\cal V}(M),\oline\Omega^1(M,\R)). $$
For a parallelizable manifold $M$, these spaces have been determined 
in \cite{BiNe07}, and below we describe the different types of cocycles
showing up. 

We also discuss the integrability of these twists to abelian 
extensions of $\Diff(M)$, or at least of the 
simply connected covering $\tilde \Diff(M)_0$ of its identity component 
$\Diff(M)_0$.

\subsection{Some cohomology of the Lie algebra of vector fields} 

The most obvious source of $2$-cocycles of ${\cal V}(M)$ with values 
in differential forms is described in the following lemma, applied to $p = 2$ 
(\cite{Ne06c}, Prop.~6): 
\begin{lemma} \label{lem:transfer} 
For each closed $(p+q)$-form $\omega \in \Omega^{p+q}(M,V)$,  
the prescription 
$$ \omega^{[p]}(X_1,\ldots, X_p) 
:= [i_{X_p} \ldots i_{X_1} \omega] \in \oline\Omega^q(M,V)$$
defines a continuous $p$-cocycle in $Z^p({\cal V}(M),\oline\Omega^q(M,V))$. 
\end{lemma}

The preceding lemma is of particular interest for $p = 0$. In this case 
it associates  Lie algebra cohomology classes to closed differential forms. 
According to Cor.~3.2 in \cite{Lec85}, we have the following theorem: 

\begin{theorem} \label{thm:lec1} Each smooth $p$-form $\omega \in \Omega^p(M,\R)$ 
defines a $p$-linear alternating continuous 
map ${\cal V}(M)^p\to C^\infty(M,\R)$, which leads to 
an inclusion of chain complexes 
$$ (\Omega^\bullet(M,\R), \dd) \into (C^\bullet({\cal V}(M), C^\infty(M,\R)), \dd_{{\cal V}(M)}) $$
inducing an algebra homomorphism 
$$ \Psi \: H^\bullet_{\rm dR}(M,\R) 
\to H^\bullet({\cal V}(M), C^\infty(M,\R)) $$ 
whose kernel is the ideal generated by the 
Pontrjagin classes of $M$, i.e., the characteristic classes 
of the tangent bundle. 
\end{theorem} 

\begin{remark} There are several classes of compact manifolds 
for which the Pontrjagin classes all vanish but the tangent bundle 
is non-trivial. According to \cite{GHV72}, II, Proposition~9.8.VI, 
this is in particular the case for all Riemannian manifolds with constant 
curvature, hence in particular for all spheres. 
\end{remark}

A second source of cocycles is the following 
(\cite{Kosz74}):  
\begin{lemma} \label{lem:koszul} 
Any affine connection $\nabla$ on $M$ defines a $1$-cocycle 
$$ \zeta \: {\cal V}(M) \to \Omega^1(M,\End(TM)), 
\quad X \mapsto {\cal L}_X\nabla, $$
where 
$({\cal L}_X\nabla)(Y)(Z) 
= [X,\nabla_Y Z] - \nabla_{[X,Y]}Z - \nabla_Y[X,Z].$
For any other affine connection $\nabla'$, the corresponding cocycle $\zeta'$ 
has the same cohomology class. 
\end{lemma} 

\begin{definition} We use the cocycle $\zeta$, 
associated to an affine connection 
$\nabla$ to define $k$-cocycles $\Psi_k \in Z^k_c({\cal V}(M), 
\Omega^k(M,\R))$. Note that $A \mapsto \Tr(A^k)$ 
defines a homogeneous polynomial of 
degree $k$ on $\gl_d(\R)$, invariant under conjugation. 
The corresponding invariant symmetric 
$k$-linear map is given by 
$$ \beta(A_1,\ldots, A_k) = \sum_{\sigma \in S_k} \Tr(A_{\sigma(1)}\cdots 
A_{\sigma(k)}), $$
and we consider it as a linear $\GL_d(\R)$-equivariant map 
$\gl_d(\R)^{\otimes k} \to \R,$
where $\GL_d(\R)$ acts trivially on $\R$. This map leads to a vector bundle map 
$$ \beta_M \: \End(TM)^{\otimes k} \to M \times \R, $$
where $M \times \R$ denotes the trivial vector bundle with fiber $\R$. 
On the level of bundle-valued differential forms, 
this in turn yields an alternating $k$-linear $\Diff(M)$-equivariant 
map $$ \beta_M^1 \: \Omega^1(M,\End(TM))^{k} \to \Omega^k(M,\R). $$
The $\Diff(M)$-equivariance implies the invariance of this map 
under the natural action of ${\cal V}(M)$, so that 
we can use $\beta^1_M$ to multiply Lie algebra cocycles 
(cf.\ \cite{Fu86} or App.~F in \cite{Ne04a}). In 
particular, we obtain for each $k \in \N$ a Lie algebra cocycle 
$$ \Psi_k \in Z^k_c({\cal V}(M), \Omega^k(M,\R)), $$
defined by 
$$ \Psi_k (X_1, \ldots, X_k) 
:= (-1)^k \sum_{\sigma\in S_k} \sgn(\sigma)
\beta_M^1(\zeta(X_{\sigma(1)}), \ldots, \zeta(X_{\sigma(k)})). $$

For any other affine connection $\nabla'$, the corresponding 
cocycle $\zeta'$, and the associated cocycles $\Psi_k'$, 
the difference $\zeta' - \zeta$ is a coboundary,  and since 
products of cocycles and coboundaries are coboundaries, 
$\Psi_k' - \Psi_k$ is a coboundary. 
Hence its cohomology class in $H^k_c({\cal V}(M), \Omega^k(M,\R))$ 
does not depend on $\nabla$. 
\end{definition}

\begin{remark} \label{rem:psi-int} The Lie algebra $1$-cocycle $\zeta$ is the differential 
of the group cocycle 
$$ \zeta_D \: \Diff(M) \to \Omega^1(M,\End(TM)), 
\quad \phi \mapsto \phi.\nabla- \nabla, $$
where 
$$ (\phi.\nabla)_X Y := \phi.(\nabla_{\phi^{-1}.X} \phi^{-1}.Y) \quad \mbox{for} \quad 
\phi.X = T(\phi) \circ X \circ \phi^{-1}. $$
Now the $k$-fold cup product of $\zeta_D$ with itself defines  
an $\Omega^k(M,\R)$-valued smooth $k$-cocycle on $\Diff(M)$ 
whose corresponding Lie algebra cocycle is 
$\Psi_k$ (cf.\ \cite{Ne04a}, Lemma~F.3). 
\end{remark}

\begin{remark} \label{rem:4.6} 
Assume that $M$ is a parallelizable $d$-dimensional 
manifold and that the $1$-form $\tau \in \Omega^1(M,\R^d)$ 
implements this trivialization in the sense that each map $\tau_m$ 
is a linear isomorphism $T_m(M) \to\R^d$. 
Then, for each $X \in {\cal V}(M)$, ${\cal L}_X\tau$ 
can be written as 
${\cal L}_X\tau = -\theta(X)\cdot \tau$ for some smooth function 
$\theta(X) \in C^\infty(M,\gl_d(\R))$ and $\theta$ is a 
{\it crossed homomorphism}, i.e., 
$$ \theta([X,Y]) = X.\theta(Y) - Y.\theta(X) + [\theta(X),\theta(Y)] $$
(\cite{BiNe07}, Example~II.3). 
Moreover, $\nabla_X Y := \tau^{-1}(X.\tau(Y))$ 
defines an affine connection on $M$ for which the parallel vector fields 
correspond to constant functions, and for the corresponding 
cocycle $\zeta$, the map 
$$ \tilde\tau \: \Gamma(\End(TM)) \to C^\infty(M,\gl_d(\R)), 
\quad \tilde\tau(\phi)(m) = 
\tau_m \circ \phi_m \circ \tau_m^{-1} $$
satisfies 
$\tilde\tau \circ \zeta(X) = - \dd\big(\theta(X)\big) \in \Omega^1(M,\gl_d(\R)).$
This leads to 
\begin{equation}
  \label{eq:psi_k}
\Psi_k (X_1, \ldots, X_k) = \sum_{\sigma\in S_k} \sgn(\sigma)
\Tr \left( \dd\theta(X_{\sigma(1)}) \wedge \ldots \wedge \dd\theta(X_{\sigma(k)}) \right). 
\end{equation}
If $TM$ is trivial, we further obtain the cocycles 
$\Phi_k \in Z_c^{2k-1} ({\cal V}(M), C^\infty(M,\R)),$
$$ \Phi_k (X_1, \ldots, X_{2k-1}) = \sum_{\sigma\in S_{2k-1}} \sgn(\sigma)
\Tr \left( \theta(X_{\sigma(1)}) \cdots \theta(X_{\sigma(2k-1)}) \right)$$
and $\Bp_k \in Z^k_c ({\cal V}(M), \oline\Omega^{k-1}(M,\R))$, defined by 
$$ \Bp_k (X_1, \ldots, X_k) = \sum_{\sigma\in S_k} \sgn(\sigma)
[\Tr \left( \theta(X_{\sigma(1)}) \dd\theta(X_{\sigma(2)}) \wedge \ldots \wedge \dd\theta(X_{\sigma(k)}) \right)]$$ 
(\cite{BiNe07}, Def.~II.7). For each $k\geq 1$ we have 
$\dd \circ \oline\Psi_k = \Psi_k.$

For $k = 1$ we have in particular
$$ \oline\Psi_1(X) = \Tr(\theta(X)) \in C^\infty(M,\R) = \oline\Omega^0(M,\R), 
\qquad
\Psi_1(X) = \Tr(\dd \theta(X)), $$
and 
$$ \oline\Psi_2(X_1, X_2) = [\Tr\big(\theta(X_1)\dd\theta(X_2)
-\theta(X_2)\dd\theta(X_1)\big)] \in \oline\Omega^1(M,\R). $$
\end{remark}

Collecting the results from \cite{BiNe07}, Section~IV, we now obtain  the 
following classification result which, for a parallelizable manifold $M$, 
provides in particular all the information on the twists of the fundamental  
central extensions of the mapping Lie algebras defined by cocycles of 
type (I)-(III). 

\begin{theorem} \label{thm:h2vect} Let 
$M$ be a compact $d$-dimensional manifold with 
trivial tangent bundle. Then the following assertions hold: 
\begin{description}
\item[\rm(1)] For $p \geq 2$ and $d \geq 2$ the map 
$$H^{p+2}_{\rm dR}(M,\R) \to H^2({\cal V}(M),\oline\Omega^p(M,\R)), \quad 
[\omega] \mapsto [\omega^{[2]}] $$
is a linear isomorphism. 
\item[\rm(2)] For $d \geq 2$ we have 
$$ H^2({\cal V}(M),\oline\Omega^1(M,\R)) \cong H^3_{\rm dR}(M,\R)
 \oplus \R [\oline\Psi_1 \wedge \Psi_1] 
\oplus \R[\oline\Psi_2], $$
where $H^3_{\rm dR}(M,\R)$ embeds via $[\omega] \mapsto [\omega^{[2]}]$. 
\item[\rm(3)] For $d \geq 2$ the map 
$$ H^2_{\rm dR}(M,\R) \oplus H^1_{\rm dR}(M,\R) \to 
H^2({\cal V}(M), C^\infty(M,\R)), \quad 
([\alpha], [\beta]) \mapsto 
[\alpha + \beta \wedge \oline\Psi_1] $$
is a linear isomorphism. 
\item[\rm(4)] For $M = \bS^1$, 
$\dim H^2({\cal V}(M), C^\infty(M,\R)) = 2$ and 
$\dim H^2({\cal V}(M), \oline\Omega^1(M,\R)) = \dim H^2({\cal V}(M), \R) 
= 1$. 
\end{description}
\end{theorem}
In addition, for any compact connected manifold $M$,  
$$ H^2({\cal V}(M), \Omega^2(M,\R)) = \R [\Psi_2] \oplus \R [\Psi_1 \wedge \Psi_1] $$
is $2$-dimensional and 
$$ H^2({\cal V}(M), \Omega^1(M,\R)) = \R [\oline\Psi_1 \wedge \Psi_1] \oplus 
\{ [\alpha \wedge \Psi_1] \: 
[\alpha] \in H^1_{\rm dR}(M,\R) \}, $$
is of dimension $1 + b_1(M)$.

\begin{remark}
The main point of the preceding theorem is the precise description of the 
cohomology spaces, which is obtained under the assumption that 
the tangent bundle is trivial. The cocycles $\Psi_k$ exist for all manifolds, 
whereas the construction of the cocycles $\Phi_k$ and 
$\oline\Psi_k$ requires at least 
the existence of a flat affine connection on $M$: 
If $\nabla$ is a flat affine connection on the $d$-dimensional manifold 
$M$, then we have a holonomy homomorphism $\rho \: \pi_1(M) \to \GL_d(\R)$, 
for which the tangent bundle $TM$ is equivalent to the associated bundle 
$\tilde M \times_\rho \R^d$. 
We then obtain a trivializing 
$1$-form $\tilde\tau \in \Omega^1(\tilde M,\R^d)$ satisfying 
$\gamma^*\tilde\tau = \rho(\gamma)^{-1} \circ \tilde\tau$ for $\gamma \in \pi_1(M)$
and a corresponding crossed homomorphism 
$$ \tilde\theta \: {\cal V}(\tilde M) \to C^\infty(\tilde M,\gl_d(\R)). $$
For $X \in {\cal V}(M)$, let $\tilde X \in {\cal V}(\tilde M)^{\pi_1(M)}$ denote 
the canonical lift of $X$ to $\tilde M$. Then 
$$ \gamma^*\tilde\theta(\tilde X) = \rho(\gamma)^{-1} \circ \tilde\theta(\tilde X) 
\circ \rho(\gamma) 
\quad \mbox{ for } \quad \gamma \in \pi_1(M), $$
which means that the cocycles $\Psi_k$, $\Phi_k$ and $\oline\Psi_k$ 
lead on vector fields of the form $\tilde X$ to $\pi_1(M)$-invariant 
differential forms, hence that they actually define cocycles on 
${\cal V}(M)$ with values in 
$\Omega^k(M,\R)$, $C^\infty(M,\R)$ and $\oline\Omega^{k-1}(M,\R)$. 
\end{remark}

\begin{problem} Compute the second cohomology of ${\cal V}(M)$ in 
$\oline\Omega^1(M,\R)$ for all (compact) connected smooth manifolds. 
The cohomology with values in $\Omega^p(M,\R)$ can be reduced with 
results of Tsujishita (\cite{Tsu81}, Thm.~5.1.6; \cite{BiNe07}, Thm.~III.1) 
to the description of 
$H^\bullet({\cal V}(M), C^\infty(M,\R))$ which depends very much on 
the topology of $M$ (\cite{Hae76}). 
\end{problem}

\subsection{Abelian extensions of diffeomorphism groups} 

We now address the integrability of the cocycles described in the preceding 
subsection. Here a key method consists in a systematic use of crossed homomorphisms 
which are used to pull back central extensions of mapping groups (cf.\cite{Bi03}). 

\begin{definition} Let $G$ and $N$ be Lie groups and 
$S \: G \to \Aut(N)$ a smooth action of $G$ on $N$. A smooth map 
$\theta \: G \to N$ is called a {\it crossed homomorphism} if 
$$ \theta(g_1g_2) = \theta(g_1) \cdot g_1.\theta(g_2), $$
i.e., if $\tilde\theta := (\theta,\id_G) \: G\to N \rtimes_S G$ is a morphism of 
Lie groups. 
\end{definition}

Crossed homomorphisms are non-abelian generalizations 
of $1$-cocycles. They have the interesting application that for any smooth 
$G$-module $V$ (considered as a module of $N \rtimes_S G$ on which $N$ 
acts trivially), we have a natural pull-back map 
\begin{equation}
  \label{eq: crossed-pull}
\tilde\theta^* = (\theta,\id_G)^* \: H^2(N \rtimes_S G,V) \to H^2(G,V). 
\end{equation}

Two crossed homomorphisms $\theta_1, \theta_2 \: G \to N$ are said to 
be {\it equivalent} if there exists an $n \in N$ with 
$\theta_2 = c_n \circ \theta_1$, $c_n(x) = nxn^{-1}$. 
Since $N$ acts trivially by conjugation 
on the cohomology group $H^2(N \rtimes_S G,V)$ (cf.\ \cite{Ne04a}, Prop.~D.6), 
equivalent crossed homomorphisms define the same pull-back map on cohomology.

\begin{remark} (a) If $\theta \: G \to N$ is a crossed homomorphism, then 
$\tilde\theta = (\theta, \id_G)$ yields a new splitting on the semidirect 
product group $N \rtimes_S G$, which leads to an isomorphism 
$N \rtimes_S G \cong N \rtimes_{S'} G$, where 
$S'(g) = c_{\theta(g)} \circ S(g)$. 

As we shall see below, it may very well happen that a Lie group extension 
of $N \rtimes_S G$ is trivial on $\{\1\} \times G$ but not on $\tilde\theta(G)$, 
which leads to a non-trivial extension of $G$ (Example~\ref{ex:vira}). 

(b) If $\theta \: G \to N$ is a crossed homomorphism of Lie groups, 
$\L(\theta) := T_\1(\theta) \: \L(G) \to \L(N)$ is a crossed homomorphism 
of Lie algebras because $(\L(\theta),\id_{\L(G)}) = \L(\tilde\theta) \: 
\L(G) \to \L(N) \rtimes \L(G)$ is a morphism of Lie algebras. 

(c) If the central extension $\hat\fn = V \oplus_\omega \fn$ is defined by a 
$\g$-invariant cocycle $\omega$ and $\tilde\omega$ is the corresponding 
cocycle of the abelian extension $\hat\fn \rtimes \g$, corresponding to the 
diagonal action of $\g$ on $\hat\fn$, then the pull-back of this cocycle 
with respect to a crossed homomorphism $\theta \: \g \to \fn$ has the simple form 
$$ (\tilde\theta^*\tilde\omega)((\theta(x),x),(\theta(y),y)) 
= \omega(\theta(x),\theta(y)) = (\theta^*\omega)(x,y). $$
\end{remark}

As we shall see below, crossed 
homomorphisms can be used to prove the integrability of many 
interesting cocycles of the Lie algebra of vector fields 
to the corresponding group of diffeomorphisms. 

\begin{example} (a) Let $q \: \bV \to M$ be a natural vector bundle over $M$ with typical fiber $V$. Then 
the action of $\Diff(M)$ on $M$ lifts to an action on $\bV$ by bundle 
automorphisms. Any trivialization of $\bV$ yields an 
isomorphism 
$$\Aut(\bV) \cong \Gau(\bV) \rtimes \Diff(M) 
\to \Aut(M \times V) \cong C^\infty(M,\GL(V)) \rtimes \Diff(M), $$
and the lift $\rho \: \Diff(M) \to \Aut(\bV)$ can now be written 
as $\rho = (\theta, \id_{\Diff(M)})$ for some crossed homomorphism 
$\theta \: \Diff(M) \to C^\infty(M,\GL(V))$.

(b) Let $\tau \in \Omega^1(M,\R^d)$ be such that each $\tau_m$ 
is invertible in each $m \in M$, so that 
$\tau$  defines a trivialization of the tangent bundle of $M$. 
In \cite{BiNe07}, Ex.~II.3, it is shown that 
${\cal L}_X\tau = -\theta(X)\cdot \tau$ defines a crossed homomorphism 
$$ \theta \: {\cal V}(M) \to C^\infty(M,\gl_d(\R))$$ 
(cf.\ Remark~\ref{rem:4.6}). 
A corresponding crossed homomorphism on the group level is given by 
$$ \Theta \: \Diff(M) \to C^\infty(M,\GL_d(\R)), \quad 
\phi.\tau = (\phi^{-1})^*\tau = \Theta(\phi)^{-1} \cdot \tau. $$
It satisfies $\L(\Theta) = \theta$, taking into account that 
${\cal V}(M)$ is the Lie algebra of the opposite group $\Diff(M)^{\rm op}$, 
which causes a minus sign. 

We shall use this crossed homomorphism to write the 
cocycles $\oline\Psi_2$ and $\oline\Psi_1 \wedge \Psi_1$ as pull-backs 
(cf.\ \cite{Bi03} for the case where $M$ is a torus). 
For the Lie algebra $\fk := \gl_d(\R)$, the space $\Sym^2(\fk,\R)^\fk$ 
is $2$-dimensional, spanned by the two forms 
$$ \kappa_1(x,y) := \tr(x)\tr(y) \quad \hbox{ and } \quad 
\kappa_2(x,y) := \tr(xy). $$
Then the relations 
\begin{align*}
[\tr(\theta(X_1))\tr(\dd \theta(X_2))-\tr(\theta(X_2))
\tr(\dd\theta(X_1))] 
&= (\oline\Psi_1 \wedge \Psi_1)(X_1, X_2), \\
[\tr(\theta(X_1)\dd \theta(X_2)-\theta(X_2)\dd\theta(X_1))] 
&= \oline\Psi_2(X_1, X_2) 
\end{align*}
imply that 
$$ 2\theta^*\omega_{\kappa_1} = \oline\Psi_1 \wedge \Psi_1
\quad \mbox{ and } \quad 
2\theta^*\omega_{\kappa_2} =\oline\Psi_2. $$
\end{example} 

With these preparations, we are now ready to prove the following: 

\begin{proposition} If the tangent bundle $TM$ is trivial, then the cocycles 
$$\oline\Psi_1 \wedge \Psi_1\quad \mbox{ and}\quad 
\oline\Psi_2 \quad \mbox{ in } \quad  Z^2({\cal V}(M), \oline\Omega_1(M,\R))$$ 
integrate to abelian Lie group extensions of some covering group of 
$\Diff(M)_0$. 
\end{proposition}

\begin{proof} Fix $i \in \{1,2\}$. Since the invariant form $\kappa_i$ 
is one of the two components of the universal invariant 
symmetric bilinear form of $\gl_d(\R)$, 
Theorem~\ref{thm:int-type1} implies that the cocycles 
$2\omega_{\kappa_i}$ of $C^\infty(M,\gl_d(\R))$ are integrable 
to a central extensions $\hat G$ of the simply connected covering group 
$\tilde G$ of $G := C^\infty(M,\GL_d(\R))_0$ by 
$Z := \oline\Omega^1(M,V)/\im(\per_{\omega_{\kappa_i}})$ 
on which $\Diff(M)$ acts smoothly by automorphisms 
(Remark~\ref{rem:3.3} and the Lifting Theorem~\ref{thm:lifting-theorem}). 
We thus obtain a semidirect product group 
$\hat H := \hat G \rtimes \tilde\Diff(M)_0$ which is an abelian extension 
of $H := \tilde G \rtimes \tilde\Diff(M)_0$ by $Z$. 

The crossed homomorphism 
$\Theta \: \Diff(M)_0 \to G$ has a unique lift to a crossed homomorphism 
$\Xi \: \tilde\Diff(M)_0 \to \tilde G$ and the pull back 
extension $\tilde\Xi^*\hat H$ integrates $2\theta^*\omega_{\kappa_i}$. 
In view of the preceding example, the assertion follows. 
\smartqed\qed
\end{proof} 

\begin{problem} It would be nice to know if we really need to pass to a 
covering group in the preceding proposition. As the argument in the proof 
shows, this is not necessary if the crossed homomorphism $\Theta$ 
lifts to a crossed homomorphism into $\tilde G$. The obstruction for that 
is the homomorphism 
$$ \pi_1(\Theta) \: \pi_1(\Diff(M)) \to \pi_1(G), $$
but this homomorphism does not seem to be very accessible. 

In all cases where we have a compact subgroup 
$C \subeq \Diff(M)$ for which $\pi_1(C) \to \pi_1(\Diff(M))$ is surjective 
(this is in particular the case for $\dim M \leq~2$; \cite{EE69}), 
we may choose the trivialization $\tau$ in a $C$-invariant fashion, so that 
$\Theta(C)$ is trivial. This clearly implies that $\pi_1(\Theta)$ vanishes. 
\end{problem}

We now take a closer look on some examples of extensions of 
$\Diff(M)_0$ defined by these the two cocycles discussed above. 

\begin{example} \label{ex:vira} (cf.\ \cite{Bi03}) (The Virasoro group) For $M = \bS^1$, 
the preceding constructions can be used in 
particular to obtain the Virasoro group, resp., a corresponding 
global smooth $2$-cocycle. 

The simply connected cover of $C^\infty(\bS^1,\GL_1(\R))_0$ is the 
group $N := C^\infty(\bS^1, \R)$ on which we have a smooth 
action of $G := \Diff(\bS^1)_0^{\rm op}$ by $\phi.f := f \circ \phi$. 
In the notation of the preceding example, we have $\kappa := \kappa_1 = \kappa_2$, and, 
identifying ${\cal V}(\bS^1)$ with $C^\infty(\bS^1, \R)$, we obtain 
the $1$-cocycle $\theta(X) = X'$. The corresponding $1$-cocycle 
on $G$ with values in $N$ is given by 
$\Theta(\phi) = \log(\phi')$. 

The cocycle $\omega := 2\omega_\kappa$ on ${\cal V}(\bS^1)$ with values in 
$\oline\Omega^1(\bS^1,\R) \cong \R$ is 
$$ \omega(\xi_1, \xi_2) = \int_0^1 \xi_1 \xi_2'- \xi_2\xi_1'\, dt 
= 2\int_{\bS^1} \xi_1 \dd \xi_2,$$ 
and $\frac{1}{2}\omega = \omega_\kappa$ is a corresponding 
group cocycle, defining the central extension 
$\hat N := \R \times_\omega N$ on which $G$ acts by 
$\phi.(z, f) := (z, \phi.f)$, so that 
the semi-direct product $\hat N \rtimes G$ is a central 
extension of $N \rtimes G$ (cf.\ \cite{Ne02}).

We now pull this extension back with $\tilde\Theta \: G \to N \rtimes G$. 
Since we have in $\hat N \rtimes G$ the relation 
$$ (0,\Theta(\phi),\phi)(0,\Theta(\psi),\psi)
= \Big(\omega_\kappa(\Theta(\phi),\phi.\Theta(\psi)),
\Theta(\phi)\phi.\Theta(\psi),\phi\psi\Big),$$
this leads to the $2$-cocycle  
\begin{align*}
&\Omega(\phi, \psi) 
= \omega_\kappa(\Theta(\phi), \phi.\Theta(\psi))
= \omega_\kappa(\Theta(\phi), \Theta(\phi\psi) - \Theta(\phi)) \\
&= \omega_\kappa(\Theta(\phi), \Theta(\phi\psi)) 
= -\omega_\kappa(\log(\psi \circ \phi)', \log\phi') 
= - \int_{\bS^1} \log(\psi \circ \phi)' \dd(\log\phi'). 
\end{align*}
This is the famous Bott--Thurston cocycle for $\Diff(\bS^1)_0^{\rm op}$, 
and the corresponding central extension is the Virasoro group. 
On the Lie algebra level, the pull-back cocycle is 
$$ (\theta^*\omega)(\xi_1,\xi_2) = \omega(\xi_1', \xi_2') 
= \int_0^1 \xi_1'\xi_2''- \xi_2'\xi_1''\, dt $$ 
(cf.\ \cite{FF01}). 
\end{example}

\begin{example} (\cite{Bi03}) The preceding example easily generalizes to tori 
$M = \T^d$. Then we can trivialize $TM$ with the Maurer--Cartan form $\tau$ 
of the group $\T^d$, which leads after identification of ${\cal V}(M)$ with 
$C^\infty(M,\R^d)$ to $\theta(X) = - X'$ (the Jacobian matrix of $X$). 
A corresponding 
crossed homomorphism is 
$$ \Theta \: \Diff(\T^d) \to C^\infty(\T^d,\GL_d(\R)), \quad 
\Theta(\phi) := (\phi^{-1})', $$
where the derivative $\phi'$ of an element of 
$\Diff(\T^d) \subeq C^\infty(\T^d,\T^d)$ is considered as a $\GL_d(\R)$-valued 
smooth function on $\T^d$. 

The main difference to the one-dimensional case is that here 
we do not have an explicit formula for the action of $\Diff(\T^d)$ 
on the central extension of the mapping group $C^\infty(M,\GL_d(\R))_0$. 
As a consequence, we do not get explicit cocycles integrating 
$\oline\Psi_2$ and $\oline\Psi_1 \wedge \Psi_1$. 
\end{example}

\subsection*{On the integration of cocycles of the form $\omega^{[2]}$} 

Now we turn to the other types of cocycles listed in 
Theorem~\ref{thm:h2vect}(2), namely those of the form $\omega^{[2]}$, where 
$\omega$ is a closed $V$-valued $3$-form on $M$. 

Let $G$ be a Lie group and  $\sigma \: M \times G \to M$ be a smooth right action 
on the (possibly infinite dimensional) smooth manifold 
$M$. Let $\omega \in \Omega^{p+q}(M,V)$ be a closed $(p+q)$-form 
with values in a Mackey complete space $V$. 
If \break $\L(\sigma) \: \g \to {\cal V}(M)$ is the corresponding homomorphism of 
Lie algebras, then 
$$ \omega_\g := \L(\sigma)^*\omega^{[p]} \in Z^p(\g,\oline\Omega^q(M,V)) $$
is a Lie algebra cocycle and we have a well-defined 
period homomorphism 
$$ \per_{\omega_\g} \: \pi_p(G) \to \oline\Omega^q(M,V)^G. $$
If $M$ is smoothly paracompact, which is in particular the case if 
$M$ is finite-dimensional, then any class $[\alpha] \in 
\oline\Omega^q(M,V)$ is determined by its integrals over smooth 
singular $q$-cycles $S$ in $M$ which permits us to determine 
$\per_{\omega_\g}$ geometrically. 

In the following we write for an oriented compact manifold $F$: 
$$ \fint_F \: \Omega^{p+q}(M \times F,V) \to \Omega^q(M,V) $$
for the fiber integral (cf.\ \cite{GHV72}, Ch.~VII). 

\begin{theorem}
  \label{thm:per-thm} {\rm(Period Formula)} For any smooth singular $q$-cycle $S$ on $M$ 
and any smooth map $\gamma \: \bS^p \to G$, we have 
the following integral formula 
\begin{equation}
  \label{eq:per-form}
\int_S \per_{\omega_\g}([\gamma]) = (-1)^{pq} \int_{S \bullet \gamma} \omega,
\end{equation}
where $S \bullet \gamma \in H_{p+q}(M)$ denotes the singular cycle 
obtained from the natural map 
$H_q(M) \otimes H_p(G) \to H_{p+q}(M), 
[\alpha] \otimes [\beta] \mapsto [\sigma \circ (\alpha \times \beta)]$
induced by the action map $M \times G \to M$. 
\end{theorem}

\begin{proof} First, let $\gamma \: \bS^p \to G$ be a smooth map and 
$\alpha \: \Delta_q \to M$ be a smooth singular simplex. We then obtain 
a smooth map 
$$ \alpha \bullet \gamma :=  \sigma \circ (\alpha \times \gamma) \: 
\Delta_q \times \bS^p \to M, \quad 
(x,y) \mapsto \alpha(x).\gamma(y). $$
Integration of $\omega$ yields with 7.12/14 in \cite{GHV72} (on pull-backs and 
Fubini's Theorem):  
\begin{align*} 
&\int_{\alpha \bullet\gamma} \omega 
= \int_{\Delta_q \times \bS^p} (\alpha \bullet\gamma)^*\omega 
= \int_{\Delta_q \times \bS^p} (\alpha \times\gamma)^*\sigma^*\omega \cr
&= \int_{\Delta_q \times \bS^p} (\alpha \times\id_{\bS^p})^*(\id_M \times \gamma)^*
\sigma^*\omega 
= \int_{\Delta_q} \fint_{\bS^p}(\alpha \times\id_{\bS^p})^*(\id_M \times \gamma)^*
\sigma^*\omega \cr
&= \int_{\Delta_q} \alpha^*\fint_{\bS^p}(\id_M \times \gamma)^*
\sigma^*\omega 
= \int_{\alpha} \underbrace{\fint_{\bS^p}(\id_M \times \gamma)^*\sigma^*\omega}
_{\in \Omega^q(M,V)}. 
\end{align*} 
It therefore remains to show that 
\begin{equation}
  \label{eq:per-ident}
\per_{\omega_\g}([\gamma]) = (-1)^{pq} 
\Big[\fint_{\bS^p}((\id_M \times \gamma)^*\sigma^*\omega)\Big]. 
\end{equation}
So let $Y_1, \ldots, Y_q \in {\cal V}(M)$ be smooth vector fields on $M$ 
and write \break $j_m^F \: F \to M \times F, x \mapsto (m,x)$ for the 
inclusion map. We then have 
\begin{align*}
&\fint_{\bS^p}((\id_M \times \gamma)^*\sigma^*\omega)(Y_1,\ldots, Y_q)_m  
= \int_{\bS^p} (j_m^{\bS^p})^*(i_{Y_q} \cdots i_{Y_1}
((\id_M \times \gamma)^*\sigma^*\omega)) \cr
&= \int_{\bS^p} (j_m^{\bS^p})^*(\id_M \times \gamma)^* (i_{Y_q} 
\cdots i_{Y_1}(\sigma^*\omega)) 
= \int_\gamma (j_m^G)^*(i_{Y_q} \cdots i_{Y_1}(\sigma^*\omega)). 
\end{align*}
We further observe that for $x_1,\ldots, x_p \in \g = T_\1(G)$ we have 
\begin{align*}
&\ \ \ \ (i_{Y_q} \cdots i_{Y_1}(\sigma^*\omega))_{(m,g)}(g.x_1, \ldots, g.x_p)\cr 
&= (\sigma^*\omega)_{(m,g)}(Y_1(m),\ldots, Y_q(m),g.x_1, \ldots, g.x_p) \cr
&= \omega_{(m.g)}(Y_1(m).g,\ldots, Y_q(m).g, m.(g.x_1), \ldots, m.(g.x_p)) \cr
&= \omega_{(m.g)}(Y_1(m).g,\ldots, Y_q(m).g, \L(\sigma)(x_1)(m.g), \cdots, 
\L(\sigma)(x_p)(m.g)) \cr 
&= (-1)^{pq} (\sigma_g^*(i_{\L(\sigma)(x_p)}\cdots i_{\L(\sigma)(x_1)}\omega))
(Y_1(m),\ldots, Y_q(m)) \cr 
&= (-1)^{pq} (g.\omega_\g(x_1,\ldots, x_p))(Y_1,\ldots, Y_q)_m \cr 
&= (-1)^{pq} \omega_\g^{\rm eq}(g.x_1, \ldots, g.x_p)
(Y_1,\ldots, Y_q)_m. 
\end{align*}
We thus obtain in $\Omega^q(M,V)$ the identity: 
$$ \fint_{\bS^p}(\id_M \times \gamma)^*\sigma^*\omega  
= (-1)^{pq} \int_\gamma \omega_\g^{\rm eq}
= (-1)^{pq} \per_{\omega_\g}([\gamma]), $$
which in turn implies (\ref{eq:per-ident}) and hence completes the proof. 
\smartqed\qed
\end{proof}

\begin{example} \label{ex:3.19} For $q =0$ and a closed $p$-form $\omega \in 
\Omega^p(M,V)$, the preceding theorem implies in particular that 
$\per_{\omega_\g} \: \pi_p(G) \to C^\infty(M,V)$ satisfies 
$$ \per_{\omega_\g}([\gamma])(m) := \int_{m \bullet \gamma} \omega 
\quad \mbox{ with } \quad (m \bullet \gamma)(x) = m.\gamma(x). $$
\end{example}

As an important consequence of the period formula in the preceding theorem, 
we derive that the $2$-cocycles $\omega^{[2]}$ of integral $(p+2)$-forms  
integrate to abelian extensions of a covering group of $\Diff(M)_0$. 
This applies in particular to $p = 1$, which leads to the 
integrability of the cocycles occurring in Theorem~\ref{thm:h2vect}(2). 

\begin{corollary} Let $\omega \in \Omega^{p+2}(M,V)$ be closed 
with discrete period group $\Gamma_\omega$. 
Then there exists an abelian extension of the 
simply connected covering group $\tilde\Diff(M)_0$ by 
the group $Z:= \oline\Omega(M,V)/\Gamma_Z$, where 
$$ \Gamma_Z := \Big\{ [\alpha] \in H^p_{\rm dR}(M,V)\: 
\int_{Z_p(M)} \alpha \subeq \Gamma_\omega\Big\} 
\cong \Hom(H_p(M), \Gamma_\omega), $$ 
whose Lie algebra 
cocycle is $\omega^{[2]} \in Z^2({\cal V}(M),\oline\Omega^p(M,V))$. 
\end{corollary}

\begin{proof} Applying the period formula in Theorem~\ref{thm:per-thm} 
to $\omega_\g = \omega^{[2]}$ and $G = \Diff(M)^{\rm op}$, which 
acts from the right on $M$, we see that the image of the period map 
is contained in the group $\Gamma_Z$, which is discrete. 
Now Theorem~\ref{thm:abext} below applies. 
\smartqed\qed
\end{proof} 

It is a clear disadvantage of the preceding theorem that it says only something 
about the simply connected covering group, and since $\pi_1(\Diff(M))$ is not 
very accessible, it would be nicer if we could say more about the flux homomorphism 
of $\omega^{[2]}$. Unfortunately it does not vanish in general, so that we 
cannot expect an extension of the group $\Diff(M)_0$ itself. However, 
the following result is quite useful to evaluate the flux: 

\begin{theorem} {\rm(Flux Formula)} Let $\omega \in \Omega^{p+2}(M,V)$ 
be a closed $(p+2)$-form and 
$\gamma \: \bS^1 \to \Diff(M)$ a smooth loop. Then 
$\eta_\gamma := \int_0^1 \gamma(t)^*\big(i_{\delta^r(\gamma)_t}\omega\big)\, dt$
is a closed $(p+1)$-form and the flux of $\omega^{[2]}$ can be calculated as 
$$F_{\omega^{[2]}}([\gamma]) = [\eta_\gamma^{[1]}] \in H^1({\cal V}(M), 
\oline\Omega^p(M,V)). $$
\end{theorem}

\begin{proof} We recall from Appendix A 
that $F_{\omega^{[2]}}([\gamma]) = [I_\gamma]
= - [I_{\gamma^{-1}}]$ for 
\begin{align*}
I_\gamma(X) 
&=  \int_0^1 \gamma(t).\omega^{[2]}(\Ad(\gamma(t))^{-1}.X, \delta^l(\gamma)_t)\, dt.
\end{align*}
Hence $I_{\gamma^{-1}}(X) \in \oline\Omega^p(M,V)$ is represented 
by the $p$-form $\tilde I_{\gamma^{-1}}(X)$, 
defined on vector fields $Y_1,\ldots, Y_p$ by 
\begin{align*}
&\tilde I_{\gamma^{-1}}(X)(Y_1, \ldots, Y_p)  \\
&=  \int_0^1 \gamma(t)^*\big(\omega(\Ad(\gamma(t)).X, \delta^l(\gamma^{-1})_t, 
\Ad(\gamma(t)).Y_1, \ldots, \Ad(\gamma(t)).Y_p)\, dt \\
&=  \int_0^1 \gamma(t)^*\Big(\big(i_{\delta^r(\gamma)_t}\omega\big)
(\Ad(\gamma(t)).X, \Ad(\gamma(t)).Y_1, \ldots, \Ad(\gamma(t)).Y_p)\Big)\, dt \\
&=  -\int_0^1 \Big(\gamma(t)^*\big(i_{\delta^r(\gamma)_t}\omega\big)\Big)
(X, Y_1, \ldots, Y_p)\, dt \\
&=  -\eta_\gamma(X,Y_1, \ldots, Y_p)
=  -\eta_\gamma^{[1]}(X)(Y_1, \ldots, Y_p).
\end{align*}
That $\eta_\gamma$ is closed follows directly from the fact that 
for any smooth path $\gamma \: I \to \Diff(M)$, we have 
\begin{align*}
\gamma(1)^*\omega - \omega 
&= \int_0^1  \gamma(t)^*\big({\cal L}_{\delta^r(\gamma)_t}\omega\big)\, dt 
= \int_0^1  \gamma(t)^*\big(\dd i_{\delta^r(\gamma)_t}\omega\big)\, dt \\
&= \dd\int_0^1  \gamma(t)^*\big(i_{\delta^r(\gamma)_t}\omega\big) \, dt 
= \dd \eta_\gamma. 
\end{align*}
\smartqed\qed
\end{proof}

In general it is not so easy to get good hold of the flux homomorphism, but 
there are cases where it factors through an evaluation homomorphism 
$\pi_1(\ev_m) \: \pi_1(\Diff(M)) \to \pi_1(M,m)$. According to Lemma~11.1 in 
\cite{Ne04a}, this is the case if $\omega$ is a volume form on $M$.

At this point we have discussed the integrability of all 
the $2$-cocycles of ${\cal V}(M)$ that are relevant to understand the 
twists of the fundamental central extensions of the mapping algebras 
with values in $\oline\Omega^1(M,V)$. 

The twists relevant for cocycles with values in $\Omega^1(M,V)$ and 
$C^\infty(M,V)$ are easier to handle. 
First we recall that 
a product of two integrable Lie algebra $1$-cocycles is 
integrable to the cup product of the corresponding group cocycles 
(\cite{Ne04a}, Lemma~F.3), so that, combined with the 
integrability of the $\Psi_k$ (Remark~\ref{rem:psi-int}), 
the relevant information in contained in the following remark. 

\begin{remark} (a) If $\alpha \in \Omega^1(M,\R)$ is a closed $1$-form, 
then Example~\ref{ex:3.19}, applied to $p =1$ provides a condition 
for integrability to a $1$-cocycle on $\Diff(M)$. Passing to a covering 
$q \: \hat M\to M$ of $M$ for which $q^*\alpha$ is exact shows that 
there exists a covering of the whole diffeomorphism group to which 
$\alpha$ integrates as a $1$-cocycle with values in $C^\infty(M,\R)$ 
(cf.\ Remark~\ref{rem:2.3}(c)). 

(b) If $\omega  \in \Omega^2(M,\R)$ is a closed $2$-form, 
then Example~\ref{ex:3.19}, applied with $p = 2$,
 shows that if 
$\omega$ has a discrete period group $\Gamma_\omega$, 
then it integrates to a group cocycle on some covering group 
of $\Diff(M)_0$. This can be made more explicit by 
the observation that in this case $Z := \R/\Gamma_\omega$ 
is a Lie group and there exists a $Z$-bundle $q \: P \to M$ 
with a connection $1$-form $\theta$ satisfying $q^*\omega = \dd\theta$. 
Then the group $\Aut(M)$ is an abelian extension of the open subgroup 
$\Diff(M)_P$ of all diffeomorphisms lifting to bundle automorphisms of 
$P$ by $C^\infty(M,Z)$ and $\omega$ is a corresponding Lie algebra 
cocycle (cf.\ \cite{Ko70}, \cite{NV03}). 

For these cocycles the flux can also be made quite explicit 
(cf.\ \cite{Ne04a}, Prop.~9.11): For 
$\alpha \in C^\infty(\bS^1,M)$ and 
$\gamma \in C^\infty(\bS^1,\Diff(M)_0)$, we have 
$$ \int_\alpha F_\omega([\gamma]) = \int_{\gamma^{-1} \bullet \alpha}\omega, 
\quad 
\mbox{ where } \quad (\gamma^{-1}\bullet \alpha)(t,s) 
:= \gamma(t)^{-1}(\alpha(s)). $$

(c) If $\mu$ is a volume form on $M$, then the cocycle 
$\oline\Psi_1$ satisfies ${\cal L}_X\mu = - \oline\Psi_1(X) \mu$ 
(\cite{BiNe07}, Lemma~III.3), so that $\oline\Psi_1(X) = - \div X$. 
Therefore $\zeta(\phi)$, defined by $\phi.\mu = \zeta(\phi)^{-1}\mu$, 
defines a $C^\infty(M,\R)$-valued group cocycle integrating $\oline\Psi_1$.

If $M$ is not orientable, then there is a $2$-sheeted covering 
$q \: \hat M\to M$ such that $\hat M$ is orientable. 
Accordingly we have a $2$-fold central extension 
$\hat\Diff(M)$ of $\Diff(M)$ by $\Z/2$ acting on $\hat M$ and 
the preceding construction provides a cocycle 
$\zeta \: \hat\Diff(M) \to C^\infty(\tilde M,\R)$. 
If $\sigma \: \hat M \to \hat M$ is the non-trivial deck transformation, 
then $\sigma$ reverses the orientation of $\hat M$ and there exists a 
volume form $\mu$ on $\hat M$ with $\sigma^*\mu = -\mu$. 
Since $\sigma$ commutes with $\hat\Diff(M)$, 
all functions $\zeta(\phi)$ are $\sigma$-invariant, 
hence in $C^\infty(\hat M,\R)^\sigma\cong C^\infty(M,\R)$. 
We also note that $\sigma^*\mu = - \mu$ implies that $\sigma$ is not 
contained in the identity component of $\hat\Diff(M)$, so that 
$\hat\Diff(M)_0 \cong \Diff(M)_0$ and $\oline\Psi_1$ also 
integrates to a cocycle on $\Diff(M)_0$. 
\end{remark}

\subsection*{Toroidal groups and their generalizations} 

For the spaces $\fz = \oline\Omega^1(M,V)$, $C^\infty(M,V)$ and 
$\Omega^1(M,V)$ the cocycles $\omega$ of type (I)-(III) 
on $\g = C^\infty(M,\fk)$ are ${\cal V}(M)$-invariant, 
hence extend trivially to cocycles on $\fh := \g \rtimes {\cal V}(M)$ 
and we have also seen how to classify the $\fz$-cocycles $\tau$ 
on ${\cal V}(M)$ that occur as twists (cf.\ Theorems~\ref{thm:h2vect},
\ref{thm:ex-sequ}). 
We write $\omega_\tau$ for the corresponding twisted 
cocycle on $\g \rtimes {\cal V}(M)$. 

Let $G := C^\infty(M,K)_0$. 
If $\omega$ is integrable to a Lie group extension, then we write 
$\hat G$ for a corresponding extension of $\tilde G$ by 
$Z \cong \fz/\Gamma_Z$. This extension is uniquely determined by 
$[\omega]$, which has the consequence that the natural action of
$\Diff(M)$ on $\hat\g$ integrates to a smooth action on $\hat G$, 
which leads to a semidirect product 
$\hat G \rtimes \Diff(M)$. 

If $p \: \hat\Diff(M) \to \tilde\Diff(M)_0$ is an abelian 
extension of $\tilde\Diff(M)_0$ by $Z$ corresponding to $\tau$, then 
$\hat G \rtimes \hat\Diff(M)$ is an extension of 
$H := G \rtimes \Diff(M)_0$ by $Z \times Z$, and the antidiagonal 
$\oline\Delta_Z$ in $Z$ is a normal split Lie subgroup, so that 
$$ \hat H :=  (\hat G \rtimes \hat\Diff(M))/\oline\Delta_Z $$
carries a natural Lie group structure, and it is now easy to verify 
that $\L(\hat H) \cong \fz \oplus_{\omega_\tau} \fh$. 

For the special case where $M = \T^d$ is a torus, $\fk$ 
is a simple complex algebra, $\kappa$ is the Cartan--Killing form, 
$V = \C$, $\omega = \omega_\kappa$ and $\tau$ is a linear combination
of the two cocycles $\oline\Psi_2$ and 
$\oline \Psi_1 \wedge \Psi_1$, we thus obtain Lie groups 
$\hat H$ whose Lie algebras are Fr\'echet completions of the 
so-called {\it toroidal Lie algebras}. These algebras (with the twists) 
have been introduced by Rao and Moody in \cite{ERMo94} and since then their 
representation theory has been an active research area 
(cf.\ \cite{Bi98}, \cite{Lar99}, \cite{Lar00}, \cite{ERC04}, \cite{FJ07}). 

\section{Central extensions of gauge groups}  \label{sec:gaugeext}

In this section we shall see how, and under which circumstances, 
the fundamental types of cocycles on $C^\infty(M,\fk)$ 
generalize to gauge algebras $\gau(P)$. Therefore it becomes a natural issue 
to understand the corresponding central extensions of the gauge group 
$\Gau(P)$, resp., its identity component. Depending on the complexity 
of the bundle $P$, this leads to much deeper questions, most of which 
are still open. 

\subsection{Central extensions of $\gau(P)$} 
\label{sec:4.1}

In this section we describe natural analogs of the 
cocycles of type (I)-(III) for gauge Lie algebras of non-trivial bundles. 

Let $q \: P \to M$ be a $K$-principal bundle and 
$\gau(P)$ its gauge algebra, which we realize as 
$$ \gau(P) \cong \{ f \in C^\infty(P,\fk) \: 
(\forall p \in P)(\forall k \in K)\ f(p.k) = \Ad(k)^{-1}.f(p)\} $$ 
which is the space of smooth sections 
of the associated Lie algebra bundle $\Ad(P) = P\times_{\Ad} \fk$ 
(cf.\ Remark~\ref{rem:2.3}(b)). 
Further, let $\hat M := P/K_0$ denote the corresponding {\it squeezed bundle}, 
which is a, not necessarily connected, covering of $M$ 
and a $\pi_0(K)$-principal bundle. The squeezed bundle is associated to the 
$\pi_1(M)$-bundle $q_M \: \tilde M \to M$ 
by the homomorphism $\delta_1 \: \pi_1(M) \to \pi_0(K)$ 
from the long exact homotopy sequence of~$P$. Sometimes it is convenient 
to reduce matters to connected structure groups, which amounts to considering  
$P$ as a $K_0$-principal bundle over the covering space $\hat M$ of $M$.  

Let $\rho \: K \to \GL(V)$ be a homomorphism with $K_0 \subeq \ker \rho$, 
so that $V$ is a $\pi_0(K)$-module. Then the associated bundle 
$$ \bV := P \times_\rho V $$
is a flat vector bundle. It is  associated to the 
squeezed bundle $P/K_0$ via the representation $\oline\rho \: \pi_0(K) \to\GL(V)$. 
We  have a natural exterior derivative on the space  
$\Omega^\bullet(M,\bV)$ of $\bV$-valued differential forms and 
we define $\oline\Omega^1(M,\bV) := \Omega^1(M,\bV)/\dd(\Gamma\bV)$. 
If $V$ is finite-dimensional, then $\dd(\Gamma \bV)$ is a closed subspace of 
the Fr\'echet space $\Omega^1(M,\bV)$, so that the quotient inherits 
a natural Hausdorff 
locally convex topology (cf.\ \cite{NeWo07}). 
The Lie algebra $\aut(P) \subeq {\cal V}(P)$ acts in a natural 
way on all spaces of sections of vector bundles associated to $P$, 
via their realization as smooth functions on $P$, hence in 
particular on the spaces 
$\Omega^p(M,\bV)$, which can be realized as $V$-valued 
differential forms on $P$. Since $\bV$ is flat, the action on 
$\Omega^p(M,\bV)$ 
factors through an action of the Lie algebra ${\cal V}(M) \cong 
\aut(P)/\gau(P)$.

\begin{remark} The passage from $M$ to the covering 
$\hat M := P/K_0$, the squeezed bundle associated to $P$, provides 
a simplification of our setting to the case where the structure 
group under consideration is connected if we consider $P$ 
as a $K_0$-bundle over $\hat M$. One has to pay for this reduction 
by changing $M$. In particular, if $\pi_0(K)$ is infinite, the manifold 
$\hat M$ is non-compact. 
\end{remark}

We now take a closer look at analogs of the cocycles of types (I)-(III) for 
gauge algebras $\g := 
\gau(P)$ of non-trivial bundles. Here the natural target spaces 
are $\Gamma\bV, \Omega^1(M,\bV)$ and $\oline\Omega^1(M,\bV)$, replacing 
$C^\infty(M,V), \Omega^1(M,V)$ and $\oline\Omega^1(M,V)$.

{\bf (I):} First we choose a connection $1$-form $\theta \in \Omega^1(P,\fk)$ 
and write $\nabla$ for the corresponding covariant derivative, which induces
covariant derivatives on all associated vector bundles such as $\Ad(P)$. 
For the flat bundle $\bV$ the corresponding covariant derivative coincides 
with the Lie derivative. 
Let \break $\kappa \: \fk \times \fk \to V$ be a $K$-invariant symmetric bilinear 
map. For $\xi \in \gau(P)$, 
the $1$-form $\nabla \xi \in \Omega^1(P,\fk)$ 
is $K$-equivariant and basic, hence defines a bundle-valued $1$-form in 
$\Omega^1(M,\Ad(P))$. Now $\kappa$ defines a $2$-cocycle 
$$ \omega_\kappa \in Z^2(\gau(P),\oline\Omega^1(M,\bV)), \quad 
\omega_\kappa(\xi_1,\xi_2) := [\kappa(\xi_1, \nabla \xi_2)] $$
with values in the trivial module $\oline\Omega^1(M,\bV)$ 
(cf.\ \cite{NeWo07}; and \cite{LMNS98} for bundles with connected $K$). 

{\bf (II):} Let $K$ be a connected Lie group with Lie algebra $\fk$ 
and $\eta \in Z^2(\fk,V)$ be a $2$-cocycle, defining a central 
Lie algebra extension $q_\fk \: \hat\fk \to \fk$ of $\fk$ by $V$ 
which integrates to a central Lie group extension $\hat K$ of $K$ by some 
quotient $Z := V/\Gamma_Z$ of $V$ by a discrete subgroup. Then $K$ acts 
trivially on $V$, so that the bundle $\bV$ is trivial. 

The $K$ action on $\hat\fk$ defines an 
associated Lie algebra bundle $\hat\Ad(P)$ 
with fiber $\hat\fk$ which is a central extension 
of the Lie algebra bundle $\Ad(P)$ by the trivial bundle $M \times V$ 
(cf.\ \cite{Ma05}). Now the Lie algebra $\hat\gau(P) := \Gamma\hat\Ad(P)$ 
of smooth sections of this bundle 
is a central extension of $\gau(P)$ by the space $C^\infty(M,V) \cong \Gamma\bV$.

Any splitting of the short exact sequence of vector bundles
$$ \0 \to M \times V \to \hat\Ad(P) \to \Ad(P) \to \0 $$
leads to a description of the central extension $\hat\gau(P)$ by a 
$2$-cocycle 
$\omega \in Z^2(\gau(P),C^\infty(M,V))$ which is $C^\infty(M,\R)$-bilinear, 
hence can be identified with a section of the associated bundle 
with fiber $Z^2(\fk,V)$.

{\bf (III):} To see the natural analogs of cocycles of type (III), 
 we first note that 
the natural analogs of the cocycles 
$\omega_{0,\eta} = \dd \circ \omega_\eta$ are cocycles of the form
$\dd \circ \omega \in Z^2(\gau(P),\Omega^1(M,V))$, where 
$\omega \in Z^2(\gau(P),C^\infty(M,V))$ corresponds to a central 
extension $\hat\Ad(P)$ of the Lie algebra bundle $\Ad(P)$ 
by the trivial bundle $M \times V$. 

Let $\kappa$ be an invariant $V$-valued symmetric bilinear form $\kappa$ 
and $\theta$ a principal connection $1$-form on $P$ with associated 
covariant derivative $\nabla$. We are looking for a 
$C^\infty(M,\R)$-bilinear alternating map 
$\tilde\eta \in Z^2(\gau(P),C^\infty(M,V))$ for which 
$$ \omega_{\kappa,\nabla,\tilde\eta}(\xi_1, \xi_2) := 
\kappa(\xi_1, \nabla \xi_2)- \kappa(\xi_2, \nabla \xi_1)
- \dd(\tilde\eta(\xi_1, \xi_2)) $$
defines a Lie algebra cocycle with values in $\Omega^1(M,V)$. 

For 
$\tilde\omega_\kappa(\xi_1, \xi_2) 
:= \kappa(\xi_1, \nabla \xi_2)- \kappa(\xi_2, \nabla \xi_1)$
we obtain as in Section~3: 
$$ \dd_\g(\tilde\omega_\kappa)(\xi_1,\xi_2, \xi_3) 
= \dd(\Gamma(\kappa)(\xi_1, \xi_2, \xi_3)). $$
Here we use that in the realization of $\g = \gau(P)$ in $C^\infty(P,\fk)$, we have
$$ \nabla \xi = \dd\xi + [\theta, \xi], $$
so that 
$$\tilde\omega_\kappa(\xi_1, \xi_2) 
= \kappa(\xi_1, \dd\xi_2)- \kappa(\xi_2, \dd\xi_1) 
+ 2 \kappa(\theta, [\xi_2, \xi_1]). $$
Therefore $\dd_\g \omega_{\kappa,\nabla,\tilde\eta} = 0$ 
is equivalent to 
$\Gamma(\kappa) - \dd_\g\tilde\eta$ having values in constant functions, 
and since this map is $C^\infty(M,\R)$-trilinear, we find the condition 
$$\Gamma(\kappa)  = \dd_\g\tilde\eta \quad \mbox{ in } \quad Z^3(\gau(P), C^\infty(M,V)). $$
This means 
that for each $m \in M$, we have fiberwise in local trivializations
$\Gamma(\kappa) = \dd_\fk\tilde\eta(m)$. In particular, $\kappa$ has to 
be exact. 

\begin{lemma} If $\Gamma(\kappa)$ is a coboundary, 
then there exists a $C^\infty(M,\R)$-bilinear alternating map 
$\tilde\eta \in Z^2(\gau(P),C^\infty(M,V))$ for which 
$\omega_{\kappa,\nabla,\tilde\eta}$ is a $2$-cocycle. 
\end{lemma}

\begin{proof} The set $C^2(\fk,V)_\kappa := \{ \eta \in C^2(\fk,V) \: \dd_\fk \eta
= \Gamma(\kappa)\}$ is an affine space on which $K$ acts by affine map. 
Now the associated bundle $A := P \times_K C^2(\fk,V)_\eta$ has a smooth 
section $\tilde\eta$, and any such section defines a 
$C^\infty(M,\R)$-bilinear element of $C^2(\g,C^\infty(M,V))$ with 
$\dd_\g(\dd \circ \tilde\eta) = \dd \circ \Gamma(\kappa) 
= \dd_\g(\tilde\omega_\kappa)$. 
\smartqed\qed
\end{proof}

In view of the preceding lemma, analogs of the type (III) cocycles 
on the Lie algebra $\gau(P)$ with values in $\Omega^1(M,V)$ always exist 
if $\kappa$ is exact.

\begin{problem} (Universal central extensions) 
Suppose that $\fk$ is a semisimple Lie algebra and 
$q \: P \to M$ a principal $K$-bundle. Find a universal central 
extension of $\gau(P)$. If $P$ is trivial, then we have a universal central 
extension by $\oline\Omega^1(M,V(\fk))$, given by the universal 
invariant symmetric bilinear form $\kappa_u$ of $\fk$ 
(Remark~\ref{rem:3.2}(d)). 
The construction described above yields a central extension 
of $\gau(P)$ by $\oline\Omega^1(M,\bV)$ for the associated bundle with 
fiber $V(\fk)$, but it is not clear that this extension is universal. 

A class of gauge algebras closest to those of trivial bundles are 
those of flat bundles which can be trivialized by a finite covering of $M$. 
It seems quite probable that the above central extension is indeed 
universal. The analogous result for multiloop algebras 
(cf.\ Section~\ref{sec:loopalg} below) 
has recently been obtained by E.~Neher (\cite{Neh07}, Thm.~2.13; 
cf.\ also \cite{PPS07}, p.147). 
\end{problem}

\subsection{Covariance of the Lie algebra cocycles} 
\label{sec:4.2}

\begin{remark} \label{rem:cov-gau} (Covariance for type (I)) 
(a) Realizing $\g = \gau(P)$ in $C^\infty(P,\fk)$, 
we have $\nabla \xi = \dd \xi + [\theta,\xi]$,so that 
$$ \omega_\kappa(\xi_1, \xi_2) 
= [\kappa(\xi_1, \nabla\xi_2)] 
= [\kappa(\xi_1, \dd\xi_2) + \kappa(\theta,[\xi_2, \xi_1])]. $$

If $X = X_\xi \in \aut(P)$ corresponds to $\xi \in \gau(P)$ in the sense that 
$\theta(X_\xi) = -\xi$, then 
${\cal L}_{X_\xi}\xi' = -[\theta(X_\xi),\xi'] = [\xi,\xi']$
and the fact that the curvature $\dd\theta + \frac{1}{2}[\theta,\theta]$ is 
a basic $2$-form leads to 
\begin{align*}
{\cal L}_{X_\xi}\theta
&= \dd(i_{X_\xi}\theta) + i_{X_\xi}\dd\theta 
= - \dd\xi  - \frac{1}{2}i_{X_\xi}[\theta,\theta] 
= - \dd\xi  - [\theta(X_\xi),\theta] \\
&= - \dd\xi  + [\xi, \theta] = -\nabla \xi, 
\end{align*}
so that $\omega_\kappa$ can also be written as 
$$ \omega_\kappa(\xi_1, \xi_2) = -[\kappa(\xi_1, {\cal L}_{X_{\xi_2}}\theta)]. $$

The group $\Aut(P)$ acts on the affine space 
${\cal A}(P) \subeq \Omega^1(P,\fk)$ of principal connection 
$1$-forms by $\phi.\theta := (\phi^{-1})^*\theta$ and on 
$\g = \gau(P)$ by $\phi.\xi = \xi \circ \phi^{-1}$. 
We then have 
$$ \phi.\nabla\xi 
= \phi.(\dd \xi + [\theta,\xi])
= \dd(\phi.\xi) + [\phi.\theta,\phi.\xi])
= \nabla(\phi.\xi) + [\phi.\theta-\theta,\phi.\xi]. $$
This leads to 
$$ (\phi.\omega_\kappa)(\xi_1, \xi_2) 
=  \phi.\omega_\kappa(\phi^{-1}.\xi_1, \phi^{-1}.\xi_2) 
=  \omega_\kappa(\xi_1, \xi_2) + [\kappa(\phi.\theta-\theta,[\xi_2, \xi_1])].$$
Note that 
$$\zeta \: \Aut(P) \to \Omega^1(M,\Ad(P)),\quad \phi \mapsto 
\phi.\theta -\theta $$
is a smooth $1$-cocycle, so that 
$$ \Psi \: \Aut(P) \to \Hom(\g,\oline\Omega^1(M,\bV)), \quad 
\Psi(\phi)(\xi) := [\kappa(\phi.\theta -\theta,\xi)] $$
is a $1$-cocycle with 
$\dd_\g(\Psi(\phi)) = \phi.\omega_\kappa - \omega_\kappa.$ 
Therefore Lemma~\ref{lem:aut-lift} below implies that 
$$ \phi.([\alpha],\xi) := ([\phi.\alpha] + [\kappa(\phi.\theta-\theta,\phi.\xi)], 
\phi.\xi) $$
defines an automorphism of the central extension $\hat\g$ and it is easy to 
see that we thus contains a smooth action of $\Aut(P)$ on $\hat\g$ 
(cf.\ \cite{NeWo07}). 

If $\phi_f \in \Gau(P)$ is a gauge transformation corresponding 
to the smooth function $f \: P \to K$, then 
$\phi_f^*\theta = \delta^l(f) + \Ad(f)^{-1}.\theta$ implies 
$$ \phi_f.\theta 
= \delta^l(f^{-1}) + \Ad(f).\theta 
= -\delta^r(f^{-1}) + \Ad(f).\theta. $$ 

(b) If $P = M \times K$ is the trivial bundle and 
$p_M \: P \to M$ and $p_K \: P \to K$ are the two projection maps, then 
the identification of gauge and mapping groups can be written as 
$$C^\infty(M,K) \to \Gau(P) \subeq C^\infty(P,K), \quad 
f \mapsto \tilde f := p_K^{-1} \cdot p_M^*(f) \cdot p_K, $$
and the canonical connection $1$-form is $\theta = \delta^l(p_K) 
\in \Omega^1(P,\fk)$. 
We further write $\tilde \alpha := \Ad(p_K)^{-1}.p_M^*\alpha$ for the 
realization of $\alpha \in \Omega^1(M,\fk)$ as an element of 
$\Omega^1(M,\Ad(P)) \subeq \Omega^1(P,\fk)$. 

If $\tilde\phi\in \Aut(P)$ is the lift of $\phi \in \Diff(M)$ defined by 
$\tilde\phi(m,k) := (\phi(m),k)$, then $p_K\circ \tilde\phi= p_K$ 
implies that $\tilde\phi.\theta = \theta$, so that 
$\zeta(\phi) = 0$. 
To determine $\zeta$ on the gauge group, we calculate: 
\begin{align*}
\phi_{\tilde f}.\theta 
&= -\delta^r(\tilde f) + \Ad(\tilde f).\theta 
= -\delta^r(\tilde f) - \Ad(\tilde f).\delta^r(p_K^{-1}) \\ 
&= -\delta^r(\tilde f\cdot p_K^{-1})  = -\delta^r(p_K^{-1} \cdot p_M^*(f)) 
=-\delta^r(p_K^{-1}) - \Ad(p_K)^{-1}.\delta^r(p_M^*(f))\\ 
&=\delta^l(p_K) - \Ad(p_K)^{-1}.(p_M^*\delta^r(f))  
= \theta  - \delta^r(f)\,\tilde.   
\end{align*}
In this sense, the 
restriction of the cocycle $\zeta$ to $\Gau(P)$ is 
$\zeta(\phi_{\tilde f}) = - \delta^r(f)\,\tilde{}$, so that it corresponds 
to the $1$-cocycles defined by the right logarithmic derivative 
$$ \delta^r \: C^\infty(M,K) \to \Omega^1(M,\fk). $$
We therefore recover the formulas 
for the actions of $C^\infty(M,K) \rtimes \Diff(M)$ on 
$\hat\g$ from Remarks~\ref{rem:3.3} and \ref{rem:covmap}. 

For later reference, we also note that the 
action of the gauge group on $1$-forms is given in these terms by 
\begin{align*}
&\phi_{\tilde f}.\tilde\alpha 
= (\phi_{\tilde f}^{-1})^*\tilde\alpha 
= (\phi_{\tilde f}^{-1})^*\big(\Ad(p_K^{-1}).(p_M^*\alpha)\big) \\
&= \Ad(p_K \cdot \tilde f^{-1})^{-1}.(p_M \circ \phi_{\tilde f}^{-1})^*\alpha 
= \Ad(p_K^{-1}\cdot p_M^*f)p_M^*\alpha \\
&= \Ad(p_K^{-1}).p_M^*(\Ad(f).\alpha) 
= \big(\Ad(f).\alpha\big)\, \tilde{} 
= \Ad(\tilde f).\tilde\alpha. 
\end{align*}
\end{remark}

\begin{remark} (Covariance for type (II)) 
In the type (II) situation we have constructed the Lie algebra
extension $\hat\g = \hat\gau(P)$ as $\Gamma\hat\Ad(P)$, and since 
$\hat\Ad(P)$ is a Lie algebra bundle associated to $P$ via an 
action of $K$ on $\hat\fk$ by automorphisms, 
the group $\Aut(P)$ and its Lie algebra $\aut(P)$ 
act naturally on $\hat\gau(P)$ by automorphisms, resp., derivations, 
lifting the action on $\gau(P)$. 
\end{remark}

\begin{problem} (Covariance for type (III)) 
It would be interesting to see if the type (III) cocycles on 
$\gau(P)$ with values in $\Omega^1(M,V)$ are $\aut(P)$-covariant 
in the sense that $\aut(P)$ acts on $\hat\gau(P)$. 
\end{problem}

\begin{problem} \label{prob:cross-mod} (Corresponding crossed modules) 
For cocycles of type (I) and (II) we have seen that the action of 
$\aut(P)$ on $\gau(P)$ lifts to an action by derivations on the 
central extension $\hat\gau(P)$, which defines a crossed module 
$$ \mu \: \hat\gau(P) \to \aut(P) $$
of Lie algebras. 

For type (I), the characteristic class of this crossed module is an element 
of $H^3({\cal V}(M),\oline\Omega^1(M,V))$ and for type (II) in 
$H^3({\cal V}(M),C^\infty(M,V))$. The vanishing of these characteristic 
classes is equivalent to the embeddability of the central extension
$\hat\gau(P)$ to an abelian extension of $\aut(P)$ 
(cf.\ \cite{Ne06b}, Thm.~III.5). 

For type (II) it has been shown in \cite{Ne06b}, Lemma~6.2, that the characteristic 
class can be represented by a closed $3$-form, so that it comes from 
an element of $H^3_{\rm dR}(M,V)$. 
It seems that a proper understanding of the corresponding situation for 
type (I) requires a descriptions of the third cohomology space \break 
$H^3({\cal V}(M), \oline\Omega^1(M,\R))$, so that we need a generalization 
of Theorem~\ref{thm:h2vect} for cohomology in degree $3$. 

For type (III) cocycles, one similarly needs the cohomology spaces \break 
$H^3({\cal V}(M), \Omega^1(M,\R))$, but they can be described with 
Tsujishita's work (\cite{Tsu81}, Thm.~5.1.6; \cite{BiNe07}, Thm.~III.1) 
which asserts that the space  
$H^\bullet({\cal V}(M), \Omega^\bullet(M,\R))$ is a free module of 
$H^\bullet({\cal V}(M), C^\infty(M,\R))$, generated by the non-zero 
products of the cohomology classes $[\Psi_k]$. This leads to 
$$H^3({\cal V}(M), \Omega^1(M,\R))
\cong H^2({\cal V}(M), C^\infty(M,\R))\cdot [\Psi_1]. $$
If $TM$ is trivial, we obtain with Theorem~\ref{thm:h2vect} 
$$H^3({\cal V}(M), \Omega^1(M,\R))
= H^1_{\rm dR}(M,\R) \cdot [\Psi_1] 
\oplus H^2_{\rm dR}(M,\R) \cdot [\oline\Psi_1 \wedge \Psi_1]. $$
\end{problem}

\subsection{Corresponding Lie group extensions} 
\label{sec:4.3}

{\bf (I):} For cocycles 
of the form $\omega_\kappa \in Z^2(\gau(P),\oline\Omega^1(M,\bV))$, 
it seems quite difficult to decide their integrability. However, if $\pi_0(K)$ is finite, we have the 
following generalization of Theorem~\ref{thm:int-type1} (\cite{NeWo07}): 

\begin{theorem} \label{thm:int-type1-gauge} 
If $\pi_0(K)$ is finite and $G := \Gau(P)_0$, then the following are equivalent: 
\begin{description}
\item[\rm(1)] $\omega_{\kappa}$ integrates for each principal 
$K$-bundle $P$ over a compact manifold $M$ to a Lie group extension of $G$. 
\item[\rm(2)] $\omega_{\kappa}$ integrates for the trivial 
$K$-bundle $P = \bS^1 \times K$ over $M = \bS^1$ 
to a Lie group extension of $G$. 
\item[\rm(3)] The image of $\per_{\kappa} \: \pi_3(K) \to V$ is discrete. 
\end{description}
These conditions are satisfied if $\kappa$ is the universal invariant 
bilinear form with values in $V(\fk)$. 
\end{theorem}

As we have already seen for trivial bundles, the cocycle 
$\dd \circ \omega_\kappa$ with values in $\Omega^2(M,\bV)$ can 
be integrated quite explicitly. We now explain the geometric background 
for that in the general case. We start with a general group theoretic remark. 

\begin{remark} Let $E$ and $F$ be vector space and 
$\omega \: E \times E \to F$ an alternating bilinear form. 

(a) We write $H := H(E,F,\omega)$ for the corresponding 
{\it Heisenberg group}, which is the product set $F \times E$, 
endowed with the product 
$$ (f_1,e_1)(f_2,e_2) = (f_1 + f_2 + \omega(e_1, e_2), e_1+ e_2). $$
The invariance group of $\omega$ is 
$$ \Sp(E,F,\omega) := \{ (\phi,\psi) \in \GL(E) \times \GL(F) \: 
\phi^*\omega = \psi \circ \omega\}, $$
and the natural diagonal action of this group on $F \times E$ preserves 
$\omega$, so that we may form the semidirect product group 
$$ \HSp(E,F,\omega) := H(E,F,\omega) \rtimes \Sp(E,F,\omega) $$
with the product 
$$ (f_1, e_1, \phi_1, \psi_1) (f_2, e_2, \phi_2, \psi_2) 
= (f_1 + \psi_1.f_2 + \omega(e_1, \phi_1.e_2), e_1 + \phi_1.e_2, \phi_1\phi_2, 
\psi_1\psi_2). $$
This group is an abelian extension of the group 
$$E \rtimes \Sp(E,F,\omega) \subeq \Aff(E) \times \GL(F) $$ 
by $F$ and the corresponding cocycle is 
$$ \Omega((e_1,\phi_1,\psi_1),(e_2,\phi_2,\psi_2)) 
= \omega(e_1, \phi_1.e_2). $$

(b) Now let $G$ be a group, 
$\rho_E \: G \to \Aff(E)$ and $\rho_F \: G \to \GL(F)$ 
group homomorphisms such that $\rho := (\rho_E, \rho_F)$ maps 
$G$ into $E \rtimes \Sp(E,F,\omega)$. Then the pull-back 
$\hat G := \rho^*\HSp(E,F,\omega)$ 
is an abelian extension of $G$ by $V$. To make the cocycle 
$\omega_G$ of this extension explicit, we write 
$\rho_E(g) = (\mu(g), \phi(g))$ and obtain 
$$ \omega_G(g_1, g_2) = \omega(\mu(g_1), \rho(g_1).\mu(g_2)). $$
This $2$-cocycle $\omega_G \in Z^2(G,F)$ is the cup-product 
$\omega_G = \mu \cup_\omega \mu$ of the 
$1$-cocycle $\mu \in Z^1(G,E)$ with itself with respect to the 
equivariant bilinear map $\omega$ (cf.\ \cite{Ne04a}, App.~F). 

If we think of $E$ as an abstract affine space without specifying an 
origin, then any point $o \in E$ leads to a $1$-cocycle 
$\mu(g) := g.o - o$ with values in the translation group $(E,+)$ and, 
accordingly, we get 
$$ \omega_G(g_1, g_2) = \omega(g_1.o - o, \rho(g_1)(g_2.o - o)). $$ 
\end{remark}

\begin{remark} On the translation space $\Omega^1(M,\Ad(P))$ of the 
affine space ${\cal A}(P)$ of principal connection $1$-forms, 
any $V$-valued $K$-invariant symmetric bilinear map $\kappa$ on $\fk$ 
defines an alternating $C^\infty(M,\R)$-bilinear map 
$$ \Omega^1(M,\Ad(P)) \times \Omega^1(M,\Ad(P)) \to \Omega^2(M,\bV),
\quad (\alpha,\beta) \mapsto \alpha \wedge_\kappa \beta, $$
so that ${\cal A}(P)$ 
is an affine space with a translation invariant 
$\Omega^2(M,\bV)$-valued $2$-form that is invariant under the affine action 
of $\Aut(P)$ by pull-backs. 
\end{remark}

\begin{proposition} The smooth $2$-cocycle 
$$ \zeta_\kappa \in Z^2(\Aut(P), \Omega^2(M,\bV)), \quad 
\zeta_\kappa(\phi_1, \phi_2) 
= (\phi_1.\theta - \theta) \wedge_\kappa \phi_1.(\phi_2.\theta - \theta) $$
defines an abelian extension of $\Aut(P)$ by $\Omega^2(M,\bV)$ whose corresponding 
Lie algebra cocycle is given by 
$$ (X_1, X_2) \mapsto 
({\cal L}_{X_1} \theta) \wedge_\kappa ({\cal L}_{X_2} \theta) 
- ({\cal L}_{X_2} \theta) \wedge_\kappa ({\cal L}_{X_1} \theta) 
= 2({\cal L}_{X_1} \theta) \wedge_\kappa ({\cal L}_{X_2} \theta). $$
\end{proposition}

\begin{remark} If $P = M \times K$ is a trivial bundle, then we  
obtain on $\Gau(P) \cong C^\infty(M,K)$ with Remark~\ref{rem:cov-gau}(b): 
$$ \zeta_\kappa(\phi_{\tilde f}, \phi_{\tilde g}) 
=  \delta^r(\tilde f) \wedge_\kappa \Ad(\tilde f).\delta^r(\tilde g)
=  \delta^l(\tilde f) \wedge_\kappa \delta^r(\tilde g)
=  \big(\delta^l(f) \wedge_\kappa \delta^r(g)\big)\, \tilde{}. $$
\end{remark}

\begin{problem} From the existence of the cocycle $\zeta_\kappa$, 
it follows that $\dd \circ c = 0$ 
holds for the obstruction class $c$ of the corresponding crossed module 
(cf.\ Problem~\ref{prob:cross-mod}).  
Therefore the long exact cohomology sequence associated to the short exact 
sequence 
$$ \0 \to H^1_{\rm dR}(M,\bV) \into \oline\Omega^1(M,\bV) \to B^2_{\rm dR}(M,\bV)\to0 $$
implies that $c$ lies in the image of the space 
$H^3({\cal V}(M), H^1_{\rm dR}(M,\bV))$. As the action of 
${\cal V}(M)$ on $H^1_{\rm dR}(M,\bV)$ is trivial, we have to consider 
$H^3({\cal V}(M),\R)$. This space vanishes for $\dim M \geq 3$ 
(cf.\ \cite{BiNe07}, Cor.~D.6, resp., \cite{Hae76}, Thm.~3/4), so that 
in this case the central extension $\hat\gau(P)$ corresponding to 
$\omega_\kappa$ embeds into an abelian extension of $\aut(P)$. 
It is an interesting problem to describe these extensions more explicitly. 
Are there any $2$-dimensional manifolds for which the obstruction 
class in $c$ is non-zero? 
\end{problem}

{\bf (II):} Now we turn to cocycles of type (II). 
As before, we fix a central Lie group extension $\hat K$ of $K$ 
by $Z \cong V/\Gamma_Z$. Then there is an obstruction class 
$\chi_{\hat K} \in \check H^2(M,\uline{Z})$ with values in the sheaf 
of germs of smooth $Z$-valued functions that vanishes if and only if 
there exists a $\hat K$-principal bundle $\hat P$ with $\hat P/Z \cong P$, 
as a $K$-bundle (\cite{GP78}). Using the isomorphisms 
$$\check H^2(M,\uline{Z}) 
\cong \check H^3(M,\Gamma_Z) 
\cong H^3_{\rm sing}(M,\Gamma_Z), $$
we see that $\chi$ maps to some class $\oline\chi \in H^3_{\rm dR}(M,V)$, 
and it has been shown in \cite{LW06} that $\oline\chi$ yields 
(up to sign) the characteristic class $\chi(\mu)$ of the 
Lie algebra crossed module 
$\mu \: \hat\gau(P) \to \aut(P).$ 
(cf.\ Theorem~\ref{thm:lec1}).

On the level of de Rham classes, the following theorem can also be found in 
\cite{LW06}; the present version is a refinement. 

\begin{theorem} 
Let $K$ be a connected finite-dimensional Lie group 
and \break $q_K \: \hat K \to K$ a central Lie group extension by $Z$ with 
$Z_0 \cong V/\Gamma_Z$ and $q \: P\to M$ a principal $K$-bundle. 
Then the following assertions hold: 
\begin{description}
\item[\rm(a)] If $Z$ is connected and $\pi_2(K)$ vanishes, which is the 
case if $\dim K < \infty$, then there is a $\hat K$-bundle 
$\hat P$ with $\hat P/Z \cong P$ and we obtain an abelian Lie group extension 
$$\1 \to C^\infty(M,Z) \to \Aut(\hat P) \to \Aut(P) \to \1 $$ 
containing the central extension
$$\1 \to C^\infty(M,Z) \to \Gau(\hat P) \to \Gau(P) \to \1 $$ 
integrating the Lie algebra $\hat\gau(P)$. 
\item[\rm(b)] If $\dim K < \infty$, then the 
obstruction class in $H^3_{\rm dR}(M,V)$ vanishes. 
\item[\rm(c)] If $K$ is $1$-connected and $\hat K$ is connected, then $Z$ is 
connected. 
\end{description}
\end{theorem}

\begin{proof} (a) That $\pi_2(K)$ vanishes if $K$ is finite-dimensional 
is due to \'E.~Cartan (\cite{CaE36}). 
If $\pi_2(K)$ vanishes and $Z$ is connected, then 
Theorem~7.12 of \cite{Ne02} implies that the central extension 
$\hat K$ of $K$ by $Z\cong V/\Gamma_Z$ can be lifted to a central extension 
$K^\sharp$ of $K$ by $V$ with $\hat K \cong K^\sharp/\Gamma_Z$. 
This implies that the obstruction class to lifting the structure 
group of the bundle $P$ to $\hat K$, which is an element of the 
sheaf cohomology group  
$\check H^2(M,\uline Z)$, is the image of an element of 
$\check H^2(M,\uline{V}) =  \{0\}$, hence trivial. 
Here we have used that the sheaf cohomology groups 
$\check H^p(M,\uline{V})$, $p > 0$, vanish because the sheaf of germs 
of smooth maps with values in the vector space $V$ is soft 
(cf.\ \cite{God73}). 

(b) Since $Z_0$ is divisible, we have a direct product decomposition 
$Z \cong Z_0 \times \pi_0(Z)$ and, accordingly, the central extension 
$\hat K$ is a Baer sum of an extension by $Z_0$ and an extension 
by $\pi_0(Z)$. According to (a), the extension by $Z_0$ does not contribute 
to the obstruction class, so that it is an element of 
$H^2(M,\pi_0(Z))$, and therefore the corresponding de Rham class in $H^3(M,V)$ 
vanishes. 

(c) The group $\hat K/Z_0$ is a connected covering group of $K$, hence trivial, 
and therefore $Z = Z_0$. 
\smartqed\qed
\end{proof}

\begin{remark} If $\hat K$ is a central $\T$-extension of some infinite-dimensional 
group for which the corresponding de Rham obstruction class in 
$H^3_{\rm dR}(M,\R)$ is non-zero, Lecomte's Theorem~\ref{thm:lec1}  
implies that the corresponding class in 
$H^3({\cal V}(M), C^\infty(M,\R))$ is also non-zero because 
the degrees on the Pontrjagin classes are multiples of $4$, so that the 
ideal they generate vanishes in degrees $\leq 3$. 
This means that the class of the corresponding crossed module 
$\mu \: \hat\gau(P) \to \aut(P)$ is non-trivial, hence 
the central extension $\hat\gau(P)$ does not embed equivariantly into an 
abelian extension of $\aut(P)$. 
\end{remark}

In a similar fashion as the group $H^2(M,\Z)$ classifies $\T$-bundles 
over $M$, the group $H^3(M,\Z)$ classifies certain geometric objects 
called {\it bundle gerbes} (\cite{Bry93}, \cite{Mu96}). 
As the preceding discussion shows, non-trivial 
bundle gerbes can only be associated to the central $\T$-extensions of 
infinite-dimensional Lie groups. For more on bundle gerbes and 
compact groups we refer to \cite{SW07} in this volume.

\section{Multiloop algebras} \label{sec:loopalg} 

In this section we briefly discuss the connection between gauge algebras of 
flat bundles and the algebraic concept of a (multi-)loop algebra. 
These are infinite-dimensional Lie algebras which 
are presently under active investigation from the algebraic point of view 
(\cite{ABP06}, \cite{ABFP07}). To simplify matters, we discuss only the 
case where $\fk$ is a complex Lie algebra. 

\subsection{The algebraic picture} 

Let $\fk$ be a complex Lie algebra, $m_1,\ldots, m_r \in \N$, 
$\zeta_i \in \T$ with $\zeta_i^{m_i} = 1$ and 
$\sigma_1,\ldots, \sigma_r \in \Aut(\fk)$ with $\sigma_i^{m_i} = \id_\fk$. 
Then we obtain an action of the 
group $\Delta := \Z/m_1 \times \ldots \times \Z/m_r$ on the Lie algebra 
$$ \fk \otimes \C[t_1,t_1^{-1}, \ldots, t_r, t_r^{-1}] $$
by letting the $i$-th generator $\hat\sigma_i$ of $\Delta$ acts by 
$$ \hat\sigma_i(x \otimes t^{\alpha}) := \sigma_i(x) \otimes \zeta_i^{\alpha_i}t^{\alpha}. $$
The algebra 
$$ M(\fk, \sigma_1, \ldots, \sigma_r) 
:= (\fk \otimes \C[t_1,t_1^{-1}, \ldots, t_r, t_r^{-1}])^\Delta $$
of fixed points of this action is called the corresponding {\it multi-loop algebra}. 
For $r = 1$ we also write 
$L(\fk,\sigma) := M(\fk,\sigma)$ and call it a {\it loop algebra}. 

The loop algebra construction can be iterated: If we start with 
$\sigma_1 \in \Aut(\fk)$ with $\sigma_1^{m_1} = \id$ and pick 
$\sigma_2 \in \Aut(L(\fk,\sigma_1))$ with $\sigma_2^{m_2} = \id$, we put 
$$L(\fk,\sigma_1, \sigma_2) := L(L(\fk,\sigma_1), \sigma_2). $$
Repeating this process, we obtain the {\it iterated loop algebras} 
$$ L(\fk,\sigma_1, \ldots, \sigma_r) 
:= L(L(\fk,\sigma_1, \ldots, \sigma_{r-1}), \sigma_r), $$
where $\sigma_j \in \Aut(L(\fk,\sigma_1,\ldots, \sigma_{j-1}))$ is assumed to satisfy 
$\sigma_j^{m_j} =\id$. Clearly, every multiloop algebra can also be described as 
an iterated loop algebra, but the converse is not true 
(cf.\ Example~\ref{ex:9.7} below). 

\subsection{Geometric realization of multiloop algebras} 

There is a natural analytic variant of the loop algebra construction. 
Here the main idea is that the analytic counterpart of the 
algebra $\C[t,t^{-1}]$ of Laurent polynomials is the Fr\'echet algebra 
$C^\infty(\bS^1,\C)$ of smooth complex-valued functions on the circle. Its dense subalgebra of 
trigonometric polynomials is isomorphic to the algebra of Laurent polynomials, 
so that it is a completion of this algebra. A drawback of this completion 
process is that the algebra of smooth function has zero divisors, but this 
in turn facilitates localization arguments. 

The analytic analog of the algebra $\fk \otimes \C[t_1,t_1^{-1},\ldots, t_r, t_r^{-1}]$ 
is the Fr\'echet--Lie algebra $C^\infty(\T^r,\fk)$ of smooth $\fk$-valued functions 
on the $r$-dimensional torus $\T^r$, endowed with the pointwise Lie bracket. 
To describe the analytic analogs of the multi-loop algebras, let 
$\Delta := \Z/m_1 \times \ldots \times \Z/m_r $ 
and fix a homomorphism $\rho \: \Delta \to \Aut(\fk)$, which is determined by 
a choice of automorphisms $\sigma_i$ of $\fk$ with $\sigma_i^{m_i} = \id_\fk$. 
Let ${\bf m} := (m_1,\ldots, m_r)$ and consider the corresponding torus 
$$ \T^r_{{\bf m}} := (\R/m_1 \Z) \times \ldots \times (\R/m_r \Z), $$
on which the group $\Delta \cong \Z^r/(m_1\Z \oplus \ldots \oplus m_r \Z)$ acts 
in the natural fashion from the right 
by factorization of the translation action of $\Z^r$ on $\R^r$. 
Now we obtain an action of $\Delta$ on the Lie algebra  
$C^\infty(\T^r_{\bf m},\fk)$ by 
$$ (\gamma.f)(x) := \gamma.f(x.\gamma) $$
Then we have a dense embedding 
$M(\fk, \sigma_1,\ldots, \sigma_r) \into C^\infty(\T^r_{\bf m},\fk)^\Delta.$
To realize this Lie algebra geometrically, we note that the 
quotient map $q \: \T^r_{\bf m} \to\T^r$ is a regular covering for which 
$\Delta$ acts as the group of deck transformations on $\T^r_{\bf m}$. 
Then the homomorphism $\rho \: \Delta \to \Aut(\fk)$ from above defines a flat Lie algebra 
bundle 
$$ \fK := (\T^r_{\bf m} \times\fk)/\Delta, \quad [t,x] \mapsto q(t), $$
over $\T^r$, where $\Delta$ acts by $\gamma.(t,x) := (t.\gamma^{-1}, \gamma.x)$, and 
the space $\Gamma\fK$ of smooth sections of this Lie algebra bundle can be 
realized as 
$$ \{ f \in C^\infty(\T^r_{\bf m},\fk) \: (\forall \gamma \in \Delta)
(\forall x \in \T^r_{\bf m})\ f(x.\gamma) = \gamma^{-1}.f(x) \} = C^\infty(\T^r_{\bf m},\fk)^\Delta. $$
In this sense, the Lie algebra $\Gamma\fK$ is a natural geometric analog 
of a multiloop algebra. 

\subsection{A generalization of multiloop algebras} 

Our geometric realization of multiloop algebras suggests the following 
geometric generalization. Let $\Delta$ be a discrete group and 
$q \: \hat M \to M$ a regular covering, where $\Delta$ acts from the right 
via $(t,\gamma) \mapsto t.\gamma$ 
on $\hat M$ as the group of deck transformations. Let $\fk$ be a locally convex 
Lie algebra and $\rho \: \Delta \to \Aut(\fk)$ an action 
of $\Delta$ by topological automorphisms of $\fk$. 
Let $K_0$ be a $1$-connected (regular or locally 
exponential) Lie group with Lie algebra $\fk$. Then the action of 
$\Delta$ integrates to an action $\rho_K \: \Delta \to \Aut(K_0)$ 
by Lie group automorphisms (\cite{GN07}), 
so that we may form the semidirect product Lie group 
$K := K_0 \rtimes \Delta$ with $\pi_0(K) = \Delta$. 

Next we form the flat Lie algebra bundle 
$$ \fK := (\hat M \times\fk)/\Delta, \qquad q_\fK \: \fK \to M,\ \ 
[t,x] := \Delta.(t,x)\mapsto q(t), $$
over $M$, where $\Delta$ acts on $M \times \fK$ 
by $\gamma.(t,x) := (t.\gamma^{-1}, \gamma.x)$. 
Each section of $\fK$ can be written as 
$s(q(t)) = [t, f_s(t)]$ with a smooth function $f \: \hat M \to \fk$, 
satisfying $f(t.\gamma) = \gamma^{-1}.f(t)$ for $t \in \hat M$, $\gamma \in \Delta$. 
We thus obtain a realization of the space $\Gamma\fK$ 
of smooth sections of this Lie algebra bundle as 
\begin{equation} \label{eq:subspace}
\{ f \in C^\infty(\hat M,\fk) \: (\forall \gamma \in \Delta)\ 
f(t.\gamma) = \gamma^{-1}.f(t) \} = C^\infty(\hat M,\fk)^\Delta, 
\end{equation}
where the $\Gamma$-action on $C^\infty(\hat M,\fk)$ is defined by 
$(\gamma.f)(t) := \gamma.f(t.\gamma)$. 

\begin{proposition} 
If $P := (\hat M \times K)/\Delta$ is the flat $K$-bundle associated
 to the inclusion $\Delta \into K$, 
then 
$\gau(P) \cong \Gamma\fK.$
\end{proposition}

Hence the Lie algebra of sections of a flat Lie algebra bundle 
can always be realized as a gauge algebra of a flat principal bundle; 
the converse is trivial. 

\begin{proof} Our assumptions on $K_0$ imply that 
$\Aut(K_0) \cong \Aut(\fk)$ (\cite{GN07}), so that 
the action of $\Delta$ on $\fk$ integrate to an action $\rho_K$ on $K_0$ 
by automorphisms. 

Clearly, $P \cong \hat M \times K_0$ is a trivial $K_0$-bundle, but it is not 
trivial as a $K$-bundle, because the group $\Delta$ acts on 
$\hat M \times K_0$ by 
$$ (\hat m,k).\gamma = (\hat m.\gamma, \rho_K(\gamma)^{-1}(k)). $$

Let $\Ad(P) = (P \times \fk)/K$ denote the corresponding 
gauge bundle with $\Gamma\Ad(P)\cong \gau(P)$. Then 
\begin{align*}
\gau(P) 
&\cong \{ f \in C^\infty(P,\fk) \: 
(\forall p \in P)(\forall k \in K)\ f(p.k) = \Ad(k)^{-1}.f(p)\} \\
&\cong \{ f \in C^\infty(\hat M,\fk) \: 
(\forall m \in \hat M)(\forall \gamma \in \Delta)\ f(\hat m.\gamma) 
= \Ad(\gamma)^{-1}.f(\hat m)\} \\
&= \{ f \in C^\infty(\hat M,\fk) \: 
(\forall m \in \hat M)(\forall \gamma \in \Delta)\ f(\hat m.\gamma) 
= \rho(\gamma)^{-1}.f(\hat m)\}. 
\end{align*}
\smartqed\qed
\end{proof}

\subsection{Connections to forms of Lie algebras over rings} 

In this subsection we assume that $\Delta$ is finite and describe briefly 
the connection to forms of Lie algebras over rings. For more details 
on this topic we refer to \cite{PPS07}. 
The Lie algebra $\Gamma\fK$ is a module of the algebra $R := C^\infty(M,\C)$ and 
the Lie bracket is $R$-bilinear. 
The pull-back bundle 
$$q^*\fK = \{ (t, x) \in \hat M \times \fK \: q(t) = q_\fK(x)\} $$
is trivial, because the map 
$\hat M \times \fk \to q^*\fK, (t,y) \mapsto (t,[t,y])$
is an isomorphism of Lie algebra bundles. Accordingly, the space 
$\Gamma(q^*\fK)$, which is a module of $S := C^\infty(\hat M,\C)$, 
is isomorphic to $\Gamma(\hat M \times \fk) \cong C^\infty(\hat M,\fk)$. 
On the level of sections, we have an $R$-linear embedding 
$$ q^* \: \Gamma\fK \to \Gamma(q^*\fK), \quad 
(q^*s)(t) := (t, s(q(t))), $$
which corresponds to the realization of $\Gamma\fK$ 
as the subspace (\ref{eq:subspace}) of $\Delta$-fixed points in 
$C^\infty(\hat M,\fk)$. 
Using the local triviality of the bundle and the finiteness of $\Delta$, 
one easily verifies that the canonical map 
$$ \Gamma\fK \otimes_R S \to \Gamma\hat\fK \cong C^\infty(\hat M,\fk), \quad 
s \otimes f \mapsto (q^*s) \cdot f $$
is a linear isomorphism (cf.\ Thm~7.1 in \cite{ABP06}). 
In this sense, the $R$-Lie algebra $\Gamma\fK$ is an $R$-form of 
the $S$-Lie algebra $C^\infty(\hat M,\fk)$.

\begin{problem} \label{prob:iter2} Let $\hat\otimes$ denote the completion 
of the algebraic tensor product as a locally convex space. 
Determine to which extent the map 
$$ \Gamma\fK \hat\otimes_{C^\infty(M,\K)} C^\infty(\hat M,\K)
\to \Gamma\hat\fK \cong C^\infty(\hat M,\fk), \quad 
s \otimes f \mapsto f \cdot (q^*s)  $$
is an isomorphism of Lie algebras if $\Delta$ is not necessarily finite. 
In the proof of Theorem~7.1 in \cite{ABP06},  
the finiteness of $\Delta$ is used 
to conclude that its cohomology with values in any rational 
representation space is trivial. 
\end{problem}

\begin{example} \label{ex:9.7} (cf.\ Example 9.7 in \cite{ABP06}) 
We consider $\fk = \fsl_n(\C)$ for $n \geq 2$ and 
the iterated loop algebra $L(\fk,\sigma_1, \sigma_2)$ defined 
as follows. As $\sigma_1 \in \Aut(\fk)$ we choose an involution 
not contained in the identity component of $\Aut(\fk)$, so that 
$L(\fk,\sigma_1)$ is a non-trivial loop algebra corresponding to the 
affine Kac--Moody algebra of type $A_{n-1}^{(2)}$ (cf.\ \cite{Ka90}, Ch.~8). 
As $\sigma_2$, we take the involution on 
$$L(\fk,\sigma_1) = \{ f\in C^\infty(\T, \fk) \: f(-t) = \sigma_1(f(t))\}, $$
where $\T$ is realized as the unit circle in $\C^\times$, 
defined by 
$$ \sigma_2(f)(t_1) := f(t_1^{-1}) $$
and put 
$$ L := L(\fk,\sigma_1, \sigma_2) = \{ f \in C^\infty(\T, L(\fk,\sigma_1)) \: :
f(-t_2) = \sigma_2(f(t_2))\}. $$
We may thus identify elements of $L$ with smooth functions 
$f \: \T^2 \to \fk$ satisfying 
$$ f(t_1, -t_2) = f(t_1^{-1}, t_2) \quad \hbox{ and } \quad 
f(-t_1, t_2) = \sigma_1(f(t_1, t_2)). $$

For $\tau_1(t_1, t_2) := (-t_1,t_2)$ and 
$\tau_2(t_1, t_2) := (t_1^{-1},-t_2)$ we thus obtain a fixed point free 
action of $\Delta := \la \tau_1, \tau_2 \ra$ on $\T^2$ and a representation 
$\rho \: \Delta \to \Aut(\fk)$, defined by $\rho(\tau_1) := \sigma_1$ and 
$\rho(\tau_2) := \id$. Then 
$$ L \cong \{ f \in C^\infty(\T^2,\fk) \: (\forall \gamma \in \Delta) 
(\forall t \in \T^2)\ f(t.\gamma) = \gamma^{-1}.f(t)\}. $$
This means that $L \cong \Gamma\fK$ for the Lie algebra bundle 
$\fK$ over $M := \T^2/\Delta$, where $M$ is the Klein bottle 
($\tau_2$ is orientation reversing). 
This implies that the centroid of $L$ is $C^\infty(M,\C)$ 
(\cite{Lec80}), 
hence not isomorphic to $C^\infty(\T^2,\C)$, so that $L$ cannot be a 
multiloop algebra. 
\end{example}

\section{Concluding remarks} 
\label{sec:6}

From the structure theoretic perspective on infinite-dimensional Lie groups, 
it is natural to ask for a classification of gauge groups of principal 
bundles and associated structural data, such as their central extensions, 
invariant forms etc. 

One readily observes that two equivalent 
$K$-principal bundles $P_i \to M$ have isomorphic gauge groups, so that a 
classification of gauge groups or automorphism groups of 
$K$-principal bundles can be achieved, in principle, 
by a classification of all $K$-principal bundles 
and then identifying those whose gauge groups are isomorphic. A key observation 
is that if a bundle $P$ is twisted by a $Z(K)$-bundle $P_Z$, then the 
corresponding gauge groups are isomorphic Lie groups, where the isomorphism 
on the Lie algebra level is $C^\infty(M,\R)$-linear. Taking this refined 
structure into account leads to the concept of 
(transitive) Lie algebroids and the perspective on $\Gau(P)$ 
as the group of sections of a Lie group bundle 
(cf.\ \cite{Ma89} for results in this direction and 
\cite{Ma05} for more on Lie groupoids). Presently the theory of Lie groupoids 
is only available for finite-dimensional structure groups $K$. It would 
be natural to develop it also for well-behaved infinite-dimensional 
groups over finite-dimensional manifolds.

If $K$ is finite-dimensional with finitely many connected components, 
then the classification of $K$-principal 
bundles can be reduced to the corresponding problem for a maximal compact 
subgroup, so that the well developed theory of bundles with 
compact structure group applies (\cite{Hu94}). 
If $K$ is an infinite-dimensional Lie group, it might not possess 
maximal compact subgroups
\begin{footnote}
{The representation theory of compact groups implies that the 
unitary group of an infinite-dimensional Hilbert space is such an 
example.}
\end{footnote}
so that a similar reduction does not work. 
However, recent results of M\"uller and Wockel 
show that each topological $K$-principal bundle over a finite-dimensional 
paracompact smooth manifold $M$ is equivalent to a smooth one 
and that two smooth $K$-principal bundles which are topologically equivalent 
are also smoothly equivalent (\cite{MW06}). 
Therefore the set $\Bun(M,K)$ of equivalence classes of smooth 
$K$-principal bundles over $M$ can be identified with the set 
$\Bun(M,K)_{\rm top}$ of topological $K$-bundles over $M$, and the latter 
set can be identified with the set $[M,BK]$ of homotopy classes of continuous maps 
of $M$ into a classifying space $BK$ of the topological group $K$ 
(\cite{Hu94}, Thms.~9.9, 12.2, 12.4). 

On the topological level, each continuous map 
$\sigma \: M \to BK$ defines an algebra homomorphism 
$$ H^\bullet(\sigma,\R) \: H^\bullet(BK,\R) \to H^\bullet(M,\R) 
\cong H^\bullet_{\rm dR}(M,\R), $$
whose range are the characteristic classes of the corresponding bundle $P$. 
The advantage of these characteristic classes is that they can be represented 
by closed differential forms on $M$, hence can be evaluated by integrals over 
compact submanifolds or other piecewise smooth singular cycles. 

These characteristic classes are quite well understood for finite-dimensional 
connected Lie groups $K$, 
because in this case essentially everything can be reduced 
to compact Lie groups, for which a powerful theory exists. 
However, for infinite-dimensional structure groups, the corresponding theory 
of characteristic classes has only been explored in very special cases. 
We refer to Morita's excellent textbook for more details 
and references (cf.\ \cite{Mo01}). 
Central extensions of gauge groups of bundles with infinite-dimensional fiber 
over the circle have recently been studied in the context of integrable 
systems in \cite{OR06}.

\section{Appendix A. Abelian extensions of Lie groups} 
\label{sec:7}

In this appendix we result some facts on the integration of Lie algebra 
$2$-cocycles from \cite{Ne04a}. They provide a general set of tools to integrate 
abelian extensions of Lie algebras to extensions of connected Lie groups.

Let $G$ be a connected Lie group and $V$ a Mackey complete $G$-module. 
Further, let $\omega \in Z^k(\g,V)$ be a $k$-cocycle and 
$\omega^{\rm eq} \in \Omega^k(G,V)$ be the corresponding left equivariant 
$V$-valued $k$-form with $\omega^{\rm eq}_\1 = \omega$. 
Then each continuous map $\bS^k \to G$ is homotopic to a smooth map, and 
$$ \per_\omega \: \pi_k(G) \to V^G, \quad 
[\sigma] \mapsto \int_\sigma \omega^{\rm eq} = \int_{\bS^k} \sigma^*\omega^{\rm eq} $$
defines the {\it period homomorphism} 
whose values lie in the $G$-fixed part of $V$ 
(\cite{Ne02}, Lemma 5.7). 

For $k = 2$ we define the {\it flux homomorphism}    
$$ F_\omega \: \pi_1(G) \to H^1(\g,V), \quad [\gamma] \mapsto [I_\gamma], $$
where we define for each piecewise smooth loop $\gamma \: \bS^1 \to G$ 
the $1$-cocycle 
$$ I_\gamma \: \g \to V,\quad I_\gamma(x) := 
\int_\gamma i_{x_r} \omega^{\rm eq} 
= \int_0^1 \gamma(t).\omega(\Ad(\gamma(t))^{-1}.x, \delta^l(\gamma)_t)\, dt, $$
where $x_r$ is the right invariant vector field on $x$ with $x_r(\1) = x$. 
If $V$ is a trivial module, then $\dd_\g V= \{0\}$, so that 
$H^1(\g,V)$ consists of linear maps $\g \to V$ and we may think of the 
flux as a map 
$F_\omega \: \pi_1(G) \times \g \to V.$

\begin{proposition} \label{prop:fluxcrit} {\rm(\cite{Ne02}, Prop.~7.6)} 
If $V$ is a trivial $G$-module, then the 
adjoint action of $\g$ on the central extension 
$\hat\g := V \oplus_\omega \g$  integrates to a
smooth action $\Ad_{\hat\g}$ of $G$ if and only if $F_\omega = 0$. 
\end{proposition}

\begin{theorem} \label{thm:abext}
Let $G$ be a connected Lie group, $A$ a smooth $G$-module 
of the form $A \cong \fa/\Gamma_A$, where $\Gamma_A \subeq \fa$ is a
discrete subgroup of the Mackey complete space $\fa$ 
and $q_A \: \fa \to A$ the quotient map. 
Then the Lie algebra extension $\hat\g := \fa \oplus_\omega \g$ of 
$\g$ defined by the cocycle $\omega$ integrates to an abelian Lie group 
extension of $G$ by $A$ if and only if 
$$ q_A \circ \per_\omega \: \pi_2(G) \to A 
\quad \mbox{ and } \quad 
F_\omega \: \pi_1(G) \to H^1(\g,\fa) $$
vanish. 
\end{theorem}

\begin{remark} If only $q_A \circ \per_\omega$ vanishes, then 
the preceding theorem applies to the simply connected covering group 
$\tilde G$ of $G$, so that we only obtain an extension of 
$G$ by a non-connected group $A^\sharp$ which itself is a central extension of 
the discrete group $\pi_1(G)$ by $A$. For a more detailed discussion of these 
aspects, we refer to Section~7 of \cite{Ne04a}. 
\end{remark}

\begin{remark} \label{rem:8.4} To calculate period and flux homomorphisms, it is often 
convenient to use related cocycles on different groups. So, let 
us consider a morphism $\phi \: G_1 \to G_2$ of Lie groups and 
$\omega_i \in Z^2(\g_i, V)$, $V$ a trivial $G_i$-module, satisfying 
$\L(\phi)^*\omega_2 = \omega_1$. Then a straight forward argument 
shows that 
\begin{equation}
  \label{eq: per-trans}
\per_{\omega_2} \circ \pi_2(\phi) = \per_{\omega_1} \: \pi_2(G_1) \to V. 
\end{equation}
For the flux we likewise obtain 
\begin{equation}
  \label{eq: flux-trans}
F_{\omega_2} \circ (\pi_1(\phi) \times \L(\phi)) = F_{\omega_1} \: 
\pi_1(G_1) \times \g_1 \to V,  
\end{equation}
if we consider the flux as a bihomomorphism $\pi_1(G) \times \g\to V$.
\end{remark}

\begin{problem} \label{prob:8.3}
It is a crucial assumption in the preceding theorem 
that the group $G$ is connected and it would be very desirable to have a 
suitable 
generalization to non-connected groups $G$. 

The generalization to non-connected Lie groups $G$ means to 
derive accessible criteria for
the extendibility of a $2$-cocycle, resp., abelian extensions, 
from the identity component $G_0$ to the whole group $G$. 
From the short exact sequence 
$G_0 \into G \onto \pi_0(G)$, we obtain maps 
$$ H^2(\pi_0(G), A) \sssmapright{I} H^2(G, A) 
\sssmapright{R} H^2(G_0, A)^{G} = H^2(G_0, A)^{\pi_0(G)}, $$
but it seems to be difficult to describe the image of the restriction 
map $R$. 

If the identity component $G_0$ of $G$ splits, i.e., if 
$G \cong G_0 \rtimes \pi_0(G)$, then it is easy to see that 
$R$ is surjective, because the invariance of the cohomology class of 
$\hat G_0$ permits to lift the conjugation action of $\pi_0(G)$ to an action 
on $\hat G_0$, and then the semidirect product 
$\hat G_0 \rtimes \pi_0(G)$ is an abelian extension of $G$ by $A$, 
containing $\hat G_0$ as a subgroup.

If $A$ is a trivial module, one possible approach is to introduce 
additional structures on 
a central extension $\hat G_0$ of $G_0$ by $A$, so that the map 
$q \: \hat G_0 \to G$ describes a crossed module of groups 
(cf.\ \cite{Ne07}), 
which requires an extension of the 
natural $G_0$-action of $G$ on $\hat G_0$, resp., its Lie algebra $\hat\g$, 
to an action of $G$. 
If we have such an action lifting the conjugation action of $G$ on $G_0$, 
then the characteristic class 
$\chi(\hat G_0)$ of this crossed module is an element of the cohomology group 
$H^3(\pi_0(G), A)$ which vanishes if there exists an $A$-extension 
$\hat G$ of $G$ into which $\hat G_0$ embeds in a $\hat G$-equivariant fashion 
(cf.\ \cite{Ne07}, Thm.~III.8). 
\end{problem}

The following lemma is easy to verify (\cite{MN03}, Lemma~V.1; 
see also \cite{Ne06b}, Lemma~II.5) for a generalization to general 
Lie algebra extensions). 

\begin{lemma} \label{lem:aut-lift} Let 
$\hat\g = V \oplus_f \g$ be a central extension of 
the Lie algebra $\g$ and $\gamma = (\gamma_V, \gamma_\g) \in \GL(V) 
\times \Aut(\g)$. For $\theta \in C^1(\g,V)$ the formula 
$$ \hat\gamma(z,x) := (\gamma_\z(z) + \theta(\gamma_\g(x)), \gamma_\g(x)), 
\quad x \in \g, z \in V,$$
defines a continuous Lie algebra automorphism of $\hat\g$ 
if and only if 
$$ \gamma.\omega - \omega = \dd_\g \theta, \quad 
\mbox{ resp.} \quad 
\dd_\fn \tilde\theta = \omega - \gamma^{-1}.\omega \quad \mbox{ for } \quad 
\tilde\theta := \theta \circ \gamma_\g.$$ 
\end{lemma} 

The following theorem can be found in \cite{MN03}, Thm.~V.9: 

\begin{theorem} \label{thm:lifting-theorem} {\rm(Lifting Theorem)}  
Let $q \: \hat G \to G$ be a central Lie group extension of the $1$-connected 
Lie group $G$ by the Lie group $Z \cong \z/\Gamma_Z$. 
Let $\sigma_G \: H \times G \to G$, resp., $\sigma_Z \: H \times Z \to Z$ 
be smooth automorphic actions of the 
Lie group $H$ on $G$, resp., $Z$ and 
$\sigma_{\hat\g}$ a smooth action of $H$ on $\hat\g$ compatible with the 
actions on $\z$ and $\g$. Then there is a unique smooth action 
$\sigma_{\hat G} \: H \times \hat G \to \hat G$ 
by automorphisms compatible with the actions on $Z$ and $G$, for which 
the corresponding action on the Lie algebra $\hat\g$ 
is $\sigma_{\hat\g}$. 
\end{theorem}

\section{Appendix B. Abelian extensions of semidirect sums} 
\label{sec:8}

We consider a semidirect product of topological Lie algebras 
$\fh = \fn \rtimes_S \g$
and a topological $\fh$-module $V$. 
We are interested in a description of the cohomology space 
$H^2(\fh,V)$ in terms of $H^2(\fn,V)$ and $H^2(\g, V)$. 

To this end, we have study the {\it inflation}, resp., the {\it restriction map}  
$$ I \: H^2(\g, V^\fn) \to H^2(\h,V), \quad \mbox{resp.,} \quad 
R_\fn \: H^2(\fh,V) \to H^2(\fn,V)^\g, $$
satisfying $RI =0$. We further have a restriction map 
$$ R_\g  \: H^2(\fh,V) \to H^2(\g,V). $$
The composition 
$R_\g \circ I \: H^2(\g,V^\fn) \to H^2(\g,V)$
is the natural map induced by the inclusion $V^\fn \into V$ of $\g$-modules. 
If we interpret the elements of $H^2(\h,V)$ as abelian extensions of 
$\h$ by $V$, then the inflation map leads to twistings of these extensions 
by extensions of $\g$ by $V^\fn$.

We now construct an exact sequence describing kernel and cokernel of the 
map $(R_\fn, R_\g)$ which provides a quite accessible description of $H^2(\fh,V)$.

\begin{definition} 
(a) We write $C^p(\fn,V)$ for the space of continuous Lie algebra 
$p$-cochains with values in $V$ and 
$$C^p(\g,C^q(\fn,V))_c \subeq C^p(\g,C^q(\fn,V))$$ 
for the subspace 
consisting of those cochains defining a continuous 
$(p+q)$-linear map $\g^p \times \fn^q \to V$. Accordingly, we 
define $Z^p(\g,C^q(\fn,V))_c$, \break $H^p(\g,C^q(\fn,V))_c$ etc. 

(b) We further write 
$H^2(\fn,V)^{[\g]} \subeq H^2(\fn,V)^{\g}$ for the subspace 
of those cohomology classes $[f]$ for which there exists a 
$\theta \in C^1(\g,C^1(\fn,V))_c$ with 
\begin{equation}
  \label{eq:2.0} 
\dd_\fn(\theta(x)) = x.f \quad \mbox{ for } x \in \g.  
\end{equation}
Because of the continuity requirement for $\theta$, this 
is stronger that the $\g$-invariance 
of the cohomology class $[f] \in H^2(\fn,V)$. 
\end{definition}

\begin{remark} \label{rem:liftact} Any continuous $\g$-module action by derivations 
on a central extension $\hat \fn = V \oplus_f \fn$, compatible with the 
actions on $V$ and $\fn$ is of the form 
\begin{equation}
  \label{eq:action}
x.(v,n) = (x.v + \theta(x)(n), x.n),
\end{equation}
where $\theta \in Z^1(\g,C^1(\fn,V))_c$ satisfies (\ref{eq:2.0}). 
In particular, the existence of such an action implies that 
$[f] \in H^2(\fn,V)^{[\g]}$. 
\end{remark}

\begin{lemma} For $[f_\fn] \in H^2(\fn,V)^{[\g]}$ we choose $\theta \in 
C^1(\g,C^1(\fn,V))_c$ 
as in (\ref{eq:2.0}). We thus obtain a well-defined linear map 
$$ \gamma \: H^2(\fn,V)^{[\g]} \to H^2(\g,Z^1(\fn,V))_c, \quad [f_\fn] 
\mapsto [\dd_\g \theta].$$
\end{lemma}

\begin{proof} Using $\theta$, we obtain a linear map 
$$ \hat\psi\: \g \to \der(\hat\fn), \quad 
\hat\psi(x)(v,n) = (x.v + \theta(x)(n),x.n) $$
which is a linear lift of the given Lie algebra homomorphism 
$$ \psi = (S,\rho_V) \: \g \to (\der(\fn) \times \gl(V)).$$ 
In view of \cite{Ne06b}, Prop.~A.7, the corresponding 
Lie algebra cocycle on $\g$ is given by 
$\tilde f_\g = \dd_\g \theta \in Z^2(\g,Z^1(\fn,V))_c.$
The corresponding extension 
$$ \0 \to Z^1(\fn,V) \to \tilde \g = Z^1(\fn,V) \times_{\tilde f_\g} \g \to \g 
\to \0 $$
acts naturally on $\hat \fn$ by 
$(\alpha,x).(v,n) = (\alpha(n)+ x.v + \theta(x)(n),x.n)$. 

Next we observe that $\gamma$ is well-defined. 
Indeed, if $\theta' \in C^1(\g,C^1(\fn,V))_c$ also satisfies 
(\ref{eq:2.0}), 
then $\theta'-\theta \in C^1(\g,Z^1(\fn,V))_c$, so that 
$[\dd_\g(\theta'-\theta)]=0$ in 
$H^2(\g,Z^1(\fn,V))_c$. 

Moreover, if $\tilde f_\fn = f_\fn + \dd_\fn \beta$ for some
$\beta \in C^1(\fn,V)$, then 
$$ x.\tilde f_\fn = x.f_\fn + \dd_\fn(x.\beta) = \dd_\fn(\theta(x) + x.\beta), $$
so that $\tilde\theta(x) :=\theta(x) + x.\beta$ leads to 
$\dd_\fn(\tilde\theta(x)) = x.\tilde f_\fn$. Then 
$\dd_\g\tilde\theta = \dd_\g\theta + \dd_\g^2\beta = \dd_\g\theta$. 
This shows that $\gamma$ is well-defined. 
The linearity is clear. 
\smartqed\qed
\end{proof}

\begin{lemma} For each $\theta \in Z^1(\g,Z^1(\fn,V))_c$, 
$x.(v,n) := (x.v + \theta(x)(n), x.n)$ 
defines a continuous action 
of $\g$ on the semidirect sum $\hat\fn := V \rtimes \fn$, so that 
$\hat\h_\theta := \hat\fn \rtimes \g$ is an extension of $\fh$ by $V$. 
We thus obtain a well-defined map 
$$ \phi \: H^1(\g, Z^1(\fn,V)) \sssmapright{\phi} H^2(\h,V), \quad 
[\theta] \mapsto [\hat\fh_\theta]. $$ 
\end{lemma}

\begin{proof} A $2$-cocycle defining the extension $\hat\h_\theta$ is given by 
\begin{align*}
&\omega_\theta((n_1,x_1), (n_2, x_2))\\ 
&= [(0,(n_1, x_1)), (0,(n_2, x_2))] 
- (0, [n_1, n_2] + x_1.n_2 - x_2.n_1, [x_1, x_2]) \\
&= (0, \theta(x_1)(n_2) - \theta(x_2)(n_1)). 
\end{align*}
If $\theta$ is a coboundary, i.e., there exists a $\beta \in Z^1(\fn,V)$, 
with $\theta(x) = x.\beta$, $x \in \g$, then 
$\theta(x)(n) = x.\beta(n) - \beta(x.n)$
leads to 
\begin{align*}
\omega_\theta((n_1, x_1), (n_2, x_2)) 
&= \theta(x_1)(n_2) - \theta(x_2)(n_1) \\
&= x_1.\beta(n_2) - x_2.\beta(n_1) + \beta(x_1.n_2 - x_2.n_1).
\end{align*}
If $\tilde\beta \: \fh \to V$ denotes the continuous linear map 
extending $\beta$ and vanishing on $\g$, then $\omega_\theta$ is a coboundary: 
\begin{align*}
&(\dd_\fh\tilde\beta)((n_1, x_1),(n_2, x_2))\\ 
&= x_1.\beta(n_2) - x_2.\beta(n_1) + \dd_\fn(\beta)(n_1, n_2) 
- \beta(x_1.n_2 - x_2.n_1)\\
&= x_1.\beta(n_2) - x_2.\beta(n_1) - \beta(x_1.n_2 - x_2.n_1)
= \omega_\theta((n_1, x_1), (n_2, x_2)). 
\end{align*}
\smartqed\qed
\end{proof}

\begin{theorem} \label{thm:ex-sequ} With the 
linear map 
$$ \eta = (\dd_\fn)_* \: H^2(\g,V) \to H^2(\g,Z^1(\fn,V))_c, 
\quad [f_\g] \mapsto [-\dd_\fn \circ f_\g], $$
the following sequence is exact 
  \begin{eqnarray*}
&& H^1(\g, Z^1(\fn,V)) \sssmapright{\phi} H^2(\h,V) 
\smapright{(R_\fn, R_\g)} H^2(\fn,V)^{[\g]} \oplus H^2(\g,V) \\
&& \qquad\qquad\qquad \qquad \ssmapright{\gamma - \eta} H^2(\g, Z^1(\fn,V))_c. 
  \end{eqnarray*}
\end{theorem}

\begin{proof} To see that the sequence is exact in $H^2(\h,V)$, 
we simply note that a $V$-extension corresponds to a class 
$[f] \in \ker (R_\fn, R_\g)$ is and only if its restrictions to both 
subalgebras $\fn$ and $\g$ are trivial. 
This is equivalent to $\hat\fh$ being of the form 
$\hat\h = (V \rtimes \fn) \rtimes \g,$
hence equivalent to $\hat \fh_\theta$ 
for some $\theta \in Z^1(\g, Z^1(\fn,V))_c$. 

It remains to verify that $\im(R_\fn,R_\g) = \ker(\gamma-\eta)$. 
If $[f_\fn] \in \im(R_\fn)$, then there exists a topologically split 
extension $\hat\h$ of $\fh$ by the topological $\fh$-module~$V$. 
Then $\hat \h$ contains $\hat \fn$ as an 
ideal, so that the adjoint action defines an action of $\hat \h$ 
on $\hat \fn$. As $V$ acts non-trivially on $\hat \fn$, this action 
does in general not factor  through an action of $\fh$. 
Let $p \: \hat \fh \to \fh$ denote the quotient map. 
Then $\hat \fn = p^{-1}(\fn)$, and we put $\hat \g := p^{-1}(\g)$. 
As we have seen in Remark~\ref{rem:liftact}, the action of $\hat\g$ 
is described by some $\tilde\theta \in Z^1(\hat\g,C^1(\fn,V))_c$ 
satisfying 
\begin{equation}
  \label{eq:2.0b} 
\dd_\fn(\tilde\theta(x)) = x.f \quad \mbox{ for } x \in \hat\g.  
\end{equation}
Writing $\hat\g$ as $V \oplus_{f_\g} \g$, it follows that 
$\theta(x) := \tilde\theta(0,x)$ satisfies (\ref{eq:2.0}), so that 
$[f_\fn] \in H^2(\fn,V)^{[\g]}$. 
For $x,y \in \g$ we then obtain the relation 
\begin{align*}
x.\theta(y)-y.\theta(x) 
&= x.\tilde\theta(0,y) - y.\tilde\theta(0,x) 
= \tilde\theta([(0,x),(0,y)]) \\
&= -\dd_\fn(f_\g(x,y)) + \theta([x,y]), 
\end{align*}
i.e., $\dd_\g \theta = - \dd_\fn \circ f_\g.$ 
Hence, any element 
$([f_\fn],[f_\g]) \in \im(R_\fn,R_\g)$ is contained in the kernel of 
$$ \gamma - \eta \: H^2(\fn,V)^{[\g]} \oplus H^2(\g,V) \to 
H^2(\g,Z^1(\fn,V))_c. $$

If, conversely, $[f_\fn]\in \ker(\gamma-\eta)$, then there exists 
$\theta \in C^1(\g,C^1(\fn,V))_c$ with $\dd_\g \theta = - \dd_\fn \circ f_\g$, 
and then $\tilde\theta(v,x) := \theta(x) - \dd_n v$ defines a 
representation $\hat\g := V \oplus_{f_\g} \g \to \der \hat\fn$, 
compatible with the actions of $\g$ on $\fn$ and $V$. 

Then the semidirect product $\hat\fn\rtimes \hat\g$ is 
an extension of $\fn \rtimes \g = \fh$ by the space 
$V \times V$. Since both $V$-factors act in the same way on $\hat \fn$, 
the anti-diagonal $\oline\Delta_V$ acts trivially, hence commutes with $\hat \fn$. 
Since $\hat \g$, resp., $\g$, acts diagonally on $V \times V$, 
it also preserves the anti-diagonal. 
We thus obtain the $V$-extension 
$\hat \h := (\hat \fn \rtimes \hat \g)/\oline\Delta_V,$
of $\fh$. This extensions contains $\hat\fn$ and $\hat\g$ as subalgebras. 
This completes the proof. 
\smartqed\qed
\end{proof}

\begin{corollary} If $\fn$ is topologically perfect and $V = V^\fn$, 
then 
$$ (R_\fn, R_\g) \:  H^2(\h,V) \to H^2(\fn,V)^{[\g]} \oplus H^2(\g,V) $$
is a linear bijection. 
\end{corollary}

\begin{proof} In this case $Z^1(\fn,V)$ vanishes, and the assertion follows 
from Theorem~\ref{thm:ex-sequ}. 
\smartqed\qed
\end{proof}

\section{Appendix C. Triviality of the group action on Lie algebra 
cohomology} 
\label{sec:9}

\begin{theorem} \label{thm:9.1} Let 
$G$ be a connected Lie group with Lie algebra $\g$ and $V$ a Mackey complete 
smooth $G$-module. Then the natural action of $G$ on the 
Lie algebra cohomology $H^\bullet(\g,V)$ is trivial. 
\end{theorem}

\begin{proof} Since the space $H^\bullet(\g,V)$ carries no natural locally 
convex topology for which the action of $G$ on this space is smooth, 
we cannot simply argue that the triviality of the $G$-action follows 
from the triviality of the $\g$-action. 

The assertion is clear for $p = 0$, because $H^0(\g,V) = V^\g = V^G$ follows 
from the triviality of the $\g$-action on the closed $G$-invariant 
subspace $V^\g$ of $V$ (\cite{Ne06a}, Rem.~II.3.7). 
We may therefore assume $p > 0$. 
Let $\omega\in Z^p(\g,V)$ be a continuous $p$-cocycle and $g \in G$. We have 
to find an $\eta \in C^{p-1}(\g,V)$ with $g.\omega - g = \dd_\g \eta$. 
Let $\gamma \: [0,1] \to G$ be a smooth curve with 
$\gamma(0) = \1$ and $\gamma(1) = g$. We write $\delta^l(\gamma)_t := 
\gamma(t)^{-1}.\gamma'(t)$ for its left logarithmic derivative and 
recall the Cartan formula 
${\cal L}_x  = i_x \circ \dd_\g + \dd_\g \circ i_x$
for the action of $\g$ on $C^\bullet(\g,V)$. 
We thus obtain in the pointwise sense on each $p$-tuple 
of elements of $\g$: 
\begin{align*} 
g.\omega - \omega 
&= \int_0^1 \frac{d}{dt} \gamma(t).\omega \, dt 
= \int_0^1 \gamma(t).\big({\cal L}_{\delta^l(\gamma)_t}\omega\big) \, dt \\
&= \int_0^1 \gamma(t).\big(
i_{\delta^l(\gamma)_t}\dd_\g\omega + \dd_\g i_{\delta^l(\gamma)_t}\omega\big) \, dt 
= \int_0^1 \dd_\g\big(\gamma(t).i_{\delta^l(\gamma)_t}\omega\big) \, dt \\
&= \dd_\g \int_0^1\big(\gamma(t).i_{\delta^l(\gamma)_t}\omega\big) \, dt.  
\end{align*}
\smartqed\qed
\end{proof}

For $p = 2$ we obtain in particular $g.\omega - \omega  = \dd_\g \eta_g$ for 
\begin{equation}
  \label{eq:pot}
\eta_g(x) := \int_0^1 
\gamma(t).\omega\big(\delta^l(\gamma)_t, \Ad(\gamma(t))^{-1}.x\big)\, dt. 
\end{equation}



\end{document}